\documentclass[reqno]{amsart}
\usepackage{amsmath,amsfonts,amssymb,amsthm}
\usepackage[usenames,dvipsnames]{xcolor}
\usepackage{graphicx}
\usepackage{enumitem}
\usepackage[colorlinks=true]{hyperref} %para PDFLaTeX
\usepackage{url}
\usepackage[utf8]{inputenc}
\hypersetup{linkcolor=Violet,citecolor=PineGreen}
\newtheorem{teor}{Theorem}[section]

\newtheorem{prop}[teor]{Proposition}
\newtheorem{coro}[teor]{Corollary}
\theoremstyle{definition}
\newtheorem{defi}[teor]{Definition}

\newtheorem{hipo}[teor]{Hypothesis}

\newtheorem{nota}[teor]{Remark}
\newtheorem{notas}[teor]{Remarks}
\numberwithin{equation}{section}
\newcommand{\R}{\mathbb R}

\newcommand{\N}{\mathbb{N}}

\newcommand{\mB}{\mathcal{B}}

\newcommand{\mU}{\mathcal{U}}
\newcommand{\mM}{\mathcal{M}}

\newcommand{\mR}{\mathcal{R}}

\newcommand{\ep}{\varepsilon}

\newcommand{\mI}{\mathcal{I}}
\newcommand{\mJ}{\mathcal{J}}

\newcommand{\W}{\Omega}
\newcommand{\w}{\omega}
\newcommand{\lb}{\lambda}

\newcommand{\wma}{\wit{\mathfrak{a}}}
\newcommand{\wmr}{\wit{\mathfrak{r}}}
\newcommand{\mac}{\mathfrak{a}_c}
\newcommand{\mbc}{\mathfrak{b}_c}
\newcommand{\mdc}{\mathfrak{d}_c}
\newcommand{\mrc}{\mathfrak{r}_c}

\newcommand{\wit}{\widetilde}

\newcommand{\n}[1]{\left\|#1\right\|}

\newcommand{\lsm}{\left[\begin{smallmatrix}}
\newcommand{\rsm}{\end{smallmatrix}\right]}

\newcommand{\rmd}{\mathrm{d}}

\begin{document}
\title[Rate-induced tipping and saddle-node bifurcation]
{Rate-induced tipping and saddle-node bifurcation for quadratic differential
equations with nonautonomous asymptotic dynamics}
\author[I.P. Longo]{Iacopo P. Longo}
\author[C. N\'{u}\~{n}ez]{Carmen N\'{u}\~{n}ez}
\author[R. Obaya]{Rafael Obaya}
\author[M. Rasmussen]{Martin Rasmussen}
\address[I.P. Longo]{Technische Universit\"{a}t M\"{u}nchen,
Forschungseinheit Dynamics, Zentrum Mathematik, M8,
Boltzmannstra{\ss}e 3, 85748 Garching bei M\"{u}nchen, Germany.}
\address[C. N\'{u}\~{n}ez and R. Obaya]{Departamento de Matem\'{a}tica Aplicada,
Universidad de Valladolid, Paseo del Cauce 59, 47011 Valladolid, Spain.}
\address[M. Rasmussen]{Department of Mathematics, Imperial College London,
180 Queen's Gate, London SW7 2AZ, United Kingdom.}
\email[Iacopo Longo]{longoi@ma.tum.de}
\email[Carmen N\'{u}\~{n}ez]{carnun@wmatem.eis.uva.es}
\email[Rafael Obaya]{rafoba@wmatem.eis.uva.es}
\email[Martin Rasmussen]{m.rasmussen@imperial.ac.uk}
\thanks{All authors were partly supported by the European Union's Horizon 2020
research and innovation programme under the Marie Skłodowska-Curie
grant agreement No 643073. I.P. Longo was also partly supported by the
European Union's Horizon 2020 research and innovation programme under
the Marie Skłodowska-Curie grant agreement No 754462. C. N\'{u}\~{n}ez and R. Obaya
were also partly supported by Ministerio de Ciencia, Innovaci\'{o}n y
Universidades under project RTI2018-096523-B-I00.}
%37B25, %Lyapunov functions and stability; attractors, repellers
%34C23 %Bifurcation
%34D20 %Stability
%37B55, %Nonautonomous dynamical systems
%37G35, %Attractors and their bifurcations
%37M20, %Computational  methods for bifurcation problems
%37M22, %Computational methods for attractors
%34D05, %Asymptotic properties
%34D09, %Dichotomy, trichotomy
%34K14, %Almost and pseudo-periodic solutions for FDEs
%34K20, %Stability theory for FDEs
%92B20, %Neural networks, artificial life and related topics
\date{}
%\date{\today}
\begin{abstract}
An in-depth analysis of nonautonomous bifurcations of saddle-node
type for  scalar differential equations $x'=-x^2+q(t)\,x+p(t)$,
where $q\colon\R\to\R$ and $p\colon\R\to\R$ are bounded and uniformly
continuous, is fundamental to explain the absence or occurrence of
rate-induced tipping for the differential equation
$y' =(y-(2/\pi)\arctan(ct))^2+p(t)$ as the rate $c$ varies on $[0,\infty)$.
A classical attractor-repeller pair, whose existence for $c=0$ is assumed,
may persist for any $c>0$, or disappear for a certain critical rate $c=c_0$,
giving rise to rate-induced tipping. A suitable example demonstrates that one can have more than one critical rate, and the existence of the classical attractor-repeller pair may return when $c$ increases.
\end{abstract}
\maketitle
%%%%%%%%%%%%%%%%%%%%%%%%%%%%%%%%%%%%%%%%%%%%%%%%%%%%%%%%%%%%%%%%%%%%%%%%%%%%%%%%%%%%%
%%%%%%%%%%%%%%%%%%%%%%%%%%%%%%%%%%%%%%%%%%%%%%%%%%%%%%%%%%%%%%%%%%%%%%%%%%%%%%%%%%%%%
%%%%%%%%%%%%%%%%%%%%%%%%%%%%%%%%%%%%%%%%%%%%%%%%%%%%%%%%%%%%%%%%%%%%%%%%%%%%%%%%%%%%%
%%%%%%%%%%%%%%%%%%%%%%%%%%%%%%%%%%%%%%%%%%%%%%%%%%%%%%%%%%%%%%%%%%%%%%%%%%%%%%%%%%%%%
\section{Introduction}
Several complex systems in nature and society have been proved susceptible to abrupt, large and irreversible transitions in their behaviour, as a consequence
of relatively small changes in parameters describing external conditions. These
often unexpected changes are commonly referred to as tipping points
(or critical transitions), and they have been reported by applied scientists
in various contexts, including epileptic seizures, ecology, earthquakes,
and climate (see Scheffer \cite{scheffer}). Recent interdisciplinary research
efforts have stimulated the foundation of a mathematical theory for the occurrence
of tipping points. It is now understood (see Ashwin, Wieczorek, Vitolo and
Cox~\cite{ashwintippinclass}) that there are three different mechanisms of
tipping: \emph{bifurcation tipping}, which can be explained using classical
bifurcation theory; \emph{noise tipping}, which involves a transition from
one to another attractor due to noisy fluctuations; and \emph{rate-induced
tipping}, which involves a fast change in the parameters, so that tracking
of an attractor is no longer possible.
\par
Rate-induced tipping can be seen as a special type of a nonautonomous
bifurcation, which manifests itself on a finite time interval, on which
the parameters change significantly and non-adiabatically. 
This encompasses various real scenarios in ecology \cite{svNhh,vwf}, climate \cite{aaJQW,aspw,walc}, biology \cite{hill} and quantum mechanics \cite{kato} among others.
Recently, a
framework has been developed that allows the analysis of this type of
finite-time bifurcation using asymptotic theory (see Ashwin, Perryman
and Wieczorek~\cite{aspw}), by using an appropriate nonautonomous transition
between two different autonomous systems that represent the past and the
future of the model. While the nonautonomous transition between the past
and future systems is defined on an infinite time interval, its speed on
a finite time interval, given by a so-called \emph{rate}, is crucial for
the dynamics of the system, which may drastically change as the rate varies,
giving rise to rate-induced tipping. The infinite-time analysis involves a
local pullback attractor, which represents the behavior in the past, and
the tipping takes place when the future behavior of the pullback attractor
changes under variation of the rate. This change of the forward limit of
the pullback attractor can be linked to a collision of the pullback attractor
with an unstable object belonging to the future system, and such a collision
has been explored analytically and numerically in several contexts so far,
in one-dimensional systems (as in \cite{aspw}), higher-dimensional systems
(as in Alkhayoun and Ashwin~\cite{alas}, Wieczorek, Xie and Jones
\cite{wieczcompact}, and Xie \cite{xiethesis}), set-valued dynamical systems
(as in Carigi \cite{carigimres}), and random dynamical systems
(as in Hartl \cite{hartlthesis}).
\par
In this paper, we investigate rate-induced tipping in a situation
where not only the transition between the past and future systems is
nonautonomous, but also both past and future system are nonautonomous.
Since, in general, nonautonomous differential equations have no constant
or periodic solutions, the local and global dynamics can be very
complicated. We work under the fundamental assumption that the past
and the future systems are described by quadratic concave differential
equations whose global dynamics are governed by the presence of a
classical attractor-repeller pair. This type of Riccati equations have been frequently used as models in mathematics and many different areas of applied sciences. In particular, the scalar and matrix Riccati equations have been extensively studied in the literature due to their theoretical interest and their important applications in calculus of variations, optimal control or dissipative control problems, among many others (see Coppel~\cite{copp2} and Johnson, Obaya, Novo, N\'{u}\~{n}ez and Fabbri~\cite{jonnf}). With specific reference to rate-induced tipping, we point out that a time-dependent parameter drift in autonomous quadratic differential equations has been treated in \cite{aspw} and in \cite{ashwintippinclass} where, in particular, a ``zero-dimensional'' global energy balance climate model (the Friedrich's model) is re-conducted to a scalar autonomous concave differential equation. Notably, some of the involved parameters (accounting for the variation of the planetary orbit or the solar constant for example) admit a natural quasi-periodic variation. Therefore, the conclusions of this work show that these nonautonomous equations also  provide an interesting insight in the analysis and understanding of the critical transitions and their implications in models of real phenomena. 

Although our results can be adapted
to more general situations, in order to increase clarity, we focus on
differential equations of the form
\begin{equation}\label{1.ecu}
 y' =-\Big(y-\frac{2}{\pi}\,\arctan(ct)\Big)^2+p(t)\,,
\end{equation}
where $p:\R \to \R$ is a bounded and uniformly continuous
function, and $c\ge0$. Due to the asymptotic behavior of $\arctan(ct)$
for $c>0$, the differential equation \eqref{1.ecu} models a transition
from the past equation
\begin{equation}\label{1.ecu-}
 y' =-(y+1)^2+p(t)
\end{equation}
to the future equation
\begin{equation}\label{1.ecu+}
 y' =-(y-1)^2+p(t)\,,
\end{equation}
and the speed of this transition is given by the rate $c>0$:
small values of $c$ describe a slow transition, while for large values
of $c$, a large part of the transition from \eqref{1.ecu-} to
\eqref{1.ecu+} takes place rapidly on a small interval around $t=0$.
\par
Under the aforementioned assumption that the past and future equation
have a classical attractor-repeller pair, we show the following:
\begin{itemize}[leftmargin=*]
\item For small $c>0$, the differential equation
\eqref{1.ecu} also has a classical attractor-repeller pair, which
connects forward and backward in time to the attractor-repeller pairs
of \eqref{1.ecu-} and \eqref{1.ecu+}. This means that the unique
attractor of \eqref{1.ecu} converges to the attractors of the limiting
equations \eqref{1.ecu-} and \eqref{1.ecu+} in the limits
$t\to-\infty$ and $t\to\infty$, respectively, and the same holds for
the unique repeller of \eqref{1.ecu}.
\item For certain functions $p$, rate-induced tipping can not occur,
which means that the aforementioned connection of the
attractor-repeller pairs holds for every $c>0$.
\item There exist functions $p$ for which the connection between
the attractor-repeller pairs breaks up at a certain rate $c=c_0>0$,
giving rise to a rate-induced tipping point. In these cases,
the attractor-repeller pair for \eqref{1.ecu} collides in the sense
that the distance between the attractor and the repeller goes to
zero in the limit $c\to c_0^-$, and for $c=c_0$, a unique bounded
solution exists that is both a local pullback attractor and a local
pullback repeller, meaning that it is attractive in the past and
repulsive in the future. This follows from the fact that this
unique bounded solution connects backward in time to the attractor
of the past equation \eqref{1.ecu-} and forward in time to the
repeller of the future equation \eqref{1.ecu+}. So while for $c<c_0$,
there were two connections (from attractor to attractor, and from
repeller to repeller), at $c=c_0$, there is only one connection,
from the attractor in the past to the repeller in the future.
If $c$ is increased further beyond $c_0$, one may observe a complete
lack of connections, which comes from the fact that there do not exist
any  bounded solutions. However, there will always exist a solution
that is a pullback attractor and a solution that is a pullback
repeller, but these solutions will not be defined for all times.
\item The point $c_0$ may not be unique, since it is possible that
there exists another tipping rate $c_1>c_0$ such that \eqref{1.ecu}
has again an attractor-repeller pair for certain $c>c_1$.
In addition, there may be several intervals in the $c$-space that can be characterized by the existence or non-existence of the connection between the attractor-repeller pairs.
\end{itemize}
\par
Apart from the theoretical analysis, we have also conducted numerical
studies that completed our understanding of the tipping phenomena
described above. In particular, we show by means of an example
that rate-induced tipping can be reversed, as mentioned above.
\par
This paper is organized as follows. The short Section \ref{2.sec}
contains basic definitions for nonautonomous differential equations
and flows. Section~\ref{3.sec} is mainly devoted to the analysis of
the concave quadratic scalar equation $x'=-x^2+q(t)\,x+p(t)$, whose
dynamical possibilities play a fundamental role in our results on
rate-induced tipping. In particular, we show that the existence of
two hyperbolic solutions is equivalent to the existence of
two bounded and uniformly separated solutions, and we provide a
global description of the dynamics in this case.
In Section~\ref{4.sec}, we come back to the differential equation
\eqref{1.ecu}, in order to obtain the results described above.
We also present our numerical observations and illustrations
that complete our understanding of the tipping phenomena.
Due to their length, Sections~\ref{3.sec} and \ref{4.sec}
 have been  divided into several subsections.
\par
We close this introduction by pointing out that the results of
Section~\ref{3.sec}, which are crucial in the description of
the rate-induced tipping in Section~\ref{4.sec}, are of independent
interest. As a matter of fact, they are (far from trivial)
generalizations of results from the papers Alonso and Obaya \cite{alob3}
and N\'{u}\~{n}ez, Obaya and Sanz~\cite{nuos4,nuos5}
to the case of a non-recurrent concave quadratic equation,
which we consider in this article. As we explain at the end of
Section~\ref{3.sec}, they constitute by themselves an analysis of
a nonautonomous bifurcation pattern of saddle-node type,
in the line of the results of N\'{u}\~{n}ez and Obaya \cite{nuob6} and
Anagnostopoulou and J\"{a}ger \cite{anja}.

%%%%%%%%%%%%%%%%%%%%%%%%%%%%%%%%%%%%%%%%%%%%%%%%%%%%%%%%%%%%%%%%%%%%%%%%%%%%%%%%%%%%%
%%%%%%%%%%%%%%%%%%%%%%%%%%%%%%%%%%%%%%%%%%%%%%%%%%%%%%%%%%%%%%%%%%%%%%%%%%%%%%%%%%%%%
%%%%%%%%%%%%%%%%%%%%%%%%%%%%%%%%%%%%%%%%%%%%%%%%%%%%%%%%%%%%%%%%%%%%%%%%%%%%%%%%%%%%%
%%%%%%%%%%%%%%%%%%%%%%%%%%%%%%%%%%%%%%%%%%%%%%%%%%%%%%%%%%%%%%%%%%%%%%%%%%%%%%%%%%%%%
\section{Some notions on nonautonomous differential equations and flows}\label{2.sec}
Let $h\colon\R\times \R\to\R$ be a continuous function such that
$\partial h/\partial x\colon\R\times\R\to\R$ exists and is continuous.
We consider the nonautonomous scalar equation
\begin{equation}\label{2.ecuh}
 x'=h(t,x)\,.
\end{equation}
Let $t\mapsto x(t,s,x_0)$ denote the maximal solution of the initial
value problem $x(s)=x_0$ for \eqref{2.ecuh}. The so-defined real-valued
mapping~$x$ is defined on an open subset of $\R\times\R\times\R$ that
contains the set $\{(s,s,x_0)\,|\;(s,x_0)\in\R\times\R\}$, and we have
the two identities $x(s,s,x_0)=x_0$ and $x(t,l,x(l,s,x_0))=x(t,s,x_0)$,
whenever all the involved terms are defined.
\par
As mentioned in the Introduction, in Section 4, we deal with the
occurrence of hyperbolic solutions that lose hyperbolicity during
a critical transition and become only locally pullback attractive
or repulsive. In order to avoid interruption of the discussion
there, we explain the required notions of hyperbolicity, attractivity
and repulsivity now, and we refer the reader for in-depth analyses
of nonautonomous attractors and repellers to \cite{klra}, \cite{rasmln} and
\cite{calr}.
\par
A globally defined solution $\wit b\colon\R\to\R$ of \eqref{2.ecuh}
is said to be {\em hyperbolic\/} if the corresponding variational
equation $z'= (\partial /\partial x)h(t, \wit b(t))\,z$ has an
exponential dichotomy on $\R$ \cite{copp1}. In this one-dimensional
context, the existence of an exponential dichotomy means that there
exist $k_b\ge 1$ and $\beta_b>0$ such that either
\begin{equation}\label{3.masi}
 \exp\int_s^t (\partial /\partial x)h(l, \wit b(l))\,\rmd l\le k_b\,e^{-\beta_b(t-s)} \quad
 \text{whenever $t\ge s$}
\end{equation}
or
\begin{equation}\label{3.menosi}
 \exp\int_s^t (\partial /\partial x)h(l, \wit b(l))\,\rmd l\le k_b\,e^{\beta_b(t-s)} \quad
 \text{whenever $t\le s$}
\end{equation}
holds. In the first case \eqref{3.masi}, the variational equation
is called {\em Hurwitz at $+\infty$}, and the hyperbolic solution
$\wit b$ is said to be {\em (locally) attractive}. In the second case
\eqref{3.menosi},  the variational equation is called {\em Hurwitz at $-\infty$},
and the hyperbolic solution $\wit b$ is said to be {\em (locally) repulsive}.
In both cases, we call $(k_b,\beta_b)$ a (non-unique)
{\em dichotomy constant pair\/} for the hyperbolic solution $\wit b$ (or for the variational
equation $z'=(\partial /\partial x)h(t, \wit b(t))\,z$).
\par
With the aim to clarify the notation as much as possible,
we write $\wit b$, $\wit a$, $\wit r$, etc., whenever we know that
these functions are hyperbolic solutions.
\par
Note that if the hyperbolic solution $\wit b$ is attractive, then
all  (non-trivial) solutions of the variational equation tend to $0$ as
$t\to\infty$, and they converge to $+\infty$ or $-\infty$ as $t\to-\infty$.
Accordingly, if the hyperbolic solution $\wit b$ is repulsive, then
all (non-trivial) solutions of the variational equation tend to $0$ as
$t\to-\infty$, and they converge to $+\infty$ or $-\infty$ as $t\to\infty$.
A proof of this well-known fact can be found, for instance,
in \cite[Proposition~1.56]{jonnf}.
\par
An attractive hyperbolic solution attracts nearby solutions forward in time,
and a repulsive hyperbolic solution attracts nearby solutions backward in time.
Since \eqref{3.masi} and \eqref{3.menosi} hold on the entire line $\R$,
this form of attraction and repulsion takes place at all times. During the
process of rate-induced tipping, attraction and repulsion is lost on a half line,
leading to attractive solutions becoming only attractive in the past, and repulsive
solutions becoming only repulsive in the future. The following notions of local
pullback attractivity and repulsivity, adapted from \cite[Section~2.3]{rasmln},
describe this behavior.
\par
A solution $\bar a\colon(-\infty,\beta)\to\R$ (with $\beta\le\infty$) of
\eqref{2.ecuh} is called {\em locally pullback attracting\/} if there exist
$s_0<\beta$ and $\delta>0$ such that if $s\le s_0$ and $|x_0-\bar a(s)|<\delta$
then $x(t,s,x_0)$ is defined for $t\in[s,s_0]$, and in addition
\[
 \lim_{s\to-\infty}\max_{x_0\in [\bar a(s)-\delta,\bar a(s)+\delta]}|\bar a(t)-x(t,s,x_0)|=0
 \quad \text{for all $t\le s_0$}\,.
\]
Note that, in our scalar case, this is equivalent to saying that if
$s\le s_0$, then the solutions $x(t,s,a(s)\pm\delta)$ are defined for
$t\in[s,s_0]$ and, in addition,
\[
 \lim_{s\to-\infty}|\bar a(t)-x(t,s,\bar a(s)\pm\delta)|=0
 \quad \text{for all $t\le s_0$}\,.
\]
\par
A solution $\bar r\colon(\alpha,\infty)\to\R$ (with $\alpha\ge-\infty$)
of \eqref{2.ecuh} is called {\em locally pullback repulsive\/} if the
corresponding solution $\bar r^*\colon(-\infty,-\alpha)\to\R$ of the
differential equation under time reversal $y'=-h(-t,y)$, given by
$\bar r^*(t)=\bar r(-t)$, is locally pullback attractive.
In other words, if there exist $s_0>\alpha$ and $\delta>0$ such that
if $s\ge s_0$, then the solutions $x(t,s,\bar r(s)\pm\delta)$ are defined for
$t\in[s_0,s]$ and, in addition,
\[
 \lim_{s\to\infty}|\bar r(t)-x(t,s,\bar r(s)\pm\delta)|=0
 \quad \text{for all $t\ge s_0$}\,.
\]
\par
We proceed by summarizing some basic concepts and properties of
topological dynamics, which are needed in the proof of one our
main results, Theorem~\ref{3.teorhyp}.
\par
A (real and continuous) {\em global flow\/} on a complete metric
space $\W$ is a continuous map $\sigma\colon\R\times\W\to\W,\;
(t,\w)\mapsto\sigma(t,\w)=:\sigma_t(\w)$, such that
$\sigma_0=\text{Id}$ and $\sigma_{s+t}=\sigma_t\circ\sigma_s$
for each $s,t\in\R$.
The flow is {\em local\/} if the map $\sigma$
is defined, continuous, and satisfies the previous properties on an
open subset of $\R\times\W$ containing $\{0\}\times\W$.
\par
Let $(\W,\sigma)$ be a global flow.
The $\sigma$-{\em orbit\/} of a point $\w\in\W$
is the set $\{\sigma_t(\w)\,|\;t\in\R\}$.
%Restricting the time to $t\ge 0$ or
%$t\le 0$ provides the definition of {\em forward\/}
%or {\em backward\/} $\sigma$-semiorbit.
A subset $\W_1\subset\W$ is {\em $\sigma$-invariant\/}
%(resp.~{\em positively $\sigma$-invariant\/}
%and {\em negatively $\sigma$-invariant\/})
if $\sigma_t(\W_1)=\W_1$ for every $t\in\R$.
%(resp.~$t\ge 0$ and $t\le 0$).
A $\sigma$-invariant subset $\W_1\subset\W$ is {\em minimal\/} if it is compact
and does not contain properly any other compact $\sigma$-invariant set;
and the flow $(\W,\sigma)$ is
{\em minimal\/} if $\W$ itself is minimal.
If the set $\{\sigma_t(\w)\,|\;t\ge 0\}$ is relatively compact, then
the {\em omega-limit set\/} of $\w_0$ is given
by those points $\w\in\W$ such that $\w=\lim_{m\to \infty}\sigma(t_m,\w_0)$
for some sequence $(t_m)\uparrow \infty$. This set is
nonempty, compact, connected and $\sigma$-invariant.
The definition and properties of the {\em alpha-limit set\/} of $\w_0$
are analogous, working now with sequences $(t_m)\downarrow-\infty$.
\par
We end this short section by introducing the notation
$\n{b}:=\sup_{t\in\R}|b(t)|$ for the supremum norm
of any bounded continuous function $b\colon\R\to\R$.
%%%%%%%%%%%%%%%%%%%%%%%%%%%%%%%%%%%%%%%%%%%%%%%%%%%%%%%%%%%%%%%%%%%%%%%%%%%%%%%%%%%%%
%%%%%%%%%%%%%%%%%%%%%%%%%%%%%%%%%%%%%%%%%%%%%%%%%%%%%%%%%%%%%%%%%%%%%%%%%%%%%%%%%%%%%
%%%%%%%%%%%%%%%%%%%%%%%%%%%%%%%%%%%%%%%%%%%%%%%%%%%%%%%%%%%%%%%%%%%%%%%%%%%%%%%%%%%%%
%%%%%%%%%%%%%%%%%%%%%%%%%%%%%%%%%%%%%%%%%%%%%%%%%%%%%%%%%%%%%%%%%%%%%%%%%%%%%%%%%%%%%
\section{The dynamics of the concave scalar equation $x'=-x^2+q(t)\,x+p(t)$}
\label{3.sec}
As explained in the Introduction, our approach to rate-induced tipping is based on
an in-depth analysis of the dynamics of the nonautonomous
concave scalar differential equation $x'=-x^2+p(t)$,
and we establish several fundamental facts about this differential equation
in this section. Since the effort is the same, we work with the more general
quadratic equation $x'=-x^2+q(t)\,x+p(t)$.
\par
In Subsection~\ref{3.sec31},
we establish the existence of two fundamental \lq\lq special\rq\rq~solutions
$a$ and $r$ for a class of scalar differential equations
$x'=h(t,x)$ that includes both \eqref{1.ecu} and $x'=-x^2+q(t)\,x+p(t)$.
The graphs of these solutions
determine the areas of initial conditions giving rise to solutions that
are bounded backward or forward in time, respectively.
\par
In Subsection~\ref{3.sec32}, we prove that $a$ and $r$ are globally defined
and hyperbolic solutions
of $x'=-x^2+q(t)\,x+p(t)$ if and only if they are uniformly separated.
We note that this is the unique possibility for the
occurrence of (exactly two) hyperbolic solutions. We also
explain in detail the dynamics in this situation.
\par
Finally, in Subsection~\ref{3.sec33},
we introduce a real parameter $\lambda$ and describe the
only three dynamical scenarios which are possible for the differential equation
$x'=-x^2+q(t)\,x+p(t)+\lb$.
These three possibilities depend on the relation between $\lb$ and a special value
$\lb^*(q,p)$
of the parameter, and correspond to a nonautonomous bifurcation pattern
of saddle-node~type.
%%%%%%%%%%%%%%%%%%%%%%%%%%%%%%%%%%%%%%%%%%%%%%%%%%%%%%%%%%%%%%%%%%%%%%%%%%%%%%%%%
%%%%%%%%%%%%%%%%%%%%%%%%%%%%%%%%%%%%%%%%%%%%%%%%%%%%%%%%%%%%%%%%%%%%%%%%%%%%%%%%%
\subsection{Some facts on the coercitive scalar equation $x'=h(t,x)$}
\label{3.sec31}
Let the function $h\colon\R^2\to\R$ be continuous, locally Lipschitz in its
second argument,  such that $\limsup_{x\to\pm\infty} h(t,x)/x^2<0$ uniformly in
$t\in\R$. We consider the nonautonomous  scalar equation
\begin{equation}\label{3.ecuh}
x'=h(t,x)\,.
\end{equation}
As in Section \ref{2.sec}, we represent by $t\mapsto x(t,s,x_0)$
the maximal solution satisfying $x(s,s,x_0)=x_0$, defined on the
 interval $\mI_{s,x_0}=(\alpha_{s,x_0},\beta_{s,x_0})$. Note that
$-\infty\le \alpha_{s,x_0}<s<\beta_{s,x_0}\le\infty$. We define
\[
\begin{split}
 \mB^-&:=\Big\{(s,x_0)\in\R^2\,\Big|\;\sup_{t\in(\alpha_{s,x_0},s]}x(t,s,x_0)
 <\infty\Big\}\,,\\
 \mB^+&:=\Big\{(s,x_0)\in\R^2\,\Big|\;\inf_{t\in[s,\beta_{s,x_0})}x(t,s,x_0)
 >-\infty\Big\}\,,\\
\end{split}
\]
which may be empty sets. Due to the above assumptions on the function $h$,
there exist $\ep>0$ and $m>0$ such that
\begin{equation}\label{3.defm}
 h(t,\pm x)\le-\ep\quad\text{for all $t\in\R$ and $x\ge m$}\,,
\end{equation}
and this implies that for all $(s,x_0)\in\R^2$,
\begin{equation}\label{3.inf}
\liminf_{t\to(\alpha_{s,x_0})^+}x(t,s,x_0)>-m
\quad\text{and}\quad
\limsup_{t\to(\beta_{s,x_0})^-}x(t,s,x_0)<m\,.
\end{equation}
Therefore, $\alpha_{s,x_0}=-\infty$ for all $(s,x_0)\in\mB^-$ and
$\beta_{s,x_0}=\infty$ for all $(s,x_0)\in\mB^+$.
It is also clear that the sets $\mB^-$ and $\mB^+$ are invariant in the
sense that $(t,x(t,s,x_0))\in \mB^-$ for all $(s,x_0)\in\mB^-$ and
$t\in \mI_{s,x_0}$, and $(t,x(t,s,x_0))\in \mB^+$ for all
$(s,x_0)\in\mB^+$ and $t\in \mI_{s,x_0}$. The (possibly empty) set
\begin{equation}\label{3.defB}
 \mB:=\mB^-\cap\mB^+
\end{equation}
is the  set of initial pairs $(s,x_0)$ giving rise to (globally defined)
bounded solutions of \eqref{3.ecuh}. In addition,
\begin{align}
 x(t,s,m)<m\quad&\text{for all $t\in(s,\beta_{s,m})$}\,,\label{3.cotasup}\\
 \lim_{t\to(\alpha_{s,x_0})^+}x(t,s,x_0)=
 \infty\quad&\text{for all $s\in\R$ and $x_0\ge m$}\,,
 \label{3.limsup}\\
 x(t,s,-m)>-m\quad&\text{for all $t\in(\alpha_{s,m},s)$}\,,\label{3.cotainf}\\
 \lim_{t\to(\beta_{s,x_0})^-}x(t,s,x_0)=
 -\infty\quad&\text{for all $s\in\R$ and $x_0\le-m$}\,.
 \label{3.liminf}
\end{align}
Let us now define
\begin{equation}\label{3.defsets}
\begin{split}
 \mR^-&:=\{s\in\R\,|\;\text{there exists $x_0$ with $(s,x_0)\in\mB^-$}\}\,,\\
 \mB^-_s&:=\{x_0\in\R\,|\;(s,x_0)\in\mB^-\} \quad \text{for $s\in\mR^-$}\,,\\
 \mR^+&:=\{s\in\R\,|\;\text{there exists $x_0$ with $(s,x_0)\in\mB^+$}\}\,,\\
 \mB^+_s&:=\{x_0\in\R\,|\;(s,x_0)\in\mB^+\} \quad \text{for\; $s\in\mR^+$}\,.
\end{split}
\end{equation}

The following theorem shows the existence of a solution
$a\colon\mR^-\to(-\infty,m)$ of \eqref{3.ecuh} which is the maximal one in
$\mB^-$ (if $\mB^-\not=\emptyset$), and a solution
$r\colon\mR^+\to(-m,\infty)$ which is the minimal one in $\mB^+$
(if $\mB^+\not=\emptyset$). We see later in Theorem~\ref{3.teorhyp}
that, in an appropriate setting, these two solutions are globally
defined and form what we have called a
{\em classical attractor-repeller pair} in the Introduction.

\begin{teor}\label{3.teoruno}
Consider the differential equation \eqref{3.ecuh}, let $m>0$ satisfy
\eqref{3.defm}, and let $\mB^\pm,\mB,\mR^\pm$ and $\mB^\pm_s$ be
the sets defined above.
\begin{itemize}
\item[(i)] If $\mB^-$ is nonempty, then the set $\mR^-$ is either $\R$ or
a negative open half-line; for each
$s\in\mR^-$, we have $\mB^-_s=(-\infty,a(s)]$, where the
map $a\colon\mR^-\to(-\infty,m)$ is a solution of \eqref{3.ecuh}. In addition,
if $s\in\mR^-$ then
$x(t,s,x_0)$ is bounded for $t\to-\infty$ if and only if $x_0\le a(s)$; and if
$\sup\mR^-<\infty$, then $\lim_{t\to(\sup\mR^-)^-} a(t)=-\infty$.
\item[(ii)] If $\mB^+$ is nonempty, then
the set $\mR^+$ is either $\R$ or a positive open half-line; for each
$s\in\mR^+$, we have $\mB^-_s=[r(s),\infty)$, where the
map $r\colon\mR^+\to(-m,\infty)$  is a solution of \eqref{3.ecuh}. In addition,
if $s\in\mR^+$ then $x(t,s,x_0)$ is bounded for $t\to\infty$
if and only if $x_0\ge r(s)$; and
if $\inf\mR^+>-\infty$, then $\lim_{t\to(\inf\mR^+)^+} r(t)=\infty$.
\item[(iii)] Let $x$ be a solution of \eqref{3.ecuh} defined on a
maximal interval $(\alpha,\beta)$.
If it satisfies $\liminf_{t\to\beta^-}x(t)=-\infty$, then $\beta<\infty$; and if
$\limsup_{t\to\alpha^+}x(t)=\infty$, then $\alpha>-\infty$.
In particular, any globally defined solution is bounded.
\item[(iv)] $\mB$ is nonempty if and only if $\mR^-=\R$ or $\mR^+=\R$, in
which case both equalities hold, $a$ and $r$ are
globally defined and bounded solutions of \eqref{3.ecuh}, and
\[
 \mB=\{(s,x_0)\in\R^2\,|\; r(s)\le x_0\le a(s)\}\,.
\]
\item[(v)] If there exists a bounded $C^1$ function $b\colon\R\to\R$
such that $b'(t)\le h(t,b(t))$ for all $t\in\R$, then
$\mB$ is nonempty, and $r(t)\le b(t)\le a(t)$ for all $t\in\R$.
And if $b'(t)< h(t,b(t))$ for all $t\in\R$, then
$r(t)<b(t)< a(t)$ for all $t\in\R$.
\end{itemize}
\end{teor}
\begin{proof}
(i) Take $(s,x_0)\in\mB^-$ and take $l>0$. Then, $(s-l,x(s-l,s,x_0))\in\mB^-$,
so that $\mR^-$ is either $\R$ or a negative half-line; and
$(s+l,x(s+l,s,x_0))\in\mB^-$ if $l$ is small enough, so that
$\mR^-$ is open. In addition, if $(s,x_0)\in\mB^-$
and $y_0<x_0$, then $x(t,s,y_0)<x(t,s,x_0)$ whenever both terms are defined,
from which it is easy to deduce that $(s,y_0)\in\mB^-$. And
\eqref{3.limsup} ensures that $(s,m)\notin\mB^-$ for any $s\in\R$. Thus,
if $s\in\mR^-$, then
$\mB^-_s$ is a negative half-line bounded from above by $m$.
We define $a(s):=\sup\mB_s^-$ for $s\in\mR^-$.
\par
Let us now prove that $a(s)$ belongs to the set $\mB^-_s$.
We take an increasing
sequence $(a_n)$ in $\mB_s^-$ with $\lim_{n\to\infty}a_n=a(s)$.
Since $(t,x(t,s,a_n))\in\mB^-$ for any $t\le s$ and
$n\in\N$, the function $t\mapsto x(t,s,a_n)$ is
defined at least in $(-\infty,s]$, where it satisfies
$x(t,s,a_n)< m$. Hence
$x(t,s,a(s))=\lim_{n\to\infty} x(t,s,a_n)\le m$
as long as the left hand term is defined.
Combined with the first inequality in
\eqref{3.inf}, we conclude that $t\mapsto x(t,s,a(s))$ is defined
and bounded in $(-\infty,s]$, which shows that $(s,a(s))\in\mB^-$,
as asserted. Note also that $a(s)<m$.
\par
Now we take $\bar s\in\mR^-$ and $\bar t\le \bar s$,
so that $\bar t\in\mR^-$. Since $(\bar s,a(\bar s))\in\mB^-$,
there exists $x(\bar t,\bar s,a(\bar s))$, and
$x(\bar t,\bar s,a(\bar s))\le a(\bar t)$. Note also that
\begin{equation}\label{3.jaleo}
 x(l,\bar t,x(\bar t,\bar s,a(\bar s)))=
 x(l,\bar s,a(\bar s)) \quad \mbox{for all }  l\in[\bar t,\bar s]\,.
\end{equation}
On the other hand,
the solution $l\mapsto x(l,\bar t,a(\bar t))$ is defined at least for
$l$ in an interval $[\bar t,\bar l]$ with $\bar t<\bar l\le\bar s$,
and it satisfies
$x(l,\bar t,a(\bar t))\ge x(l,\bar t,x(\bar t,\bar s,a(\bar s)))=
x(l,\bar s,a(\bar s))$ for $l\in[\bar t,\bar l]$,
so that $x(l,\bar t,a(\bar t))$ is bounded from below. The inequality
\eqref{3.cotasup} shows that it is also bounded from above. Therefore,
we can take $\bar l=\bar s$, which means that
$x(l,\bar t,a(\bar t))$ exists and is bounded
for (at least) $l\in[\bar t,\bar s]$,
and hence it exists and is bounded for $l\in(-\infty,\bar s]$.
In addition, $(\bar s,x(\bar s,\bar t,a(\bar t)))\in\mB^-$ (and hence
$x(\bar s,\bar t,a(\bar t))\le a(\bar s)$): as long as
$x(l,\bar s,x(\bar s,\bar t,a(\bar t)))$ is defined, it coincides
with $x(l,\bar t,a(\bar t))$; and this last solution exists
on $(-\infty,\bar s]$ and is bounded, as just seen.
The already obtained inequalities
$x(\bar t,\bar s,a(\bar s))\le a(\bar t)$ and
$x(\bar s,\bar t,a(\bar t))\le a(\bar s)$, and
the equality \eqref{3.jaleo} for $l=\bar s$, yield
\[
 a(\bar s)=x(\bar s,\bar s,a(\bar s))=x(\bar s,\bar t,x(\bar t,\bar s,a(\bar s)))
 \le x(\bar s,\bar t,a(\bar t))\le a(\bar s)\,,
\]
from where we deduce, first, that $x(\bar s,\bar t,a(\bar t))=a(\bar s)$
and, second, that $x(\bar t,\bar s,a(\bar s))=a(\bar t)$.
The conclusion is that $x(t,s,a(s))= a(t)$ whenever $s,t\in\mR^-$, as asserted.
\par
The last two assertions of (i) follow easily from the definition
of $a$ and from the fact that it is a solution of the
equation.
\smallskip\par
(ii) The proof is analogous to that of (i), making now use of
\eqref{3.liminf}, the second inequality in \eqref{3.inf},
and \eqref{3.cotainf}.
\smallskip\par
(iii) Note that $\lim_{t\to\beta^-}x(t)=-\infty$ if
$\limsup_{t\to\beta^-}x(t)=-\infty$. Let us assume that this is the case and,
for contradiction, that $\beta=\infty$. The conditions initially assumed
on $h$ ensure the existence of $\mu>0$ and $s_0>\alpha$ with
$x(s_0)\ne 0$ and such that
$h(t,x(t))/x^2(t)\le-\mu<0$ for all $t\ge s_0$. Then,
$x'(t)/x^2(t)\le -\mu$, so that
\[
 -\frac{1}{x(t)}+\frac{1}{x(s_0)}\le-\mu\,(t-s_0)
\]
for $t\ge s_0$,
and we get the contradiction by taking limit as $t\to\infty$.
The second assertion in (iii) is proved in the analogous way.
And the last one follows from these properties and \eqref{3.inf}.
\smallskip\par
(iv) The definition \eqref{3.defB} and the previous properties
prove the assertions in (iv).\hspace{-1cm}~
\smallskip\par
(v) Let us assume that $b'(s)\le h(s,b(s))$ for any $s\in\R$,
and fix any $s_0\in\R$. Then,
standard comparison results for first order scalar differential equations ensure that
$b(t)\le x(t,s_0,b(s_0))$ if $t>s_0$ and $x(t,s_0,b(s_0))$ exists, and
$b(t)\ge x(t,s_0,b(s_0))$ if $t<s_0$ and $x(t,s_0,b(s_0))$ exists.
It follows from the first inequality (and the global
existence of $b(t)$) that $(s_0,b(s_0))$ belongs to $\mB^+$, and
from the second one that $(s_0,b(s_0))$ belongs to $\mB^-$,
so that $\mB$ is not empty and the maps $a$ and $r$ are globally defined.
In addition, $a(t)=a(t,t-1,a(t-1))\ge x(t,t-1,b(t-1))\ge b(t)$ and
$r(t)=x(t,t+1,r(t+1))\le x(t,t+1,b(t+1))\le b(t)$ for any $t\in\R$.
This proves (v) in the first case.
Finally, if the initial inequality
is strict, so are those derived from the comparison results,
and we conclude that $r(t)<b(t)<a(t)$ for any $t\in\R$.
This completes the proof of (v).
\end{proof}
%%%%%%%%%%%%%%%%%%%%%%%%%%%%%%%%%%%%%%%%%%%%%%%%%%%%%%%%%%%%%%%%%%%%%%%%%%%%%%%%%
%%%%%%%%%%%%%%%%%%%%%%%%%%%%%%%%%%%%%%%%%%%%%%%%%%%%%%%%%%%%%%%%%%%%%%%%%%%%%%%%%
\subsection{Hyperbolic solutions for
%$x'=-x^2+p(t)$.}\label{3.sec32}
$x'=-x^2+q(t)\,x+p(t)$.}\label{3.sec32}
Let $q\colon\R\to\R$ and $p\colon\R\to\R$ be bounded and uniformly
continuous functions.
In Section \ref{4.sec}, we will apply Theorem~\ref{3.teoruno}
mainly (but not only) to scalar concave equations of the type
\begin{equation}\label{3.ecucon}
 x'=-x^2+q(t)\,x+p(t)\,.
\end{equation}
In fact, we will work there mainly with $x'=-x^2+p(t)$.
As mentioned in the Introduction, we need to extend
part of the properties of recurrent concave equations proved in
\cite{alob3}, \cite{nuos4} and \cite{nuos5}
(see Remark \ref{A.notarec} below).
The proofs of the main results in this section
rely deeply on those of \cite{nuos4}.
\par
We first establish robustness of hyperbolicity for the differential
equation~\eqref{3.ecucon}. A proof of this well-known property in a
more general setting can be found in \cite[Theorem~3.8]{potz},
but we include a direct proof for the reader's convenience.
\begin{prop}\label{3.proppersiste}
Assume that \eqref{3.ecucon} has an attractive (resp.~repulsive)
hyperbolic solution $\wit b_{q,p}$. Then this hyperbolic solution
is persistent in the following sense. For $\ep>0$, there exists
$\delta_\ep>0$ such that, if $\bar q\colon\R\to\R$ and
$\bar p\colon\R\to\R$ are bounded and continuous functions
with $\n{\bar q-q}<\delta_\ep$ and $\n{\bar p-p}<\delta_\ep$,
then also the perturbed differential equation
\[
 x'=-x^2+\bar q(t)\,x+\bar p(t)
\]
has an attractive  (resp.~repulsive) hyperbolic solution
$\wit b_{\bar q,\bar p}$ that satisfies
$\|\wit b_{q,p}-\wit b_{\bar q,\bar p}\|<\ep$. In addition,
there exists a common
dichotomy constant pair for all the variational equations
$z'=(-2\,\wit b_{\bar q,\bar p}(t)+\bar q(t))\,z$,
where, as above, the functions $\bar q$ and $\bar p$ satisfy
$\n{\bar q-q}<\delta_\ep$ and
$\n{\bar p-p}<\delta_\ep$.
\end{prop}
\begin{proof}
The proof is based on that of Lemma 3.3 in \cite{aloo}.
Let us assume that the variational equation
$z'=(-2\,\wit b_{q,p}(t)+q(t))\,z$ is Hurwitz
at $+\infty$ with dichotomy
constant pair $(k_b,\beta_b)$, where $k_b\ge 1$ and $\beta_b>0$,
see Section~\ref{2.sec}.
It is easy to deduce from this definition that there
exists $\bar\delta>0$ such that if
$\bar q,\,y_0\in C(\R,\R)$ satisfy $\n{y_0}<\bar\delta$ and
$\n{\bar q-q}<\bar\delta$,
then $z'=(-2\,\wit b_{q,p}(t)+\bar q(t)+y_0(t))\,z$ is also
Hurwitz at $+\infty$, and that there
exists a common dichotomy constant pair that is valid for all the perturbed
equations corresponding to $\bar q$ and $y_0$ with $\n{y_0}<\bar\delta$ and
$\n{\bar q-q}<\bar\delta$.
We assume without restriction that $\bar\delta\le
\min(\beta_b/(3k_b),1)$.
\par
The change of variables $y=x-\wit b_{q,p}$ takes
$x'=-x^2+\bar q(t)\,x+\bar p(t)$ to
\begin{equation}\label{3.ecuz}
 y'=(-2\,\wit b_{q,p}(t)+q(t))\,y-y^2+(\bar q(t)-q(t))\,y+s(t)\,,
\end{equation}
where $s(t):=(\bar q(t)-q(t))\,\wit b_{q,p}(t)+\bar p(t)-p(t)$.
The results of \cite[Lecture 3]{copp1} (see also \cite[Theorem 7.7]{fink})
ensure that for any $y_0\in C(\R,\R)$ with $\n{y_0}\le\bar\delta\le 1$,
there exists a unique bounded solution $Ty_0$ of
\[
 y'=(-2\,\wit b_{q,p}(t)+q(t))\,y-y_0^2(t)+(\bar q(t)-q(t))\,y_0(t)+s(t)\,,
\]
given by
$Ty_0(t):=\int_{-\infty}^t u(t)\,u^{-1}(l)\,
(-y_0^2(l)+(\bar q(l)-q(l))\,y_0(l)+s(l))\,\rmd l$
for $u(t):=\exp\int_{0}^t(-2\wit b_{q,p}(l)+q(l))\,\rmd l$.
Therefore,
$\n{Ty_0}\le(k_b/\beta_b)(\n{y_0}^2+\n{\bar q-q}+\n{s})$
and $\n{Ty_1-Ty_2}\le (k_b/\beta_b)\n{y_1+y_2+q-\bar q}\n{y_1-y_2}$.
Recall that $0<\bar\delta<\beta_b/(3k_b)$. It is easy
to check that if $\n{\bar q-q}+\n{s}\le(\beta_b\bar\delta/(2k_b))
(1-2k_b\bar\delta/\beta_b)$, if $\n{\bar q-q}\le\bar\delta$,
and if $\n{y_i}\le\bar\delta$ for $i\in\{0,1,2\}$, then
\[
 \n{Ty_0}\le\frac{\bar\delta}{2}
 \quad\text{and}\quad
 \n{Ty_1-Ty_2}\le\frac{3\,k_b\,\bar\delta}{\beta_b}\,\n{y_1-y_2}\,.
\]
Since $3\,k_b\,\bar\delta/\beta_b<1$, the map
$T\colon C(\R,[-\bar\delta,\bar\delta])\to
C(\R,[-\bar\delta,\bar\delta])$ is a contraction and thus has a fixed
point $\bar y_{\bar q,\bar p}$. Clearly,
$\bar y_{\bar q,\bar p}$  solves \eqref{3.ecuz}, so that
$\wit b_{\bar q,\bar p}:=\bar y_{\bar q,\bar p}+\wit b_{q,p}$ solves
$x'=-x^2+\bar q(t)\,x+\bar p(t)$.
\par
Let us take $\ep>0$ and assume from the beginning,
without restriction, that $\bar\delta<2\ep$ (so that $\bar\delta$ may depend
on $\ep$). Now we choose a positive
$\delta_\ep<\min(\beta_b\bar\delta/(12k_b),\bar\delta)$, as
well as continuous functions $\bar q$ and $\bar p$
such that $\n{\bar q-q}\le\delta_\ep$
(so that $\n{\bar q-q}\le\bar\delta$, as previously required)
and $\n{s}<\delta_\ep$. Then,
$\n{\bar q-q}+\n{s}\le2\delta_\ep\le
(\beta_b\bar\delta/(2k_b))(1-2k_b\bar\delta/\beta_b)$,
and hence, the fixed point $\bar y_{\bar q,\bar p}$ exists and satisfies
$\n{\bar y_{\bar q,\bar p}}<\bar\delta/2$. Therefore, since
$\n{-2\bar y_{\bar q,\bar p}}\le\bar\delta$ and $\n{\bar q-q}\le\bar\delta$,
we can assert that the differential equation
$z'=(-2\,\wit b_{\bar q,\bar p}(t)+\bar q(t))\,z=
(-2\,\wit b_{q,p}(t)-2\bar y_{\bar q,\bar p}+\bar q(t))\,z$ is Hurwitz
at $+\infty$; and $\|\wit b_{q,p}-\wit b_{\bar q,\bar p}\|
=\n{\bar y_q}\le\bar\delta/2<\ep$.
This completes the proof in this first case, and the other case,
when the variational equation is Hurwitz at $-\infty$,
can be proved in analogously.
\end{proof}
The following basic properties for equation \eqref{3.ecucon}
play a role in the rest of the section. From the solution identity
$(\partial/\partial t)\,x(t,s,x_0)=-x^2(t,s,x_0)+q(t)\,x(t,s,x_0)+p(t)$, we get
$(\partial/\partial x_0)\,x(t,s,x_0)=\exp\int_s^t(-2\,x(l,s,x_0)+q(l))\,\rmd l$.
Therefore, the map $x_0\mapsto x(t,s,x_0)$ is strictly concave if $t>s$
(since its derivative with respect to $x_0$ decreases strictly in $x_0$, as does
$-2\,x(l,s,x_0)$) and convex if $t<s$ (since its derivative with respect to $x_0$
increases strictly with $x_0$). That is,
\begin{equation}\label{3.concx}
 x(t,s,\rho\,x_1+(1-\rho)\,x_2)>\rho\,x(t,s,x_1)+(1-\rho)\,x(t,s,x_2)
\end{equation}
whenever $t>s$, $x_1,x_2\in\R$ and $\rho\in(0,1)$
as long as all the involved terms are defined. Note that the sign of the
inequality changes for $t<s$.
Let us now take $x_1\le x_2$. From
\[
 x(t,s,x_2)-x(t,s,x_1)=\int_0^1 \frac{\partial}{\partial x_0}\,
 x(t,s,\lb\,x_2+(1-\lb)\,x_1)(x_2-x_1)\,\rmd\lb\,,
\]
and from $x(l,s,x_1)\le x(l,s,\lb\,x_2+(1-\lb)\,x_1)\le x(l,s,x_2)$
for $\lb\in[0,1]$, we deduce that
\begin{align}
\exp\int_s^t&(-2\,x(l,s,x_2)+q(l))\,\rmd l
\le\frac{x(t,s,x_2)-x(t,s,x_1)}{x_2-x_1}\nonumber\\
&\le\exp\int_s^t(-2\,x(l,s,x_1)+q(l))\,\rmd l\qquad
\text{for $x_1<x_2$ and $t\ge s$}\,,\label{3.tmays}
\end{align}
and
\begin{align}
\exp\int_s^t&(-2\,x(l,s,x_1)+q(l))\,\rmd l\le
\frac{x(t,s,x_2)-x(t,s,x_1)}{x_2-x_1}\nonumber\\
&\le\exp\int_s^t(-2\,x(l,s,x_2)+q(l))\,\rmd l\qquad
\text{for $x_1<x_2$ and $t\le s$}\,.\label{3.tmens}
\end{align}
All these inequalities hold as long as all the involved terms
are defined.
\par
The following proposition describes the dynamical behavior
of the differential equation \eqref{3.ecucon} in the vicinity
of hyperbolic solutions.
\begin{prop}\label{3.prophyp}
Let $\wit b$ be a hyperbolic solution
of equation \eqref{3.ecucon}, and let $(k_b,\beta_b)$
be a dichotomy constant pair for $\wit b$.
\begin{itemize}
\item[\rm (i)] If $\wit b$ is
attractive, then for all initial times $s\in\R$ and initial values
$x_0\ge\wit b(s)$, the solution $x(t,s,x_0)$ is defined for any
$t\ge s$, and
\[
 |\wit b(t)-x(t,s,x_0)|\le k_b\,e^{-\beta_b\,(t-s)}|\wit b(s)-x_0|
 \qquad\text{for $t\ge s$}\,;
\]
and, given any $\bar\beta_b\in(0,\beta_b)$,
there exists $\rho>0$ such that, if $s\in\R$ and
$x_0\in[\wit b(s)-\rho,\wit b(s)]$,
then $x(t,s,x_0)$ is defined for any $t\ge s$, and
\[
 |\wit b(t)-x(t,s,x_0)|\le k_b\,e^{-\bar\beta_b\,(t-s)}|\wit b(s)-x_0|
 \qquad\text{for $t\ge s$}\,.
\]
\item[\rm (ii)]
If $\wit b$ is repulsive,
then for any initial time $s\in\R$ and initial value $x_0\le\wit b(s)$,
the solution $x(t,s,x_0)$ is defined for any $t\le s$, and
\[
 |\wit b(t)-x(t,s,x_0)|\le k_b\,e^{\beta_b\,(t-s)}|\wit b(s)-x_0|
 \qquad\text{for $t\le s$}\,;
\]
and given any $\bar\beta_b\in(0,\beta_b)$
there exists $\rho>0$ such that, if $s\in\R$ and
$x_0\in[\wit b(s),\wit b(s)+\rho]$,
then $x(t,s,x_0)$ is defined for any $t\le s$ and
\[
 |\wit b(t)-x(t,s,x_0)|\le k_b\,e^{\bar\beta_b\,(t-s)}|\wit b(s)-x_0|
 \qquad\text{for $t\le s$}\,.
\]
\end{itemize}
\end{prop}
\begin{proof}
%(i)\&(ii)
Let us prove (i). We fix $s\in\R$ and $x_0\ge\wit b(s)$. Since
the function $r$ from Theorem~\ref{3.teoruno} is globally defined and
$\wit b\ge r$, Theorem \ref{3.teoruno}(ii) ensures that
$x(t,s,x_0)$ is defined and bounded (at least) on $[s,\infty)$. Therefore,
the first inequality in (i) follows from the definition of hyperbolicity
and the second inequality in
\eqref{3.tmays}. To prove the second one, we make
the change of variables $z=x-\wit b(t)$, which takes \eqref{3.ecucon} to
\[
 z'=(-2\,\wit b(t)+q(t))\,z-z^2\,.
\]
Let $z(t,s,z_0)$ be the solution of this transformed equation satisfying
$z(s,s,z_0)=z_0$. According to the First Approximation Theorem
(see \cite[Theorem~III.2.4]{hale} and its proof), if
$\bar\beta_b\in(0,\beta_b)$, then
there exists $\rho>0$ such that if $|z_0|\le\rho$, then
$z(t,s,z_0)$ is defined and satisfies $|z(t,s,z_0)|
\le k_b\,e^{-\bar\beta_b\,(t-s)}|z_0|$ for $t\ge s$. The second inequality
in (i) follows from this, since $x(t,s,x_0)=\wit b(t)+z(t,s,x_0-\wit b(s))$.
The proof of (ii) is analogous.
\end{proof}
Proposition~\ref{3.prophyp} shows that hyperbolicity of solutions of
\eqref{3.ecucon} (which is defined by means of the linear variational equation)
has strong implications for the nonlinear differential equation \eqref{3.ecucon}.
If the hyperbolic solution $\wit b$ is attractive, then there exists a
neighbourhood of size $\rho>0$ around the solution curve that is attracted
exponentially uniformly for all times. It follows that this solution is also
locally pullback attractive as defined in Section~\ref{2.sec}, which is a
weaker form of attractivity and requires attraction only in the past and not
uniformly for all times. Similarly, repulsivity of $\wit b$ implies that
$\wit b$ is locally pullback repulsive; see also the discussion in
Section~\ref{2.sec}.
\par
It also follows from Proposition~\ref{3.prophyp} that if the solutions $a$ and $r$
from Theorem~\ref{3.teoruno} are globally defined and hyperbolic, then they must be
uniformly separated in the following sense.
\begin{defi}\label{3.defus}
We say that two globally defined solutions $x_1(t)$ and $x_2(t)$ of \eqref{3.ecucon}
with $x_1\le x_2$ are {\em uniformly separated\/} if $\inf_{t\in\R}(x_2(t)-x_1(t))>0$.
\end{defi}
Conversely, the following theorem shows that if the solutions $a$ and $r$ from
Theorem~\ref{3.teoruno} are globally defined and uniformly separated,
then they are hyperbolic and determine the global dynamics of \eqref{3.ecucon}.
As mentioned before, the pair of solutions $(a,r)$ corresponds to what we have
described in the Introduction as an attractor-repeller pair.
\begin{teor}\label{3.teorhyp}
Suppose that \eqref{3.ecucon} has bounded solutions,
and that the (globally defined) functions $a$ and $r$ provided by
Theorem~\ref{3.teoruno} are uniformly separated. Then,
\begin{itemize}
\item[\rm(i)] the two solutions are hyperbolic, with
$a$ attractive and $r$ repulsive.
\item[\rm(ii)] Let $(k_a,\beta_a)$ and $(k_r,\beta_r)$ be dichotomy constant pairs for
the hyperbolic solutions $a$ and $r$, respectively, and let us choose
any $\bar\beta_a\in(0,\beta_a)$ and any $\bar\beta_r\in(0,\beta_r)$.
Then, given $\ep>0$, there exist $k_{a,\ep}\ge 1$ and
$k_{r,\ep}\ge 1$ (depending also on the choice of $\bar\beta_a$ and
of $\bar\beta_r$, respectively) such that
\[
\begin{split}
 &\quad\qquad |a(t)-x(t,s,x_0)|\le k_{a,\ep}\,e^{-\bar\beta_a(t-s)}|a(s)-x_0|
 \quad\text{if $x_0\ge r(s)+\ep$ and $t\ge s$}\,,\\
 &\quad\qquad |r(t)-x(t,s,x_0)|\le k_{r,\ep}\,e^{\bar\beta_r(t-s)}|r(s)-x_0|
 \quad\text{if $x_0\le a(s)-\ep$ and $t\le s$}\,.
\end{split}
\]
In addition,
\[
\begin{split}
 &\quad\qquad |a(t)-x(t,s,x_0)|\le k_a\,e^{-\beta_a(t-s)}|a(s)-x_0|
 \quad\text{if $x_0\ge a(s)$ and $t\ge s$}\,,\\
 &\quad\qquad |r(t)-x(t,s,x_0)|\le k_r\,e^{\beta_r(t-s)}|r(s)-x_0|
 \quad\text{if $x_0\le r(s)$ and $t\le s$}\,.
\end{split}
\]
\item[\rm(iii)] Equation \eqref{3.ecucon} does not have more
hyperbolic solutions, and
$a$ and $r$ are the unique bounded solutions of
\eqref{3.ecucon} which are uniformly separated.
\end{itemize}
\end{teor}
In the proof of Theorem \ref{3.teorhyp}, we employ a standard technique
for nonautonomous differential equations: the definition of a flow by means
of the hull construction, which allows us to use
techniques from topological dynamics. The proof is quite long and
technical, and the arguments are rather different from the rest of
those used in the paper. That is the reason why we prefer to postpone
this proof until Appendix~\ref{a.appendix}, where we give the definition
of the hull $\W_{q,p}$ of $(q,p)$ and of the continuous flows defined on
$\W_{q,p}$ and on $\W_{q,p}\times\R$.
\par
We point out again that Theorem \ref{3.teorhyp},
which holds for any bounded and uniformly continuous function
$(q,p)\colon\R\to\R\times\R$,
extends previously known properties for the case of {\em recurrent\/}
$(q,p)$; that is, for the case when the flow on $\W_{q,p}$ is minimal.
These properties are proved in \cite{alob3} and \cite{nuos4}.
In fact, our proof is deeply based on the results of \cite{nuos4}:
Theorem~\ref{3.teorhyp} shows that the dynamical description of the
flow on $\W_{q,p}\times\R$ given in \cite{nuos4} for the case
of minimal $\W_{q,p}$ is also valid for the general case of a compact
metric hull $\W_{q,p}$.
\par
We also point out that the paper \cite{anja} develops a
theory for discrete-time skew-product flows over a compact base.
Part of its conclusions are equivalent to many of the
statements made in Theorem~\ref{3.teorhyp} for the
continuous-time case.
%%%%%%%%%%%%%%%%%%%%%%%%%%%%%%%%%%%%%%%%%%%%%%%%%%%%%%%%%%%%%%%%%%%%%%%%%%%%%%%%%
%%%%%%%%%%%%%%%%%%%%%%%%%%%%%%%%%%%%%%%%%%%%%%%%%%%%%%%%%%%%%%%%%%%%%%%%%%%%%%%%%
\subsection{One-parametric variation of the global dynamics.}\label{3.sec33}
Let us now consider the one-parametric family of equations
\begin{equation}\label{3.ecuconl}
 x'=-x^2+q(t)\,x+p(t)+\lb\,,
\end{equation}
where $\lb\in\R$. We will add the subscript $\lb$ to the
previously established notation to refer to the differential
equation (\eqref{3.ecuconl}$_\lb$),
its solutions ($x_\lb(t,s,x_0)$),
the possibly empty set of bounded solutions ($\mB_\lb$), and the
solutions determined in Theorem~\ref{3.teoruno} ($a_\lb$ and $r_\lb$).
The next result shows that \eqref{3.ecuconl} undergoes a bifurcation
at a certain value $\lb=\lb^*$: for $\lb<\lb^*$,
there are no bounded solutions,
while for $\lb>\lb^*$, there exist two bounded hyperbolic solutions.
\begin{teor}\label{3.teorlb*}
There exists a unique $\lb^*=\lb^*(q,p)
\in[-\|q^2/4+p\|,\|p\|]$ such that
\begin{itemize}
\item[\rm(i)] $\mB_\lb$ is empty if and only if $\lb<\lb^*$.
\item[\rm(ii)] If $\lb^*\le\lb_1<\lb_2$,
then $\mB_{\lb_1}\varsubsetneq\mB_{\lb_2}$. More precisely,
\begin{equation}\label{3.chain}
 r_{\lb_2}<r_{\lb_1}\le a_{\lb_1}<a_{\lb_2}\,.
\end{equation}
In addition,
$\lim_{\lb\to\infty}a_\lb(t)=\infty$ and
$\lim_{\lb\to\infty}r_\lb(t)=-\infty$ uniformly on $\R$.
\item[\rm(iii)] If $\lb=\lb^*$, then
$\inf_{t\in\R}(a_{\lb^*}(t)-r_{\lb^*}(t))=0$, and there are
no hyperbolic solutions.
\item[\rm(iv)] If $\lb>\lb^*$, then $a_\lb$ and $r_\lb$ are
uniformly separated and the unique hyperbolic solutions, and
the asymptotic dynamics of
\eqref{3.ecuconl}$_\lb$ is described by Theorems~\ref{3.teoruno}
and \ref{3.teorhyp}.
\item[\rm(v)] $\lb^*(q,p+\lb)=\lb^*(q,p)-\lb$ for any $\lb\in\R$.
\end{itemize}
\end{teor}
\begin{proof}
(i)\&(ii) Note that \eqref{3.ecuconl} reads as
$x'=-(x-q(t)/2)^2+q^2(t)/4+p(t)+\lb$.
It is clear that if $\lb<-\n{q^2/4+p}$
(so that $q^2(t)/4+p(t)+\lb\le-\mu<0$ for
all $t\in\R$), then \eqref{3.ecuconl}$_\lb$
does not any have bounded solution. On the other hand, if $\lb>\n{p}$
(so that $p(t)+\lb\ge\mu>0$ for all $t\in\R$), then
the constant function $b(t)\equiv 0$
satisfies the condition in Theorem \ref{3.teoruno}(v), so that $\mB_\lb$ is
not empty. Let us take $\lb_1<\lb_2$.
Theorem~\ref{3.teoruno}(v) also shows that, if \eqref{3.ecuconl}$_{\lb_1}$
has a bounded solution (so that $r_{\lb_1}$ and $a_{\lb_1}$
are globally defined), then also
$r_{\lb_2}$ and $a_{\lb_2}$
are globally defined and \eqref{3.chain} holds. These facts ensure
that $\mJ:=\{\lb\in\R\,|\;\mB_\lb\text{ is nonempty}\}$ is a nonempty
positive half-line.
Let us define $\lb^*:=\inf\mJ$ and observe that
$\lb^*\in[-\n{q^2/4+p},\n{p}]$. To check now that $\mB_{\lb^*}$ is nonempty,
we define
$\mI_\lb:=\{x_0\in\R\,|\;r_\lb(0)\le x_0\le a_\lb(0)\}$
for $\lb>\lb^*$. The set $\mI_\lb$ is a compact interval,
and $\mI_{\lb_1}\subseteq\mI_{\lb_2}$ if $\lb^*<\lb_1\le\lb_2$,
so that $\mI^*:=\bigcap_{\lb>\lb^*}\mI_\lb$ contains
at least the point $x^*=\lim_{\lambda\to(\lambda^*)^+}a_{\lb}(0)$.
Then, $x_{\lb^*}(t,0,x^*)$ is globally defined and bounded.
To prove this, we take a constant $m$
satisfying simultaneously the condition of \eqref{3.defm}
for all the equations \eqref{3.ecuconl}$_\lb$
corresponding to $\lb\in[\lb^*,\lb^*+1]$, deduce from
\eqref{3.limsup} and \eqref{3.liminf} that
$a_\lambda(t)\in[-m,m]$ for all $\lb\in(\lb^*,\lb^*+1]$,
and note that $x_{\lb^*}(t,0,x^*)=
\lim_{\lambda\to(\lambda^*)^+}x(t,0,a_{\lb}(0))=
\lim_{\lambda\to(\lambda^*)^+}a_{\lb}(t)\in[-m,m]$.
The conclusion is that $\lb^*\in\mJ$, which completes the proof of (i).
\par
To prove the remaining assertion in (ii), we take $n\in\N$ and look
for $\lb(n)$ large enough to guarantee that
$(\pm n)'=0<-n^2\pm q(t)\,n+p(t)+\lb$ if $\lb\ge\lb(n)$.
Theorem \ref{3.teoruno}(v) ensures that $r_\lb(t)<-n<n<a_\lb(t)$
for any $t\in\R$ if $\lb\ge\lb(n)$, which proves the statement.
\smallskip\par
(iii) Since $\mB_\lb$ is empty for $\lb<\lb^*$, Proposition~\ref{3.proppersiste}
precludes the positility of hyperbolic solutions for \eqref{3.ecuconl}$_{\lb^*}$,
and hence Theorem \ref{3.teorhyp} guarantees
that $a_{\lb^*}$ and $r_{\lb^*}$ are not uniformly separated
(even if they are different). This proves (iii).
\smallskip\par
(iv) We assume for contradiction that
$\inf_{t\in\R}(a_\lb(t)-r_\lb(t))=0$ for a $\lb>\lb^*$, which we fix.
It follows from (ii) that
$r_\lb(t)<x_\lb(t,s,a_{\lb^*}(s))<a_\lb(t)$ for any $s,t\in\R$, and
from a standard comparison result that
$d_s(t):=x_\lb(t,s,a_{\lb^*}(s))-a_{\lb^*}(t)\ge 0$ for any
$s\in\R$ and $t\ge s$. We look for a constant $\kappa>0$ such that
$x_\lb(t,s,a_{\lb^*}(s))+a_{\lb^*}(t)-q(t)\le\kappa$ for all $s,t\in\R$.
Then,
\[
 d'_s(t)=-d_s(t)\big(x_\lb(t,s,a_{\lb^*}(s))+a_{\lb^*}(t)-q(t)\big)+\lb-\lb^*
 \ge -\kappa\,d_s(t)+\lb-\lb^* \quad\text{for $t\ge s$}\,,
\]
and hence
\[
 d_s(t)\ge \frac{\lb-\lb^*}{\kappa}\:\big( 1-e^{-\kappa(t-s)}\big)
 \quad \text{whenever $s\in\R$ and $t\ge s$}\,.
\]
In particular, there exists $l>0$ such that
$d_s(t)\ge (\lb-\lb^*)/2\kappa=:\wit \kappa$
whenever $s\in\R$ and $t\ge s+l$.
\par
Now we look for $t_0\in\R$ such that
$a_\lb(t_0)-r_\lb(t_0)<\wit\kappa$. But then, by (ii),
\[
 \wit\kappa>x_\lb(t_0,t_0-l,a_\lb(t_0-l))-r_\lb(t)>
 x_\lb(t_0,t_0-l,a_{\lb^*}(t_0-l))-a_{\lb^*}(t)=d_{t_0-l}(t_0)\,,
\]
which provides the required contradiction.
\smallskip\par
(v) This last assertion follows easily, for instance, from (i).
\end{proof}

We focus now on the differential equation \eqref{3.ecucon},
which coincides with \eqref{3.ecuconl} for $\lb=0$.
\begin{coro}\label{3.coro}
The differential equation \eqref{3.ecucon} has either no hyperbolic
solutions or two hyperbolic solutions, given by the functions
$a$ and $r$ from Theorem~\ref{3.teoruno}. The solutions $a$ and $r$
are uniformly separated and describe the global dynamics for
\eqref{3.ecucon} as explained in Theorems~\ref{3.teoruno} and
\ref{3.teorhyp}.
\end{coro}
\begin{proof}
Theorem~\ref{3.teorlb*} shows the absence of hyperbolic solutions
if $0\le\lb^*(q,p)$ as well as the existence of exactly two of them
satisfying the stated properties if $0>\lb^*(q,p)$.
\end{proof}
We complete this section with two remarks concerning the dynamics
of the parametric family of equations \eqref{3.ecuconl},
as described in Theorem~\ref{3.teorlb*}.
\begin{notas}
(a)~There are well-known examples of equations of the form
\eqref{3.ecuconl} for which the set of bounded solutions is nonempty
and its maximum and minimum solutions $a$ and $r$ are not
uniformly separated, which according to
Theorem~\ref{3.teorlb*} must correspond to $\lb^*(q,p)=0$.
The most classical ones appear in \cite{vino} and \cite{mill,mill2}.
A very precise description of a similar situation
and of the complexity of the induced dynamics can be found
in \cite[Section~8.7]{jonnf}.
\par
(b)~Theorem~\ref{3.teorlb*} can be understood as a result
on nonautonomous saddle-node bifurcation, on the line with those
in \cite{anja} and \cite{nuob6}. In these papers, the skew-product
formalism is used to analyze nonautonomous bifurcations patterns
for one-parametric families of differential equations that are more
general than \eqref{3.ecuconl}. The bifurcation occurs when the set
of bounded solutions is empty at one side of a certain value of the
parameter and contains a classical attractor-repeller pair on the
other side. Detailed descriptions of this situation for some cases of
nonautonomous quadratic differential equations are given in
\cite{jnot} and \cite{fuhr}, where the possibility of occurrence of
strange nonchaotic attractors is carefully analyzed.
\end{notas}
%%%%%%%%%%%%%%%%%%%%%%%%%%%%%%%%%%%%%%%%%%%%%%%%%%%%%%%%%%%%%%%%%%%%%%%%%%%%%%%%%%%%%
%%%%%%%%%%%%%%%%%%%%%%%%%%%%%%%%%%%%%%%%%%%%%%%%%%%%%%%%%%%%%%%%%%%%%%%%%%%%%%%%%%%%%
%%%%%%%%%%%%%%%%%%%%%%%%%%%%%%%%%%%%%%%%%%%%%%%%%%%%%%%%%%%%%%%%%%%%%%%%%%%%%%%%%%%%%
%%%%%%%%%%%%%%%%%%%%%%%%%%%%%%%%%%%%%%%%%%%%%%%%%%%%%%%%%%%%%%%%%%%%%%%%%%%%%%%%%%%%%
\section{Rate-induced tipping}\label{4.sec}
The first goal in this section is to describe in detail the three
(very different) possibilities for the global dynamics induced by
the equation
\begin{equation}\label{4.ecuini}
 y'=-\Big(y-\frac{2}{\pi}\arctan (ct)\Big)^2+p(t)\,,
\end{equation}
where $p\colon\R\to\R$ is a bounded and uniformly continuous function
and $c>0$ (we refer to this differential equation as \eqref{4.ecuini}$_c$).
We study \eqref{4.ecuini}$_c$ under a fundamental hypothesis
(see Hypothesis \ref{4.hipo} below) for the differential equation
\begin{equation}\label{4.ecup}
 x'=-x^2+p(t)\,.
\end{equation}
This differential equation is important for \eqref{4.ecuini}$_c$,
since it relates to the two limit equations of \eqref{4.ecuini}$_c$,
given by  the past equation
\begin{equation}\label{4.ecu-}
 y' =-(y+1)^2+p(t)
\end{equation}
and  the future equation
\begin{equation}\label{4.ecu+}
 y' =-(y-1)^2+p(t)\,.
\end{equation}
(Note that $\lim_{t\to\pm\infty}(\pi/2)\arctan(ct)=\pm 1$ for $c>0$.)
The dynamics induced by the last three equations is the same:
\eqref{4.ecu-} and \eqref{4.ecu+} are obtained from \eqref{4.ecup}
by trivial changes of variable.
\par
The second goal in this section is to analyze the possibility of
occurrence of rate-induced tipping: we study the existence of
critical values $c=c_0$ at which the global dynamics in a
neighborhood of $c_0$ changes from
one of the three previously
described cases (for $c<c_0$ in the neighborhood) to another one
(for $c>c_0$).
It turns out that this rate-induced tipping mechanism is closely
related to the bifurcation analysis performed in Section~\ref{3.sec}.
\par
Let us formulate the aforementioned fundamental hypothesis.
\begin{hipo}\label{4.hipo}
The differential equation \eqref{4.ecup} has exactly two
hyperbolic solutions $\wit a$ and $\wit r$.
\end{hipo}
Corollary \ref{3.coro} ensures that these solutions are uniformly
separated, coincide with those provided by Theorem \ref{3.teoruno}
when applied to the differential equation \eqref{4.ecup}, and
determine the global dynamics for this equation
as described in Theorems~\ref{3.teoruno} and~\ref{3.teorhyp}:
$(\wit a,\wit r)$ is a classical
    attractor-repeller pair for \eqref{4.ecup}.
This implies that, if
\begin{equation}\label{4.solhyp}
 \wma^-:=\wit a - 1\,,\quad\;
 \wmr^-:=\wit r - 1\,,\quad\;
 \wma^+:=\wit a + 1\,,\quad\;\text{and}\quad\;
 \wmr^+:=\wit r + 1\,,
\end{equation}
then $(\wma^-,\wmr^-)$ and $(\wma^+,\wmr^+)$
are classical attractor-repeller pairs for the limit equations
\eqref{4.ecu-} and \eqref{4.ecu+}. Note
finally that if Hypothesis \ref{4.hipo} is not satisfied for a
certain bounded and uniformly continuous function $p$, then we
can replace $p$ with $p_\lb=p+\lb$ for $\lb>\lb^*(0,p)$ in
order to satisfy this hypothesis. This follows from (v) and (iii)
of Theorem~\ref{3.teorlb*}, where the existence of $\lb^*(0,p)$
is established.
\par
The remainder of Section~4 is divided into four parts:
\begin{itemize}
\item[-] In Subsection~\ref{4.sec41}, we formulate the main results
concerning the aforementioned three dynamical possibilities,
and explain them with the help of some descriptive figures.
\item[-] Subsection \ref{4.sec42}, quite technical, is devoted
to prove the main results.
\item[-] In Subsection \ref{4.sec43}, we define a continuous map
$\lb_*\colon[\,0,\infty)\to\R$ such that the global dynamical
behavior for \eqref{4.ecuini}$_c$ with $c>0$ is determined by
the sign of $\lb_*(c)$: the rate-induced tipping mechanism is hence
characterized by changes of sign of $\lb_*$.
\item[-] Finally, Subsection \ref{4.sec44} is devoted to a basic
numerical study of some of the questions treated in the paper.
\end{itemize}
%%%%%%%%%%%%%%%%%%%%%%%%%%%%%%%%%%%%%%%%%%%%%%%%%%%%%%%
%%%%%%%%%%%%%%%%%%%%%%%%%%%%%%%%%%%%%%%%%%%%%%%%%%%%%%%
\subsection{Global dynamics: the main results}\label{4.sec41}
In the description of the three dynamical
possibilities for \eqref{4.ecuini}
(formulated in Theorems \ref{4.teorA}, \ref{4.teorB} and
\ref{4.teorC} below), the most important role is
played by the solutions $\mac\colon \mR^-_c\to\R$ and
$\mrc\colon \mR^+_c\to\R$ of \eqref{4.ecuini}$_c$
determined by Theorem~\ref{3.teoruno}. They exist for any $c>0$ under
Hypothesis~\ref{4.hipo}, which follows from two auxiliary results,
Theorem \ref{4.teorlim} and Proposition \ref{4.propari}.
In the formulation of these results, the four functions
$\wma^\pm$ and $\wmr^\pm$ defined by \eqref{4.solhyp},
that exist under Hypothesis \ref{4.hipo}, will be used.
Recall that $(\wma^-,\wmr^-)$ and
$(\wma^+,\wmr^+)$ are classical attractor-repeller pairs
for the limit equations \eqref{4.ecu-} and \eqref{4.ecu+}, respectively.
The maximal solution of \eqref{4.ecuini}$_c$ is denoted
by $t\mapsto y_c(t,s,y_0)$.
\par
The three aforementioned dynamical possibilities are given by
{\sc Cases A, B} and~{\sc C}:
\begin{defi}\label{4.deficasos}
Given $c>0$, we say that the differential equation \eqref{4.ecuini}$_c$ is
\begin{itemize}
\item[-] in \hypertarget{CA}{{\sc Case A}} if
it has two different hyperbolic solutions,
\item[-] in \hypertarget{CB}{{\sc Case B}} if
it has exactly one bounded solution, and
\item[-] in \hypertarget{CC}{{\sc Case C}} if
it has no bounded solutions.
\end{itemize}
\end{defi}
\begin{teor}\label{4.teorcasos}
Assume that Hypothesis \ref{4.hipo} holds. Then,
\hyperlink{CA}{\sc Cases A}, \hyperlink{CB}{\sc B} and
\hyperlink{CC}{\sc C} exhaust the dynamical possibilities for the
differential equation \eqref{4.ecuini}$_c$ if $c>0$.
\end{teor}
We postpone the proof of this theorem and the following ones until
Section \ref{4.sec42}.
\par
\begin{figure}[ht]
\caption{{\sc Case C}:
the trajectories of $\mac$ (solid red line) and $\mrc$ (long-dashed blue line),
the \lq\lq left half\rq\rq~of the classical attractor-repeller pair
given by ($\wma^-,\,\wmr^-$) for the
limit equation \eqref{4.ecu-} (green short-dashed lines), and
the \lq\lq right half\rq\rq~of the classical attractor-repeller pair
given by ($\wma^+,\,\wmr^+$) for the
limit equation \eqref{4.ecu+} (green dash-dotted lines).}
\label{4.CasoA-limite}
\includegraphics[width=300pt]{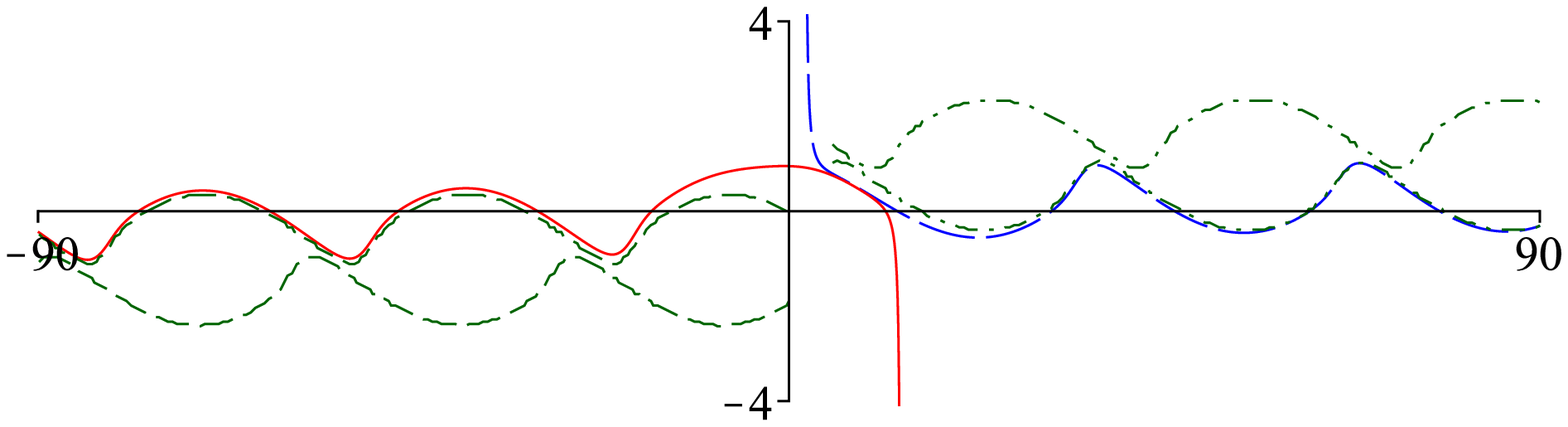}
\end{figure}
\begin{figure}[ht]
\caption{{\sc Case C}: the trajectories
of $\mac$ (solid red line), of $\mrc$ (long-dashed blue line),
two other bounded solutions (indigo short-dashed lines),
and four unbounded solutions approaching $\mrc$ as $t\to-\infty$
or $\mac$ as $t\to\infty$ (black dotted and green dash-dotted lines).}
\label{4.CasoA}
\includegraphics[width=300pt]{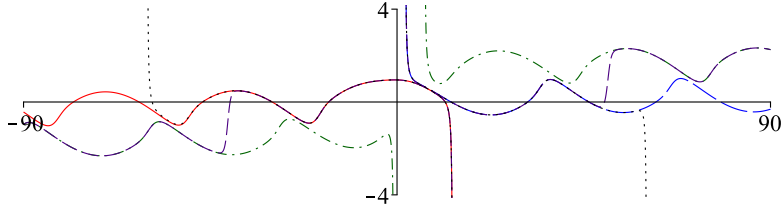}
\end{figure}
\begin{figure}[ht]
\caption{{\sc Case B}:
the trajectories of $\mac=\mrc$ (solid red line),
the \lq\lq left half\rq\rq~of the classical attractor-repeller pair
given by ($\wma^-,\,\wmr^-$) for the
limit equation \eqref{4.ecu-} (green short-dashed lines), and
the \lq\lq right half\rq\rq~of the classical attractor-repeller pair
given by ($\wma^+,\,\wmr^+$) for the
limit equation \eqref{4.ecu+} (green dash-dotted lines).}
\label{4.CasoB-limite}
\includegraphics[width=300pt]{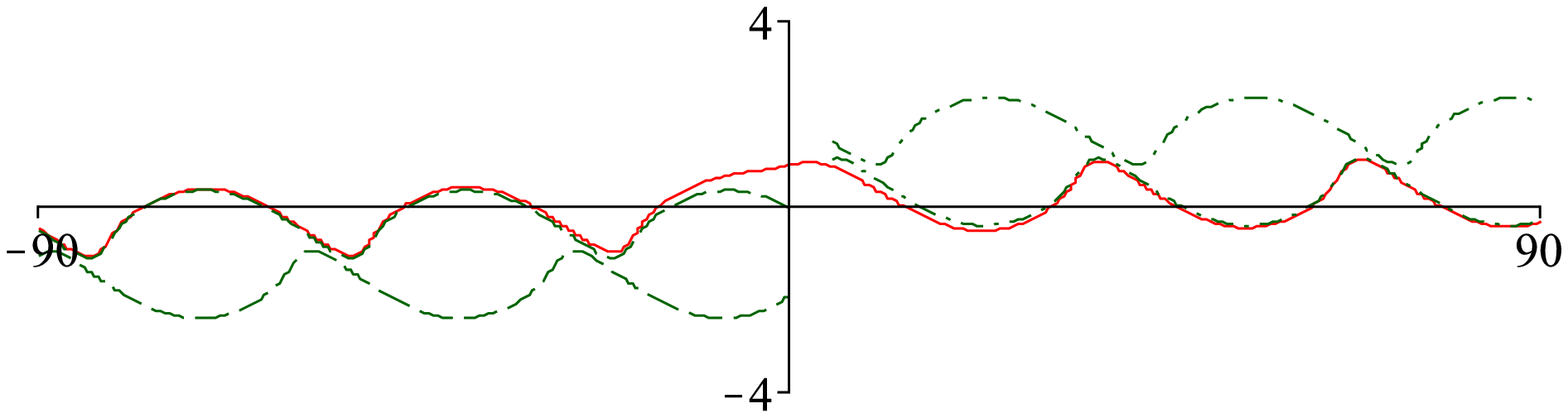}
\end{figure}
\begin{figure}[ht]
\caption{{\sc Case B}:
the trajectories of $\mac=\mrc$ (solid line),
and four solutions left-bounded or right-bounded, which approach
$\wmr^-$ as $t\to-\infty$ or $\wma^+$ as $t\to\infty$.}
\label{4.CasoB}
\includegraphics[width=300pt]{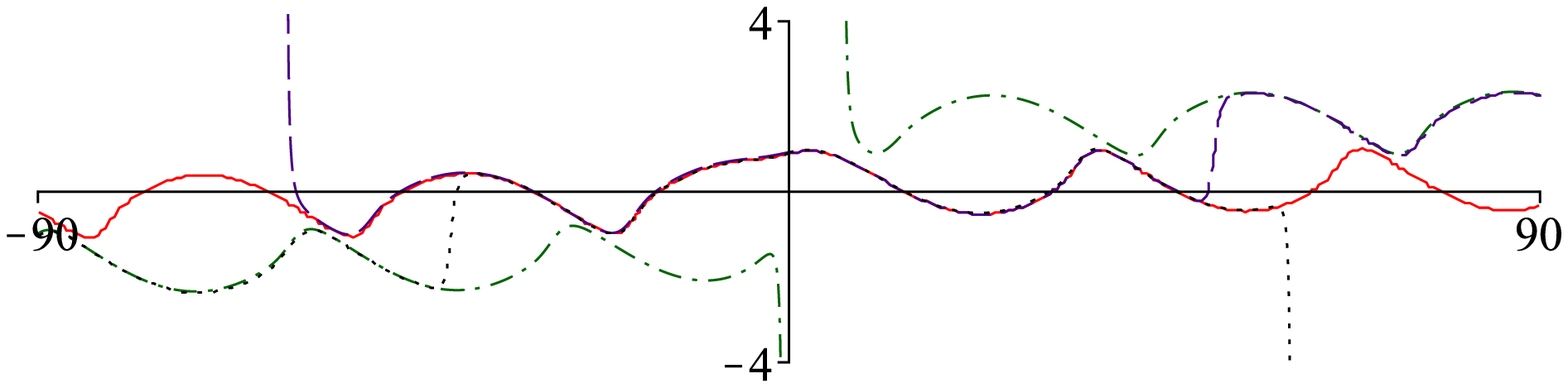}
\end{figure}
\begin{figure}[ht]
\caption{{\sc Case A}:
the trajectories of $\mac$ (solid red line) and $\mrc$ (long-dashed blue line),
the \lq\lq left half\rq\rq~of the classical attractor-repeller pair
given by ($\wma^-,\,\wmr^-$) for the
limit equation \eqref{4.ecu-} (green short-dashed lines), and
the \lq\lq right half\rq\rq~of the classical attractor-repeller pair
given by ($\wma^+,\,\wmr^+$) for the
limit equation \eqref{4.ecu+} (green dash-dotted lines).}
\label{4.CasoC-limite}
\includegraphics[width=300pt]{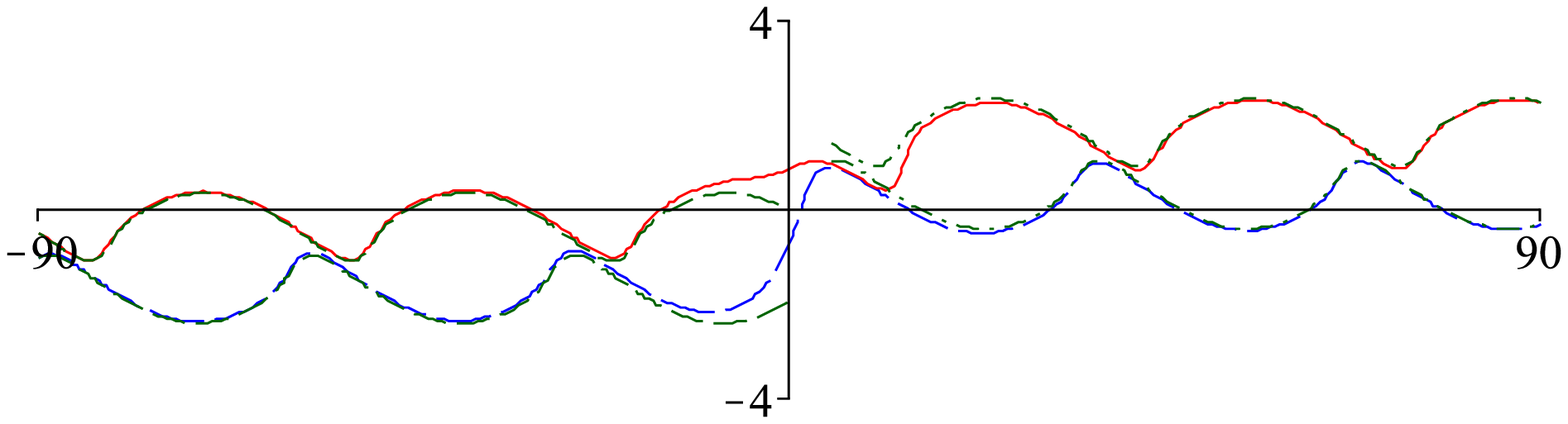}
\end{figure}
\begin{figure}[ht]
\caption{{\sc Case A}:
the trajectories of $\mac$ (solid red line),
$\mrc$ (long-dashed blue line), two solutions unbounded both to
the right and to the left (black dotted lines),
and four solutions left-bounded or right-bounded, which approach
$\wmr^-$ as $t\to-\infty$ or $\wma^+$ as $t\to\infty$
(indigo short-dashed and green
dash-dotted lines).}
\label{4.CasoC}
\includegraphics[width=300pt]{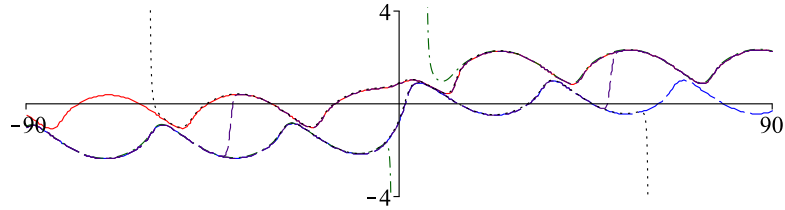}
\end{figure}

In Definition \ref{4.deficasos}, the focus of the classification
is put on the existence of bounded and/or hyperbolic solutions.
But there are many other ways to distinguish the three cases.
Surely, the most meaningful one refers to the behavior of the pair
$(\mac,\mrc)$. These functions are always locally pullback attractive
and repulsive solutions of \eqref{4.ecuini}$_c$,
respectively (see Section \ref{2.sec} for the definition).
After Theorems \ref{4.teorA}, \ref{4.teorB} and \ref{4.teorC}, and
with the help of the Figures \ref{4.CasoA-limite}-\ref{4.CasoC},
it will be clear that, if $c>0$,
\begin{itemize}
\item[-] \hyperlink{CA}{\sc Case A} holds if and only if
$(\mac,\mrc)$ is a classical attractor-repeller pair,
which turns out to be equivalent to saying that $\mac$ and
$\mrc$ are globally defined and different, and in which case
this pair provides the connection between the pairs
$(\wma^-,\wmr^-)$ and $(\wma^+,\wmr^+)$ we referred to in
the Introduction;
\item[-] \hyperlink{CB}{\sc Case B} holds if and only if $\mac=\mrc$,
hence giving  rise to the unique bounded solution, and in which
case this solution connects $\wma^-$ with $\wmr^+$;
\item[-] and \hyperlink{CC}{\sc Case C} holds if and only if
no solution is globally defined, which implies that none of the
above connections can occur.
\end{itemize}
Much more information concerning the limit behavior of all the
solutions of \eqref{4.ecuini}$_c$ is provided by the next three theorems.
Recall that the unbounded solutions are never globally defined: see Theorem
\ref{3.teoruno}.
\begin{teor} \label{4.teorA}
Assume that Hypothesis \ref{4.hipo} holds, and take $c>0$.
In \hyperlink{CA}{{\sc Case A}},
\begin{itemize}
\item[A1] the hyperbolic solutions are given by
$\mrc$ and $\mac$.
\item[A2] If $y_0>\mac(s)$, then $y_c(t,s,y_0)$ is unbounded as $t$
decreases,
and $\lim_{t\to\infty}|\mac(t)-y_c(t,s,y_0)|=0$.
\item[A3] If $\mrc(s)<y_0<\mac(s)$, then $y_c(t,s,y_0)$ is globally
defined and bounded, and
$\lim_{t\to\infty}|\mac(t)-y_c(t,s,y_0)|=
\lim_{t\to-\infty}|\mrc(t)-y_c(t,s,y_0)|=0$.
\item[A4] If $y_0<\mrc(s)$, then $y_c(t,s,y_0)$ is unbounded as
$t$ increases, and
$\lim_{t\to-\infty}|\mrc(t)-y_c(t,s,y_0)|=0$.
\item[A5] $\lim_{t\to\pm\infty}|\wma^\pm(t)-\mac(t)|=0$ and
$\lim_{t\to\pm\infty}|\wmr^\pm(t)-\mrc(t)|=0$.
\end{itemize}
\end{teor}
\begin{teor}\label{4.teorB}
Assume that Hypothesis \ref{4.hipo} holds, and take $c>0$.
In \hyperlink{CB}{{\sc Case B}}, if $\mbc$ is the unique
bounded solution of \eqref{4.ecuini},
\begin{itemize}
\item[B1] $\mrc$ and $\mac$ coincide with $\mbc$.
\item[B2] The solution $\mbc$ is locally
pullback attractive and locally pullback repulsive.
\item[B3] $\lim_{t\to-\infty}|\wma^-(t)-\mbc(t)|=0$ and
$\lim_{t\to\infty}|\wmr^+(t)-\mbc(t)|=0$.
\item[B4] The following three conditions are equivalent: $y_0>\mbc(s)$;
$y_c(t,s,y_0)$ is unbounded as $t$ decreases; and
$\lim_{t\to\infty}|\wma^+(t)-y_c(t,s,y_0)|=0$.
\item[B5] The following three conditions are equivalent:
$y_0<\mbc(s)$; $y_c(t,s,y_0)$ is unbounded as $t$ increases; and
$\lim_{t\to-\infty}|\wmr^-(t)-y_c(t,s,y_0)|=0$.
\end{itemize}
\end{teor}
\begin{teor} \label{4.teorC}
Assume that Hypothesis \ref{4.hipo} holds, and take $c>0$.
In \hyperlink{CC}{{\sc Case C}},
\begin{itemize}
\item[C1] $\mR_c^-=(-\infty,l^-_c)$ for $l^-_c\in\R$ and
$\mR^+_c=(l^+_c,\infty)$ for $l^+_c\in\R$, so that
$\lim_{t\to(l^-_c)^-}\mac(t)=-\infty$ and
$\lim_{t\to(l^+_c)^+}\mrc(t)=\infty$.
\item[C2] The solutions $\mac$ and $\mrc$ of
equation \eqref{4.ecuini}$_c$ are respectively
locally pullback attractive and locally pullback repulsive.
\item[C3] $\lim_{t\to-\infty}|\wma^-(t)-\mac(t)|=0$ and
$\lim_{t\to\infty}|\wmr^+(t)-\mrc(t)|=0$.
\item[C4] Let us take $s\in\mR_c^-$. Then,
$y_c(t,s,y_0)$ is bounded for $t\to-\infty$ if and only if
$y_0\le \mac(s)$;
and $\lim_{t\to-\infty}|\wmr^-(t)-y_c(t,s,y_0)|=0$ if and only if
$y_0<\mac(s)$.
\item[C5] Let us take $s\in\mR_c^+$. Then,
$y_c(t,s,y_0)$ is bounded for $t\to\infty$ if and only if
$y_0\ge \mrc(s)$;
and $\lim_{t\to\infty}|\wma^+(t)-y_c(t,s,y_0)|=0$ if and only if
$y_0>\mrc(s)$.
\item[C6] There exist solutions $y_c(t,s,y_0)$ which are unbounded
at both endpoints of their maximal intervals of definition.
More precisely, this situation corresponds to those points $(s,y_0)$
such that either $s<l_c^-$ and $y_0> \mac(s)$,
or $l_c^-\le s\le l_c^+$ and $y_0\in\R$ (if $l_c^-\le l_c^+$), or
$l_c^+<s$ and $y_0<\mrc(s)$.
\end{itemize}
\end{teor}
Figures \ref{4.CasoA-limite} and \ref{4.CasoA} depict the situation
in \hyperlink{CC}{{\sc Case C}} for  $c=0.15$ and the periodic
function $p(t)=0.9-\sin(t/5)$. Figure \ref{4.CasoA-limite}
shows the curves $\mac$ and $\mrc$ of \eqref{4.ecuini}$_c$,
as well as (\lq\lq a half\rq\rq~of)
the classical attractor-repeller pairs $(\wma^-,\wmr^-)$
and $(\wma^+,\wmr^+)$ for the limit equations \eqref{4.ecu-}
and \eqref{4.ecu+}. Here we observe that $(\mac,\mrc)$
provides the connection between the attractor-repeller
pairs for the limit equations.
Figure~\ref{4.CasoA} shows $\mac$, $\mrc$
and six more solutions of the equation \eqref{4.ecuini}$_c$: two of them
are above $\mac$, so that they are unbounded for $t\to-\infty$ and
approach $\mac$ (and hence $\wma^+$) for  $t\to\infty$; two
of them are below $\mrc$, so that they are unbounded for $t\to\infty$ and
approach $\mrc$ (and hence $\wmr^-$) for $t\to-\infty$; and the remaining
two, globally bounded, approach $\mrc$ (and hence $\wmr^-$) for
$t\to-\infty$ and $\mac$ (and hence $\wma^+$) for $t\to\infty$.
(For our visualization, we have chosen a periodic function $p$
for simplicity. More details concerning the method used to
numerically approximate these curves are given in Section \ref{4.sec44}.)
\par
Figures \ref{4.CasoB-limite} and \ref{4.CasoB} depict the situation in
\hyperlink{CB}{\sc Case B}, for $p(t)=0.9-\sin(t/5)$ and for a value of
$c\approx 0.22609301$, so that $\lb_*(c)=0$. It was numerically
not possible to determine $c$ exactly, and for this reason, the two
figures have been edited to represent \hyperlink{CB}{\sc Case B}:
a small difference in the tenth decimal causes a jump from
\hyperlink{CA}{\sc Case C} to \hyperlink{CC}{\sc Case A}. In Figure
\ref{4.CasoB-limite} we observe how the attractor
$\wma^-$ at $-\infty$, is connected to the repeller $\wmr^+$
at $+\infty$ by the orbit of
the unique bounded solution $\mbc$. Four more solutions are
depicted in Figure \ref{4.CasoB}: those starting above $\mbc$ are
unbounded at $-\infty$ and approach $\wma^+$ at $+\infty$;
and those starting below $\mbc$ are
unbounded at $+\infty$ and approach $\wmr^-$ at $-\infty$.
The locally pullback attractive and repulsive properties of
$\mbc$ are also demonstrated in Figure \ref{4.CasoB}.
\par
Finally, Figures \ref{4.CasoC-limite} and \ref{4.CasoC} depict
\hyperlink{CC}{\sc Case A}, for $p(t)=0.9-\sin(t/5)$ and for
$c=0.15$, following the same scheme of Figures \ref{4.CasoA-limite}
and \ref{4.CasoA}. In this case, the connection between the
attractor-repeller pairs is broken, but still $\wmr^-$ determines the
limit behavior of the solutions bounded for $t\to-\infty$
(except $\mac$, which approaches $\wma^-$); and $\wma^+$
determines the limit behavior of the solutions bounded for
$t\to\infty$ (except $\mrc$, which approaches $\wmr^+$).
The locally pullback attractive (resp.~repulsive) character of $\mac$
(resp.~$\mrc$) is also demonstrated in Figure \ref{4.CasoC}.
\begin{nota}
In the three dynamical possibilities,
if $s\in\mR^+_c$ and $y_0>\mrc(s)$, then the solution $y_c(t,s,y_0)$ is {\em locally
forward attractive}, since $\lim_{t\to\infty}|y_c(t,s,y_0)-y_c(t,s,y_1)|=0$
whenever $y_1>\mrc(s)$. This is due to the fact that
$\lim_{t\to\infty}|y_c(t,s,y_0)-\wma^-(t)|=0$ whenever $y_0>\mrc(s)$,
which is a common property in \hyperlink{CA}{\sc Cases A},
\hyperlink{CB}{B} and \hyperlink{CC}{\sc C}. Similarly,
if $s\in\mR^-_c$ and $y_0<\mac(s)$, then the solution $y_c(t,s,y_0)$ is {\em locally
forward repulsive}, since $\lim_{t\to-\infty}|y_c(t,s,y_0)-y_c(t,s,y_1)|=0$
whenever $y_1<\mac(s)$. This forward behavior can be observed in
Figures \ref{4.CasoA}, \ref{4.CasoB} and \ref{4.CasoC}.
\end{nota}
%%%%%%%%%%%%%%%%%%%%%%%%%%%%%%%%%%%%%%%%%%%%%%%%%%%%%%%%%%%%%%%%%%%%%%%%%%%%%%%%%
%%%%%%%%%%%%%%%%%%%%%%%%%%%%%%%%%%%%%%%%%%%%%%%%%%%%%%%%%%%%%%%%%%%%%%%%%%%%%%%%%
\subsection{Proofs of the main results.}\label{4.sec42}
To analyze the general dynamical properties for the family of
equations \eqref{4.ecuini} is the starting point for the proofs of
Theorems~\ref{4.teorcasos}, \ref{4.teorA}, \ref{4.teorB} and
\ref{4.teorC}. For each $c\ge0$ we make the change of variables
\begin{equation}\label{4.changeofvar}
 x=y-\frac{2}{\pi}\arctan (ct)\,,
\end{equation}
which transforms the differential equation \eqref{4.ecuini}$_c$ to
\begin{equation}\label{4.ecuc}
 x'=-x^2+p(t)-q_c(t)\,,
\end{equation}
where the function $q_c\colon\R\to\R$ is given by
\[
 q_c(t):=\frac{2c}{\pi(c^2t^2+1)}\,.
\]
As usual, \eqref{4.ecuc}$_c$ is the equation of the family \eqref{4.ecuc}
corresponding to a particular value of $c$; and
$t\mapsto x_c(t,s,x_0)$ is the maximal solution with value
$x_0$ at $t=s$. Note that \eqref{4.ecuc}$_0$ coincides with
\eqref{4.ecup}, whose maximal solutions are denoted by $x(t,s,x_0)$.
\par
It is clear that
\[
 \lim_{t\to\pm\infty} q_c(t)=0 \quad\text{and}\quad
 \lim_{c\to 0}q_c(t)=0 \;\;\text{uniformly on $\R$}\,.
\]
This implies that the differential equation \eqref{4.ecup} plays
now three roles: it coincides with \eqref{4.ecuc}$_0$ for $c=0$,
and at the same time it coincides with the past and future limit
equations of \eqref{4.ecuc}$_c$ for any $c>0$, i.e.~for $t\to\infty$
and $t\to-\infty$. The dynamical
properties of \eqref{4.ecuc}$_c$ will be analyzed in Theorems \ref{4.teorlim}
and \ref{4.teorsep}, which will be fundamental
for the proofs of the main results stated in Section~\ref{4.sec41}.
\par
Let $\wit a$ and $\wit r$ be the hyperbolic solutions
from Hypothesis \ref{4.hipo}, and fix $n_0\in\N$ such that
\begin{equation}\label{4.defnl}
 \frac{2}{n_0}\le\,\inf_{t\in\R}(\wit a(t)-\wit r(t))\,.
\end{equation}
The following construction is made for each $n>n_0$ and for any
fixed $c>0$. Proposition \ref{3.proppersiste} applied to
equation \eqref{4.ecup} and to $\ep=1/n$ provides a constant
$\delta_n$, which we also fix, such that $\n{k}\le\delta_n$ ensures that
$x'=-x^2+p(t)+k(t)$ has also two hyperbolic solutions $\wit a_k$ and $\wit r_k$,
with $\wit a_k$ attractive and satisfying $\n{\wit a_k-\wit a}<1/n$,
and with $\wit r_k$ repulsive and satisfying $\n{\wit r_k-\wit r}<1/n$.
In particular, $\wit a_k$ and $\wit r_k$ are different solutions.
In fact, they are uniformly separated.
Assuming without restriction that the sequence $(\delta_n)$ is
decreasing and $\delta_n\le 1/n$, we take
$\delta_{c,n}:=\min(\delta_n,2c/\pi)$
and define $t_{c,n}>0$ as the unique positive value of
time where $q_c(\pm t_{c,n})=\delta_{c,n}$, which implies
\begin{equation}\label{4.proptcn}
 t_{c,n}^2=\frac{2}{\pi\,c\,\delta_{c,n}}-\frac{1}{c^2}
 \qquad\text{and}\qquad
 \lim_{n\to\infty}t_{c,n}=\infty\,.
\end{equation}
And we also define
\[
 \wit q_{c,n}\colon\R\to\R\,,\quad t\mapsto\left\{\begin{array}{ll}
 q_c(t)&\text{if $\;|t|\ge t_{c,n}$}\,,\\
 q_c(t_{c,n})&\text{if $\;|t|\le t_{c,n}$}\,.\end{array}\right.
\]
Hence $\n{\wit q_{c,n}}\le\delta_{c,n}\le\delta_n$, so that the equation
\begin{equation}\label{4.ecukn}
 x'=-x^2+p(t)-\wit q_{c,n}(t)
\end{equation}
(to which we will refer as \eqref{4.ecukn}$_{c,n}$,
and whose maximal solutions will be represented by $x_{c,n}(t,s,x_0)$)
has two different hyperbolic solutions, $\wit a_{c,n}$ and $\wit r_{c,n}$,
with
\begin{equation}\label{4.cerca}
 \inf_{t\in\R}(\wit a_{c,n}-\wit r_{c,n})>0\,,\quad
 \n{\wit a_{c,n}-\wit a}\le \frac{1}{n}\quad\text{and}\quad
 \n{\wit r_{c,n}-\wit r}\le \frac{1}{n}\,.
\end{equation}
Now we define $a^-_{c,n}$ and $r^+_{c,n}$
as the unique (possibly locally defined) solutions of \eqref{4.ecuc}$_c$
with $a^-_{c,n}(-t_{c,n})=\wit a_{c,n}(-t_{c,n})$
and $r^+_{c,n}(t_{c,n})=\wit r_{c,n}(t_{c,n})$, and observe that
$a^-_{c,n}$ is at least defined on $(-\infty,-t_{c,n}\,]$
and $r^+_{c,n}$ is at least defined on $[\,t_{c,n},\infty)$, and
\begin{equation}\label{4.lados}
\begin{split}
 &a^-_{c,n}(t)=\wit a_{c,n}(t)\quad\text{for $\,t\le -t_{c,n}$}\,,\\
 &r^+_{c,n}(t)=\wit r_{c,n}(t)\quad\text{for $\,t\ge t_{c,n}$}\,,\\
 &x_c(t,s,x_0)=x_{c,n}(t,s,x_0) \quad\text{if $\,s,t\ge t_{c,n}\,$ or
 $\,s,t\le-t_{c,n}$}\,.
\end{split}
\end{equation}
Finally, we denote by $a_c\colon\mR^-_c\to\R$ and $r_c\colon\mR^+_c\to\R$
the solutions associated to \eqref{4.ecuc}$_c$ by Theorem~\ref{3.teoruno},
and we note that the domains
$\mR^-_c$ and $\mR^+_c$ (as defined in \eqref{3.defsets}) are nonempty,
due to the existence and properties of $a^-_{c,n}$ and $r^+_{c,n}$.
Recall that these sets contain at least a negative and a positive half line,
respectively, and that they are given by the entire line $\R$ if and only if
\eqref{4.ecuc}$_c$ has globally bounded solutions (see Theorem~\ref{3.teoruno}).
Note also that $a_c$ and $r_c$  are closely related to the above mentioned solutions
$\mac$ and $\mrc$ of \eqref{4.ecuini}$_c$, which is explained in
Proposition~\ref{4.propari} below.
In the statements and the proofs of the following results, there appear
the following eight functions (six in fact, as Theorem \ref{4.teorlim} shows):
\begin{itemize}
\item[-] $\wit a$ and $\wit r$ (hyperbolic
solutions of \eqref{4.ecup}=\eqref{4.ecuc}$_0$),
\item[-] $a_c$, $a^-_{c,n}$, $r_c$ and $r^+_{c,n}$
(solutions of \eqref{4.ecuc}$_c$, perhaps locally defined),
\item[-] and $\wit a_{c,n}$ and $\wit r_{c,n}$
(hyperbolic solutions of \eqref{4.ecukn}$_{c,n}$).
\end{itemize}
Recall that the definitions of locally pullback attractive and repulsive
solutions are given in Section \ref{2.sec}.
\begin{teor}\label{4.teorlim}
Assume that Hypothesis \ref{4.hipo} holds.
Let us take $n_0$ satisfying \eqref{4.defnl} and fix some $n>n_0$. For any $c>0$,
\begin{itemize}
\item[\rm(i)] $\wit r<\wit r_{c,n}<\wit a_{c,n}<\wit a$.
\item[\rm(ii)] $a^-_{c,n}=a_c$ and $r^+_{c,n}=r_c$.
\item[\rm(iii)] $\lim_{t\to-\infty}|\wit a(t)-a_c(t)|=0$ and
$\lim_{t\to\infty}|\wit r(t)-r_c(t)|=0$.
\item[\rm(iv)] Let us take $s\in\mR_c^-$. Then,
$x_c(t,s,x_0)$ is bounded for  $t\to -\infty$ if and only if $x_0\le a_c(s)$;
and $\lim_{t\to-\infty}|\wit r(t)-x_c(t,s,x_0)|=0$ if and only if
$x_0<a_c(s)$.
\item[\rm(v)] Let us take $s\in\mR_c^+$. Then,
$x_c(t,s,x_0)$ is bounded for $t\to\infty$ if and only if $x_0\ge r_c(s)$;
and $\lim_{t\to\infty}|\wit a(t)-x_c(t,s,x_0)|=0$ if and only if
$x_0>r_c(s)$.
\item[\rm(vi)] The solutions $a_c$ and $r_c$ of equation \eqref{4.ecuc}$_c$
are respectively locally pullback attractive and locally
pullback repulsive.
\end{itemize}
\end{teor}
\begin{proof}
(i) The second inequality follows, for instance, from \eqref{4.cerca}
and \eqref{4.defnl}, since $n>n_0$.
Theorem~\ref{3.teoruno}(v) ensures the other ones, since
$0<\wit q_{c,n}$.
\smallskip\par
(ii) The last equality in \eqref{4.lados} yields
$x_c(t,-t_{c,n},x_0)=x_{c,n}(t,-t_{c,n},x_0)$ for $t\le-t_{c,n}$.
Theorem \ref{3.teoruno}(i) applied to equation \eqref{4.ecukn}$_{c,n}$
and to its solution $\wit a_{c,n}$, and the first and last equalities in
\eqref{4.lados}, yield the equivalence
\begin{multline*}
 \sup_{t\in(-\infty,-t_{c,n}]} x_c(t,-t_{c,n},x_0)=
 \sup_{t\in(-\infty,-t_{c,n}]} x_{c,n}(t,-t_{c,n},x_0)<\infty\\
  \Leftrightarrow\;x_0\le\wit a_{c,n}(-t_{c,n})=a^-_{c,n}(-t_{c,n})\,.
\end{multline*}
But this is exactly the definition of $a_c(-t_{c,n})$, as
Theorem \ref{3.teoruno}(i) guarantees. Therefore,
$a^-_{c,n}(t)=a_c(t)$, since both of them solve \eqref{4.ecuc}$_c$
and they coincide at $t=-t_{c,n}$. An analogous argument shows that
$r^+_{c,n}(t)=r_c(t)$.
\smallskip\par
(iii) Property (ii) and the first equality in \eqref{4.lados}
guarantee that $a_c$ and $\wit a_{c,n}$ coincide on
$(-\infty,-t_{c,n}]$. Since $(t_{c,n})\uparrow\infty$ as $n\to\infty$
(see \eqref{4.proptcn}),
the first assertion in (iii) is a consequence of \eqref{4.cerca}.
The argument is analogous for the second one.
\smallskip\par
(iv) We fix $s\in\mR_c^-$.
Theorem \ref{3.teoruno}(i) proves the first assertion in (iv).
Let us take $x_0<a_c(s)$, and any $\ep>0$.
We look for $n>n_0$ (with $n_0$ defined by \eqref{4.defnl})
such that $1/n\le\ep/2$ and such that $-t_{c,n}<s$
(see again \eqref{4.proptcn}). Then,
$x_c(-t_{c,n},s,x_0))<a_c(-t_{c,n})=\wit a_{c,n}(-t_{c,n})$.
In addition, the last equality in \eqref{4.lados} ensures that
$x_c(t,s,x_0)=x_{c,n}(t,-t_{c,n},x_c(-t_{c,n},s,x_0))$.
Now we write
\[
\begin{split}
 |\wit r(t)-x_c(t,s,x_0)|
 &\le|\wit r(t)-\wit r_{c,n}(t,s,x_0)|\\
 &\quad +|\wit r_{c,n}(t,s,x_0)-x_{c,n}(t,-t_{c,n},x_c(-t_{c,n},s,x_0))|
\end{split}
\]
and apply Theorem \ref{3.teorhyp}(ii) to the hyperbolic solutions of
\eqref{4.ecukn}$_{c,n}$ and the last bound in \eqref{4.cerca}
in order to conclude the existence of $t_n\le -t_{c,n}$ such that
$|\wit r(t)-x_c(t,s,x_0)|<\ep$ for all $t\le t_n$. Therefore,
$\lim_{t\to-\infty}|\wit r(t)-x_c(t,s,x_0)|=0$. This limit behavior
is precluded for $x_0=a_c(s)$ by the first property in (iii)
and the uniform separation of $\wit a$ and $\wit r$; and for
$x_0>a_c(s)$ by the fact that the maximal domain of
$x_c(t,s,x_0)$ is bounded below, according to Theorem~\ref{3.teoruno}(iii).
\smallskip\par
(v) The proof of these assertions is similar to that of point (iv).
\smallskip\par
(vi) We observe that the first assertion in (iii) and properties
\eqref{4.cerca} allow us to choose $s_0\in\mR^-_c$
such that $\ep:=(1/2)\inf_{s\in(-\infty,s_0]}(a_c(s)-\wit r_{c,n}(s))>0$.
Hence, Theorem \ref{3.teoruno}(ii) applied to equation \eqref{4.ecukn}$_{c,n}$
ensures that its solutions $x_{c,n}(t,s,a_c(s)\pm\ep)$ are defined for any
$t\ge s$ if $s\le s_0$.
We assume without restriction that $s_0\le -t_{c,n}$.
\par
Now we fix $t\le s_0$ and take $s\le t$.
If $l\in[s,t]$, then
$a_c(l)=\wit a_{c,n}(l)$ (due to (ii) and the first equality in
\eqref{4.lados}), and this implies $x_c(l,s,a_c(s)\pm\ep)=
x_{c,n}(l,s,a_c(s)\pm\ep)=x_{c,n}(l,s,\wit a_{c,n}(s)\pm\ep)$
(using also the last equality in \eqref{4.lados}).
Therefore, Theorem \ref{3.teorhyp}(ii) applied to the
hyperbolic solutions of \eqref{4.ecukn}$_{c,n}$, using the $\ep$ defined above,
provides $k_0\ge 1$ and $\beta_0>0$ (independent of $s$) with
\[
 |a_c(t)-x_c(t,s,a_c(s)\pm\ep)|
 =|\wit a_{c,n}(t)-x_{c,n}(t,s,\wit a_{c,n}(s)\pm\ep)|
 \le k_0\,e^{-\beta_0(t-s)}\ep\,,
\]
which is as small as desired if $-s$ is large enough.
This proves (vi) in the case of $a_c$, and
the argument is similar for $r_c$.
\end{proof}
\par
Much more information can be given in the case of
global existence of $a_c$ and $r_c$, mainly if they are different.
Recall that $\lb^*(0,p-q_c)$ is associated to
equation \eqref{4.ecuc}$_c$ by Theorem \ref{3.teorlb*}.
\begin{teor}\label{4.teorsep}
Assume that Hypothesis \ref{4.hipo} holds.
Let us take $n_0$ satisfying \eqref{4.defnl} and fix some $n>n_0$.
Assume also that $a_c$ and $r_c$ are globally defined and bounded. Then,
if $c>0$,
\begin{itemize}
\item[\rm(i)] $\wit r<\wit r_{c,n}\le r_c\le a_c\le\wit a_{c,n}<\wit a$.
\end{itemize}
If, in addition, $a_c$ and $r_c$ are different, then
\begin{itemize}
\item[\rm(ii)] $\lim_{t\to\infty}|\wit a(t)-a_c(t)|=0$ and
$\lim_{t\to-\infty}|\wit r(t)-r_c(t)|=0$.
\item[\rm(iii)] $\inf_{t\in\R}(a_c-r_c)>0$,
$a_c$ and $r_c$ are hyperbolic solutions, and $\lb^*(0,p-q_c)<0$.
\end{itemize}
Consequently,
\begin{itemize}
\item[\rm(iv)] if $\lb^*(0,p-q_c)=0$, then equation \eqref{4.ecuc}$_c$
has only one bounded solution.
\end{itemize}
\end{teor}
\begin{proof}
(i) Theorem~\ref{3.teoruno}(v) ensures the chain of inequalities, since
$0<\wit q_{c,n}\le q_c$ and since the two hyperbolic
solutions $\wit a_{c,n}$ and $\wit r_{c,n}$ of \eqref{4.ecukn}$_{c,n}$
delimit the set of initial data of bounded solutions for
this equation (see Corollary \ref{3.coro}), as
$\wit a $ and $\wit r$ do for \eqref{4.ecup}.
\smallskip\par
(ii) Let us take $\ep>0$ and fix $n_1>\max(2/\ep,n_0)$. Since $a_c>r_c$,
it follows from (i) that
$a_c(t_{c,n_1})>\wit r_{c,{n_1}}(t_{c,n_1})$.
Therefore, the last relation in \eqref{4.lados} and Theorem
\ref{3.teorhyp}(ii)
ensure that $|a_c(t)-\wit a_{c,{n_1}}(t)|\le\ep/2$
if $t$ is large enough. Using \eqref{4.cerca}, we conclude that,
for these values of $t$, $|\wit a(t)- a_c(t)|\le 1/n_1+\ep/2\le\ep$,
which proves the result for $a_c$. The proof is analogous for $r_c$.
\smallskip\par
(iii) The uniform separation between $a_c$ and $r_c$
follows from Theorem \ref{4.teorlim}(iii), assertion (ii), and
the uniform separation between $\wit a$ and $\wit r$ guaranteed
by Theorem~\ref{3.teorhyp}(iii).
Therefore, Theorem~\ref{3.teorhyp}(i)
ensures that $a_c$ and $r_c$
are hyperbolic solutions of \eqref{4.ecuc}.
And Theorem \ref{3.teorlb*}(i)\&(iii) ensure that
$\lb^*(0,p-q_c)<0$.
\smallskip\par
(iv) Theorem~\ref{3.teorlb*}(i) ensures the existence of
at least one bounded solution, and (iii) precludes the
existence of two different bounded solutions.
\end{proof}
The next results add some information about how the domains $\mR_c^-$ and
$\mR_c^+$ of $a_c$ and $r_c$ depend on $c>0$. This is helpful
in order to justify the accuracy of the numerical simulations
performed at the end of the paper.
\begin{prop}\label{4.proptstar}
Assume that Hypothesis \ref{4.hipo} holds.
There exists $t^*>0$ independent of $c>0$
such that, if $q^*_c(t)$ is defined
as $q_c(t)$ for $|t|>t^*$ and as $q_c(t^*)$ on $[-t^*,t^*]$,
then the differential equation $x'=-x^2+q^*_c(t)+p(t)$ admits two different
hyperbolic solutions $\wit a^*_c$ and $\wit r^*_c$ such that
$a_c=\wit a^*_c$ on $(-\infty,-t^*]$ and
$r_c=\wit r^*_c$ on $[t^*,\infty)$. In particular,
$(-\infty,-t^*]\subset\mR_c^-$ and $[t^*,\infty)\subset\mR^+_c$
for any $c>0$.
\end{prop}
\begin{proof}
Using the trivial inequality $|2\,\alpha\beta|/(\alpha^2+\beta^2)\le 1$,
we check that $|2\,t\,q_c(t)|\le 2/\pi$ for all $t\in\R$ and $c>0$,
so that $0< q_c(t)\le 1/(\pi|t|)$ for all $t\ne 0$ and $c>0$.
Proposition \ref{3.proppersiste} allows us to choose
$t^*>0$ large enough to guarantee that, if $\n{q}\le 1/(\pi\,t^*)$, then
$x'=-x^2+q(t)+p(t)$ has two different hyperbolic solutions (as close to
those of $x'=-x^2+p(t)$ as desired). The function $q_c^*$
of the statement satisfies this condition. We
call the corresponding upper and lower hyperbolic
solutions $\wit a_c^*$ and $\wit r_c^*$,
and repeat the arguments leading to the first and second
equalities in \eqref{4.lados} and to Theorem \ref{4.teorlim}(ii)
in order to complete the proof.
\end{proof}
\begin{notas}\label{4.notatstar}
(a)~As a consequence of the previous result, we can assert that
$a_c$ and $r_c$ are globally defined hyperbolic solutions
if and only if they are respectively defined at least
on $(-\infty,t^*]$ and $[-t^*,\infty)$ and
satisfy $a_c>r_c$ on $[-t^*,t^*]$. The \lq\lq only if\rq\rq~is trivial,
and to check the \lq\lq if\rq\rq, we must just realize
that this situation precludes the existence of
a vertical asymptote for $a_c$ or $r_c$, since they should
respectively correspond to values of $t>t^*$ or $t<-t^*$
(and hence the graphs of $a_c$ and $r_c$ would intersect).
In fact, if is enough to find a point $t\in[-t^*,t^*]$
at which $a_c(t)>r_c(t)$.
\par
(b)~The convergence of the solution $x_c^*(t,s,x_0)$ of
$x'=-x^2+q_c^*(t)+p(t)$ to the hyperbolic solution
$\wit a_c^*$ (or $\wit r_c^*$)
is exponentially fast for $t\to\infty$ if $x_0\ge\wit a^*(s)$
(or for $t\to-\infty$ if $x_0\le\wit r^*(s)$), as
Theorem \ref{3.teorhyp}(ii) states. Therefore,
$x_c^*(-t^*,-s,x_0)$ will approach
$\wit a_c^*(-t^*)=a_c(-t^*)$ as close as required by choosing
$-s$ much smaller than $-t^*$ and $x_0\ge\wit a_c^*(-s)=a_c(-s)$.
In fact, a computer does not distinguish between
$x_c^*(-t^*,-s,x_0)$ and $a_c(-t^*)$ if $s-t^*$ is
large enough due to limited precision. A similar situation applies to
$x_c^*(t^*,s,x_0)$ and $r_c(t^*)$ if $s$ is
much larger than $t^*$ and $x_0\le r_c(s)$.
\end{notas}
To prove the main theorems, we need to introduce more notation
and go deeper to understand  the relation
between \eqref{4.ecuini}$_c$ and \eqref{4.ecuc}$_c$.
Let $y_c(t,s,y_0)$ be the solution of
\eqref{4.ecuini}$_c$ with $y_c(s,s,y_0)=y_0$. It is immediate to check that
\begin{equation}\label{4.xy}
 y_c(t,s,y_0)=x\big(t,s,y_0-(\pi/2)\arctan(cs)\big)+(\pi/2)\arctan(ct)\,.
\end{equation}
Recall that the functions $a_c$ and $r_c$, associated to \eqref{4.ecuc}$_c$,
are defined on the sets $\mR^-_c$ and $\mR^+_c$, respectively.
\begin{prop}\label{4.propari}
Consider the original equation
\eqref{4.ecuini}$_c$, and let $\mac$ and $\mrc$ be the solutions
provided by Theorem~\ref{3.teoruno}. Then,
\begin{itemize}
\item[\rm(i)] the map $\mac$ is defined on $\mR^-_c$ and
$\mac(t)=a_c(t)+(\pi/2)\arctan(ct)$, and
\item[\rm(ii)] the map $\mrc$ is defined on $\mR^+_c$ and
$\mrc(t)=r_c(t)+(\pi/2)\arctan(ct)$.
\end{itemize}
\end{prop}
\begin{proof}
Let us define $\mdc(t)=a_c(t)+(\pi/2)\arctan(ct)$ and observe that it
is a solution of \eqref{4.ecuini}$_c$ defined exactly on $\mR^-_c$.
In addition, according to \eqref{4.xy} and Theorem~\ref{3.teoruno}(i),
$y(t,s,y_0)$ is bounded for $t\to-\infty$ if and only if
$y_0-(\pi/2)\arctan(cs)\ge a_c(s)$; that is, if and only
if $y_0\ge\mdc$. Using again Theorem~\ref{3.teoruno}(i),
we conclude that $\mdc=\mac$, which completes the proof of (i).
The proof of (ii) is analogous.
\end{proof}
Most of the work is done now, and the proofs of the main theorems,
stated in Section \ref{4.sec41}, follow directly from the previous results.
\medskip\par
\noindent{\em Proof of Theorem \ref{4.teorcasos}}.~It follows
from the definition (see Section \ref{2.sec}) that $\wit b$ is a
hyperbolic solution for \eqref{4.ecuc}$_c$ if and only if
$\wit{\mathfrak{b}}$ is a hyperbolic solution for
\eqref{4.ecuini}$_c$, where $\wit{\mathfrak{b}}(t)=
\wit b(t)+(\pi/2)\arctan(ct)$. An analogous assertion holds when
\lq\lq hyperbolic\rq\rq~is replaced by \lq\lq bounded\rq\rq.
Therefore, the classification of Definition \ref{4.deficasos}
is equivalent for both equations.
According to Theorem~\ref{3.teorlb*}, we are in
\hyperlink{CC}{\sc Case A} if $\lb^*(0,p-q_c)<0$, and in
\hyperlink{CA}{\sc Case C} if $\lb^*(0,p-q_c)>0$.
Theorem~\ref{4.teorsep}(iv) shows
\hyperlink{CB}{\sc Case B} holds if $\lb^*(0,p-q_c)=0$,
and this completes the proof.\hfill{\qed}
\medskip\par
We prove Theorem \ref{4.teorC} prior to Theorems \ref{4.teorB} and
\ref{4.teorA}.
\medskip\par
\noindent{\em Proof of Theorem \ref{4.teorC}}.~
C1 follows from Proposition \ref{4.propari} and
Theorem~\ref{3.teoruno}(iv),(i)\&(ii).
To prove C2, note that it follows from Proposition~\ref{4.propari}
and \eqref{4.xy} that
\[
 |\mac(t)-y_c(t,s,\mac(s)\pm\delta)|=
 | a_c(t)-x_c(t,s,a_c(s)\pm\delta)|
\]
and
\[
 |\mrc(t)-y_c(t,s,\mrc(s)\pm\delta)|=
 |r_c(t)-x_c(t,s,r_c(s)\pm\delta)|\,,
\]
so that Theorem~\ref{4.teorlim}(vi) proves
this assertion. Similarly,
C3, C4 and C5 follow from the statements (iii), (iv) and (v) of
Theorem~\ref{4.teorlim}
combined with Proposition~\ref{4.propari} and equalities
\eqref{4.solhyp} and \eqref{4.xy}.
Finally, C6 follows easily from Theorem \ref{3.teoruno}
applied to \eqref{4.ecuini}$_c$.
\hfill{\qed}
\medskip\par
\noindent{\em Proof of Theorem \ref{4.teorB}}.~
Property B1 is trivial, and using this,
we can repeat the arguments of Theorem~\ref{4.teorC} to prove B2 and B3.
The fist equivalence in B4 is proved by Theorem \ref{3.teoruno}(i)
and the equality $\mbc=\mac$; the second one by Theorem~\ref{4.teorlim}(v)
and the equalities $\mbc=\mrc$, \eqref{4.xy} and \eqref{4.solhyp}.
The proof of B5 is analogous.\hfill{\qed}
\medskip\par
\noindent{\em Proof of Theorem \ref{4.teorA}}.~
It follows from Theorem \ref{3.teorlb*} that the existence of
two hyperbolic solutions for \eqref{4.ecuc}$_c$
corresponds to the case $\lb^*(0,p-q_c)<0$, in which case
these solutions are $\mac$ and $\mrc$. This fact,
Proposition \ref{4.propari} and Remark~\ref{4.notaequiv} prove A1.
\par
Once this is established, Theorem \ref{3.teoruno}(i)\&(ii) and
Theorem~\ref{4.teorlim}(iv)\&(v)
combined with Proposition \ref{4.propari} and
\eqref{4.xy} prove A2, A3 and A4; and Theorems \ref{4.teorlim}(iii)
and \ref{4.teorsep}(ii),
Proposition \ref{4.propari} and equalities \eqref{4.solhyp} prove A5.
\hfill{\qed}
%%%%%%%%%%%%%%%%%%%%%%%%%%%%%%%%%%%%%%%%%%%%%%%%%%%%%%%%%%%%%%%%%%%%%%%%%%%%%%%%%
%%%%%%%%%%%%%%%%%%%%%%%%%%%%%%%%%%%%%%%%%%%%%%%%%%%%%%%%%%%%%%%%%%%%%%%%%%%%%%%%%
\subsection{The bifurcation curve $\lb_*$ and rate-induced tipping}\label{4.sec43}
In this subsection, we analyse rate-induced tipping occurring in
\eqref{4.ecuini}$_c$ under variation of the rate $c>0$. We note that
\eqref{4.ecuini}$_c$ is linked to the differential equation
\eqref{4.ecuc}$_c$ by means of the change of variables
\eqref{4.changeofvar}. In particular, this implies that for any
fixed value of $c> 0$, the differential equations \eqref{4.ecuini}$_c$
and \eqref{4.ecuc}$_c$ share the same dynamics by being either in
\hyperlink{CA} {\sc Case A}, \hyperlink{CB}{B} or \hyperlink{CC}{C}
from Definition \ref{4.deficasos}.
\par
For this reason, rate-induced tipping in \eqref{4.ecuini}$_c$ occurs
if and only if \eqref{4.ecuc}$_c$ admits a bifurcation. In other words,
if the absence of bounded solutions for \eqref{4.ecuc}$_c$ gives rise
to the presence of an attractor-repeller pair as $c$ increases or decreases.
According to Theorem \ref{3.teorlb*}, this corresponds to a change of the
sign of $\lb^*(0,p-q_c)$, which is exactly the bifurcation point in $\lb$
of the differential equation
\begin{equation}\label{4.ecucl}
 x'=-x^2+p(t)-q_c(t)+\lb\,.
\end{equation}
To analyze this in the context of rate-induced tipping for \eqref{4.ecuini}$_c$,
we define
\[
 \lb_*\colon[0,\infty)\to\R\,,\qquad c\mapsto \lb_*(c):=\lb^*(0,p-q_c)
\]
and note that $\lb_*(0)=\lb^*(0,p)$.
We first establish elementary properties of  $\lb_*$ under Hypothesis
\ref{4.hipo}. Note that the existence of hyperbolic solutions ensured
by this hypothesis and Theorem ~\ref{3.teorlb*} guarantee that $\lb_*(0)<0$.
\begin{teor}\label{4.teorcont}
Assume that Hypothesis \ref{4.hipo} holds.
\begin{itemize}
\item[\rm(i)] The map $\lb_*$ is Lipschitz continuous, with
$|\lb_*(c_2)-\lb_*(c_1)|\le (2/\pi)|c_2-c_1|$ for $c_1\ge 0$ and $c_2\ge 0$,
and it takes values in the interval $[\,\lb_*(0)\,,\n{p}+1\,]$.
\item[\rm(ii)] If the equation $x'=-x^2+p(t)+\lb_*(0)$ has
only one bounded solution, then $\lb_*(c)>\lb_*(0)$ for any $c>0$.
\end{itemize}
\end{teor}
\begin{proof}
Theorem \ref{3.teorlb*}(i) ensures that there exists at
least a bounded solution $b_c$ for  $x'=-x^2+p(t)-q_c(t)+\lb_*(c)$.
Let us take $c>0$. Since $q_c>0$, we have $b_c'(t)<-b_c^2(t)+p(t)+\lb_*(c)$,
and hence the last assertion in Theorem \ref{3.teoruno}(v)
ensures that the equation $x'=-x^2+p(t)+\lb_*(c)$ has two different
bounded solutions. This fact combined with Theorem~\ref{3.teorlb*}
proves that $\lb_*(c)\ge \lb^*(0,p)=\lb_*(0)$, and that the inequality
is strict under the additional condition assumed in (ii), which proves (ii).
\par
The change of variables $y=x+(2/\pi)\arctan(ct)$ takes equation
\eqref{4.ecucl}
to $y'=-(y-(2/\pi)\arctan(ct))^2+p(t)+\lb$, and clearly preserves the
property of occurrence or absence of bounded solutions: in this regard,
the role of $\lb_*(c)$ is the same for both equations.
In particular, Theorem \ref{3.teorlb*} shows that
$\lb_*(c)\le\sup_{t\in\R}|p(t)-(4/\pi^2)\arctan^2(ct)|
\le\n{p}+1$.
\par
In order to prove the Lipschitz continuity, we fix $c_1\ge0$ and $c_2\ge 0$,
and take a bounded solution $b_{c_1}$ for the equation
$x'=-x^2+p(t)-q_{c_1}(t)+\lb_*(c_1)$. Then,
$b_{c_1}'(t)\le -b_{c_1}^2(t)+p(t)-q_{c_2}(t)+\n{q_{c_1}-q_{c_2}}+\lb_*(c_1)$,
so that Theorem \ref{3.teoruno}(v) and Theorem~\ref{3.teorlb*}(i) ensure that
$\lb_*(c_2)\le \n{q_{c_1}-q_{c_2}}+\lb_*(c_1)$, that is,
$\lb_*(c_2)-\lb_*(c_1)\le \n{q_{c_1}-q_{c_2}}$. Interchanging the
roles of $c_1$ and $c_2$ we find
$\lb_*(c_1)-\lb_*(c_2)\le \n{q_{c_1}-q_{c_2}}$, so that
$|\lb_*(c_2)-\lb_*(c_1)|\le \n{q_{c_1}-q_{c_2}}$.
It is very easy to check that $(\partial/\partial c) q_c(t)\le 2/\pi$.
Altogether, we conclude that $|\lb_*(c_2)-\lb_*(c_1)|\le (2/\pi)|c_2-c_1|$,
which proves the assertion.
\end{proof}
\begin{nota}\label{4.notaequiv}
Since the change of variables
$y=x-(2\pi)\arctan(ct)$ does not change the possible boundedness
or hyperbolicity of the solutions, the
dynamical possibilities for the equation
\begin{equation}\label{4.ecuclini}
 y'=-(y-(2/\pi)\arctan(ct))^2+p(t)+\lb
\end{equation}
are those three described by Theorem~\ref{3.teorlb*},
and they depend on the relation between $\lb$ and $\lb^*(0,p-q_c)=\lb_*(c)$.
\par
Therefore, the graph of the map $\lb_*\colon[0,\infty)\to\R$
is the bifurcation curve for the two-parametric families of equations
\eqref{4.ecucl} and \eqref{4.ecuclini}: for pairs $(c,\lb)$ above the
graph, two hyperbolic solutions exist; for pairs $(c,\lb)$ below the
graph, no bounded solutions exist; and for the points of the graph,
there exist bounded solutions, but none of them is hyperbolic.
These assertions follow from Theorem \ref{3.teorlb*}.
(In fact, Theorem \ref{4.teorsep}(iv) ensures the existence
of only one bounded solution for each point on the graph if $c>0$.)
\end{nota}
The following proposition is a reformulation of Theorem~\ref{4.teorcasos}
and shows that the sign of $\lb_*$ describes in which of the three cases
the differential equation \eqref{4.ecuini}$_c$ is. The statement follows
from the proof of Theorem~\ref{4.teorcasos}.
\begin{prop}\label{4.propmainth}
Assume that Hypothesis \ref{4.hipo} holds. Then, for $c>0$, the differential
equation  \eqref{4.ecuini}$_c$ is
\begin{itemize}
\item[\rm(i)] in \hyperlink{CA}{\sc Case A}
if and only if $\lb_*(c)<0$,
\item[\rm(ii)]  in \hyperlink{CB}{\sc Case B}
if and only if $\lb_*(c)=0$,
\item[\rm(iii)]  in \hyperlink{CC}{\sc Case C}
if and only if $\lb_*(c)>0$.
\end{itemize}
\end{prop}
The proposition implies that under Hypothesis \ref{4.hipo},
a change of sign of the function $\lambda_*$ describes rate-induced
tipping. Note that it is possible that the function $\lambda_*$ takes
only strictly negative values, in which case tipping does not
occur and we are in \hyperlink{CA}{{\sc Case A}} for any rate
$c>0$.
\par
Let us suppose now that rate-induced tipping is possible in the
sense that there exists a $c>0$ with $\lambda_*(c)>0$. Due to
$\lambda_*(0)<0$ and continuity of $\lambda_*$, this implies
the existence of $c_0>0$ with $\lambda_*(c_0)=0$. We suppose
that rate-induced tipping occurs \emph{transversally} at $c=c_0$,
i.e.~there exists $\delta_0>0$ such that  $\lb_*(c)<0$
for $c\in[c_0-\delta_0,c_0)$ and $\lb_*(c)>0$
for $c\in(c_0,c_0+\delta_0]$. Using Proposition~\ref{4.propmainth},
transversal rate-induced tipping means that
\begin{itemize}
\item[-] for $c\in[c_0-\delta_0,c_0)$, the differential equation
\eqref{4.ecuini}$_c$ is in \hyperlink{CA}{{\sc Case A}},
\item[-] for $c=c_0$, the differential equation \eqref{4.ecuini}$_c$
is in \hyperlink{CB}{{\sc Case B}}, and
\item[-] for $c\in(c_0,c_0+\delta_0]$, the differential equation
\eqref{4.ecuini}$_c$ is in  \hyperlink{CC}{{\sc Case C}}.
\end{itemize}
The following theorem shows that in this situation, tipping can
be described by a collision of the attractor-repeller pair $(\mac,\mrc)$.
\begin{teor}\label{4.teortip}
We suppose that Hypothesis \ref{4.hipo} holds, and that there
exists a $c_0>0$ such that \eqref{4.ecuini}$_c$ admits a transversal
rate-induced tipping at $c=c_0$.
\begin{itemize}
\item[\rm(i)] As $c$ increases to $c_0$, the upper and lower bounds
of the set of bounded solutions collide in the sense that
\begin{equation}\label{4.ulti}
\begin{split}
 &\lim_{c\to c_0^-}\mac=\mathfrak{b}_{c_0}\quad
 \text{uniformly on the negative half-lines of $\R$}\,,\\
 &\lim_{c\to c_0^-}\mrc=\mathfrak{b}_{c_0}\quad
 \text{uniformly on the positive half-lines of $\R$}\,.\\
\end{split}
\end{equation}
\item[\rm(ii)] As $c$ decreases to $c_0$, the
half-line domains $\mR^-_c\!$ and $\mR^+_c\!$ of
$\mac$ and $\mrc$ satisfy
\[
 \lim_{c\to c_0^+}\sup\mR_c^-=\infty \qquad\text{and}
 \lim_{c\to c_0^+}\inf\mR_c^+=-\infty\,,
\]
and the limits \eqref{4.ulti} hold also as $c\to c_0^+$.
\end{itemize}
\end{teor}
\begin{proof}
We will prove the results for the auxiliary equations \eqref{4.ecuc},
and apply Proposition \ref{4.propari} to deduce them for \eqref{4.ecuini},
since $\arctan(ct)$ converges to $\arctan(c_0t)$ uniformly on $\R$ as
$c\to c_0$. Let $n_0$ satisfy \eqref{4.defnl},
and let us fix $n>n_0$ and consider the equations
\eqref{4.ecukn}$_{c,n}$ for $c\in[c_0-\delta_0,c_0+\delta_0]$,
for which the classical attractor-repeller pairs $(\wit a_{c,n},\wit r_{c,n})$
exist. It follows from \eqref{4.proptcn}
that $\lim_{c\to c_0} t_{c,n}=t_{c_0,n}$. We set $t_0=t_{c_0,n}+1$. Then,
$\wit a_{c,n}(t)=a_c(t)$ for $t\le -t_0$ and $c\in[c_0-\delta,c_0+\delta_0]$
(see Theorem  \ref{4.teorlim}(ii) and the first equality in \eqref{4.lados}).
Recall also that $a_{c_0}=b_{c_0}$. Proposition \ref{3.proppersiste}
allows us to assert that
$\lim_{c\to c_0} a_{c,n}=a_{c_0,n}$ uniformly on $\R$, so that
$\lim_{c\to c_0} a_{c}=b_{c_0}$ uniformly on $(-\infty,-t_0]$.
\par
On the other hand, $b_{c_0}$ is also defined on $[-t_0,\infty)$.
Let us fix $k\in\N$. The theorem of
continuous dependence with respect to initial conditions and parameters
provides $\delta_k\in(0,\delta]$ such that if $c\in[c_0,c_0+\delta_k]$ then
$a_c(t)=x_c(t,-t_0,a_c(-t_0))$ is defined on $[-t_0,k]$, and in addition
$\lim_{c\to c_0^+}a_c=b_{c_0}$ uniformly on $[-t_0,k]$. This shows that
$\sup\mR^-_c$ tends to $\infty$ as $c\to c_0^+$ (recall that
$\mR_c^-=\mR_c^+=\R$ for $c\in[c_0-\delta_0,c_0]$), and that
$\lim_{c\to c_0}a_c=b_{c_0}$ uniformly on $(-\infty,k]$. The assertions
of (i) and (ii) concerning $a_c$ (and hence $\mac$) are hence proved.
And, as usual, the proof is analogous for $\mrc$.
\end{proof}
\par
We can reach similar conclusions if $c_0$ is a point at which the graph of
$\lb_*$ crosses transversally the horizontal axis in a decreasing sense:
that is, $\lb_*(c_0)=0$, and there exists $\delta_0\in(0,c_0)$
such that $\lb_*(c)>0$ for $c\in[c_0-\delta_0,c_0)$ and $\lb_*(c)<0$
for $c\in(c_0,c_0+\delta_0]$. The difference is
that now the situation changes from
\hyperlink{CC}{{\sc Case C}}  for $c<c_0$
to \hyperlink{CA}{{\sc Case A}} for $c>c_0$.
\par
On the other hand, if $0$ is a strict local maximum of $\lambda_*$ at $c_0$,
then \hyperlink{CA}{{\sc Case A}} holds for $c\not=c_0$ in a neighborhood
of $c_0$, with the same limit behavior for $|\mac-\mrc|$ as
that described in Theorem \ref{4.teortip}(i)
(from both sides). Accordingly, if $0$ is a strict local minimum attained
at $c_0$, then no bounded solutions exist for values of $c\ne c_0$ close
to $c_0$, and the limit behavior for $|\mac-\mrc|$ is that of
Theorem \ref{4.teortip}(ii) (from both sides).
Finally, we point out that the four cases that we
have mentioned do not exhaust the possibilities for the set of zeros of $\lb_*$.
%%%%%%%%%%%%%%%%%%%%%%%%%%%%%%%%%%%%%%%%%%%%%%%%%%%%%%%%%%%%%%%%%%%%%%%%%%%%%%%%%%%%%
%%%%%%%%%%%%%%%%%%%%%%%%%%%%%%%%%%%%%%%%%%%%%%%%%%%%%%%%%%%%%%%%%%%%%%%%%%%%%%%%%%%%%
%%%%%%%%%%%%%%%%%%%%%%%%%%%%%%%%%%%%%%%%%%%%%%%%%%%%%%%%%%%%%%%%%%%%%%%%%%%%%%%%%%%%%
%%%%%%%%%%%%%%%%%%%%%%%%%%%%%%%%%%%%%%%%%%%%%%%%%%%%%%%%%%%%%%%%%%%%%%%%%%%%%%%%%%%%%
\subsection{Numerical simulations}\label{4.sec44}
In this final subsection, we provide a numerical analysis of some of the questions
treated in this paper for the differential equation \eqref{4.ecuini}, where
\begin{equation}\label{4.eqpt}
 p(t):=0.895-\sin(t/2)-\sin(\sqrt{5}\,t)\,.
\end{equation}
We have chosen the value $0.895$ to capture the possible scenario in which, upon the increase of the parameter $c$, the system recovers the connection between the attractor-repeller pairs of the past and future equation.
We used the MATLAB function {\tt ode45} for numerical approximations of
all the involved equations, with the options on the relative and
absolute tolerance {\tt RelTol=1e-9} and {\tt AbsTol=1e-9}.
\par
We intend to provide a illustration of a number of rate-induced
tipping phenomena, and to do so, we have to work under two fundamental
assumptions, for which we have consistent numerical evidences
(see end of the section): firstly, that there exists
an attractor-repeller pair $(\wit a_*,\wit r_*)$ for the modified
differential equation
\begin{equation}\label{4.ecupm}
 x'=-x^2+p(t)-0.03\,;
\end{equation}
and secondly, that, for a dichotomy pair $(k,\beta)$ which is simultaneously
valid for the two hyperbolic solutions,
\begin{equation}\label{4.cond2}
 4\,k\,e^{-950\,\beta}< 10^{-16}\,.
\end{equation}
Note that the first condition means that $\lb_*(0)<-0.03$, so that
Hypothesis \ref{4.hipo} is fulfilled. That is, there exists
an attractor-repeller
pair $(\wit a,\wit r)$ for \eqref{4.ecup}. In turn, under this hypothesis,
Theorem \ref{4.teorlim} (resp.~Proposition \ref{4.propari}) ensures the
existence of the possibly locally defined solutions
$a_c$ and $r_c$ (resp.~$\mac$ and $\mrc$) associated by
Theorem~\ref{3.teoruno} to \eqref{4.ecuc}$_c$ (resp.~to
\eqref{4.ecuini}$_c$), for any $c\ge 0$. In addition, since
the constant $m=2$ satisfies the condition \eqref{3.defm} for
all the differential equations~\eqref{4.ecuc}$_c$, we know that $a_c(t)<2$,
$\mac(t)<3$, $r_c>-2$ and $\mrc(t)>-3$ on their respective
domains, see Theorem \ref{3.teoruno} and Proposition \ref{4.propari}.
The second assumption \eqref{4.cond2}
is used to find suitable initial conditions (initial time, initial value) to
obtain suitable approximations of the solutions $\mac$ and $\mrc$.
\par
We proceed with the representation of $\mac$ and $\mrc$, using the
idea given in Remark~\ref{4.notatstar}(b).
The first point is showing
that, for the procedure described in Proposition \ref{4.proptstar},
we can take $t^*=50$.
Note first that $1/(50\pi)<0.01$, which ensures that
$q_c(t)<0.01$ if $|t|\ge 50$ for any $t\ge 0$
(see the proof of Proposition \ref{4.proptstar}).
Now we construct $q_c^*$ as in the statement of Proposition \ref{4.proptstar},
so that $0\le q_c^*(t)< 0.01< 0.03$, and consider
\begin{equation}\label{4.ecu*}
 x'=-x^2-q^*_c(t)+p(t)\,,
\end{equation}
whose maximal solution is denoted by $x_c^*(t,s,x_0)$. We can use
Theorem~\ref{3.teoruno}(v) to compare \eqref{4.ecu*} with
\eqref{4.ecupm}, which combined with the fact that $m=2$
satisfies the condition \eqref{3.defm} for
all the equations \eqref{4.ecu*}$_c$,
allows us to assert that the solutions
$\wit r_c^*$ and $\wit a_c^*$ of \eqref{4.ecu*} from
Theorem \ref{3.teoruno} are globally defined and satisfy
$-2<\wit r_c^*<\wit r_*<\wit a_*<\wit a_c^*<2$
for any $c\ge 0$.
Since $\wit a_*$ and $\wit r_*$ are uniformly separated, so are $\wit a_c^*$
and $\wit r_c^*$, and hence, Theorem \ref{3.teorhyp} ensures that
$(\wit a_c^*,\wit r_c^*)$ is an attractor-repeller pair for \eqref{4.ecu*}.
In addition, the former inequalities make it easy
to check that the dichotomy constant pair
$(k,\beta)$ is also valid for $(\wit a_c^*,\wit r_c^*)$, for any $c\ge 0$.
\par
Having in mind these facts, the information
provided by Theorem \ref{3.teorhyp}(ii), and assumption \eqref{4.cond2},
we observe that the computer (working with double precision)
distinguishes neither $\wit a_c^*(-50)$
from $x_c^*(-50,-1000,2)$, nor $\wit r_c^*(50)$ from $x_c^*(50,1000,-2)$.
(Note that $4$ is a bound for $|2-\wit a_c^*(-1000)|$ and for
$|-2-\wit r_c^*(1000)|$.) Recall now that the solutions of
\eqref{4.ecuc}$_c$ are $x_c(t,s,x_0)$, and note that
$x^*_c(t,-1000,2)=x_c(t,-1000,2)$ and $a_c(t)=\wit a^*_c(t)$ for
$t\le -50$, and $x^*_c(t,1000,-2)=x_c(t,1000,2)$ and
$r_c(t)=\wit r^*_c(t)$ for $t\ge 50$: this can be proved as
\eqref{4.lados} and Theorem~\ref{4.teorlim}(ii).
Altogether, we can assert that the computer distinguishes neither $a_c(-50)$
from $x_c(-50,-1000,2)$, nor $r_c(50)$ from $x_c^*(50,1000,-2)$.
\par
All this information can be immediately transferred to the
differential equations \eqref{4.ecuini}$_c$: their
solutions $y_c(t,-1000,3)$ and $y_c(t,1000,-3)$ are respectively
adequate to obtain representations of
$\mac(t)$ for $t\ge -50$ and of $\mrc(t)$ for $t\le 50$.
\par
The way to proceed is clear now. If there is $t_a>-50$ such that
$y_c(t_a,-1000,3)<-3$, then there are not bounded solutions:
the graph of any possible bounded solution lies in
$\R\times [-3,3]$, and hence, such a bounded solution would intersect
the graph of the solution $y_c(t_a,-1000,-3)$.
Therefore, the dynamics is given by  \hyperlink{CC}{\sc Case C}.
And if we can continue the solution $y_c(t,-1000,3)$ at least
until $t=50$, and observe that
$y_c(50,-1000,3)>y_c(50,1000,-3)$, then we
are in \hyperlink{CA}{\sc Case A} (see Remark~\ref{4.notatstar}(a)).
\par
\begin{figure}
\caption{
Approximation of $\mac$ (in solid red) and $\mrc$
(in dashed blue) of \eqref{4.ecuini}$_c$ with $p$ given as
in \eqref{4.eqpt}, for ten values of $c$ between $0$ and $9$.
Three tipping points have been detected with a precision of
five digits by an approximation of the solutions  $\mac$ and $\mrc$
before and after the bifurcation. The pictures for $c=0.5$, $c=1.2$
and $c=9$ demonstrate how $\mac$ and $\mrc$ are separated when the
value of $c$ is not close to a tipping point.}
\includegraphics[trim={3cm 0.4cm 2.6cm -0.1cm},clip,width=\textwidth]{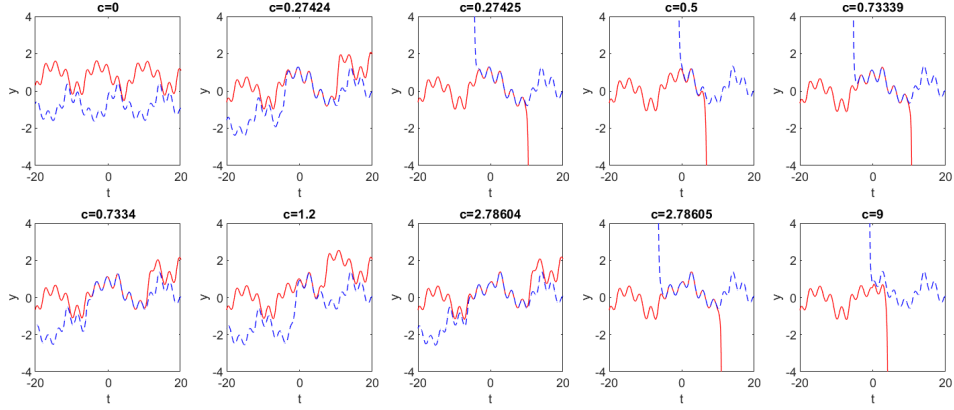}
\label{fig:c-sweep}
\end{figure}

In Figure~\ref{fig:c-sweep}, the solutions $\mac$ and $\mrc$ are plotted
for certain increasing values of~$c$ in the $(t,y)$-plane.
As expected, these two solutions are bounded for small values $c\ge0$,
with $\mac>\mrc$, since the system is in \hyperlink{CA}{{\sc Case A}}.
We observe rate-induced tipping as $c$ increases, when
the solutions $\mac$ and $\mrc$ become unbounded and hence
the system passes to be in \hyperlink{CC}{{\sc Case C}}.
Interestingly, we observe that a further increase in $c$ brings the
system back into \hyperlink{CA}{{\sc Case A}}, and then another transition
into \hyperlink{CC}{{\sc Case~C}} can be observed.
Note also that, thanks to Theorems~\ref{4.teorA},
\ref{4.teorB} and \ref{4.teorC}, we know that
$\lim_{t\to-\infty}|\wma^-(t)-\mac(t)|=
\lim_{t\to\infty}|\wmr^+(t)-\mrc(t)|=0$ holds
for all $c>0$. This limit behavior can
also be observed in Figure~\ref{fig:c-sweep}, as well
as the fact that $\lim_{t\to\infty}|\wma^+(t)-\mac(t)|=
\lim_{t\to-\infty}|\wmr^+(t)-\mrc(t)|=0$
when $c>0$ and $\mac$ and $\mrc$ are globally defined and
uniformly separated.
\par
\begin{figure}
\caption{Computation of the graph of $\lambda_*(c)$ for
\eqref{4.ecuini}
with $p$ given as in \eqref{4.eqpt} for $c\in[0,3.5]$.
The picture on the right is a magnification of the picture
on the left close to $0$. It seems very plausible that
a slight change of the constant 0.895
appearing in the definition of $p$ could affect the number of
possible tipping points for this system.}
\includegraphics[trim={0.3cm 0cm 0.9cm -0.2cm},clip,width=0.49\textwidth]{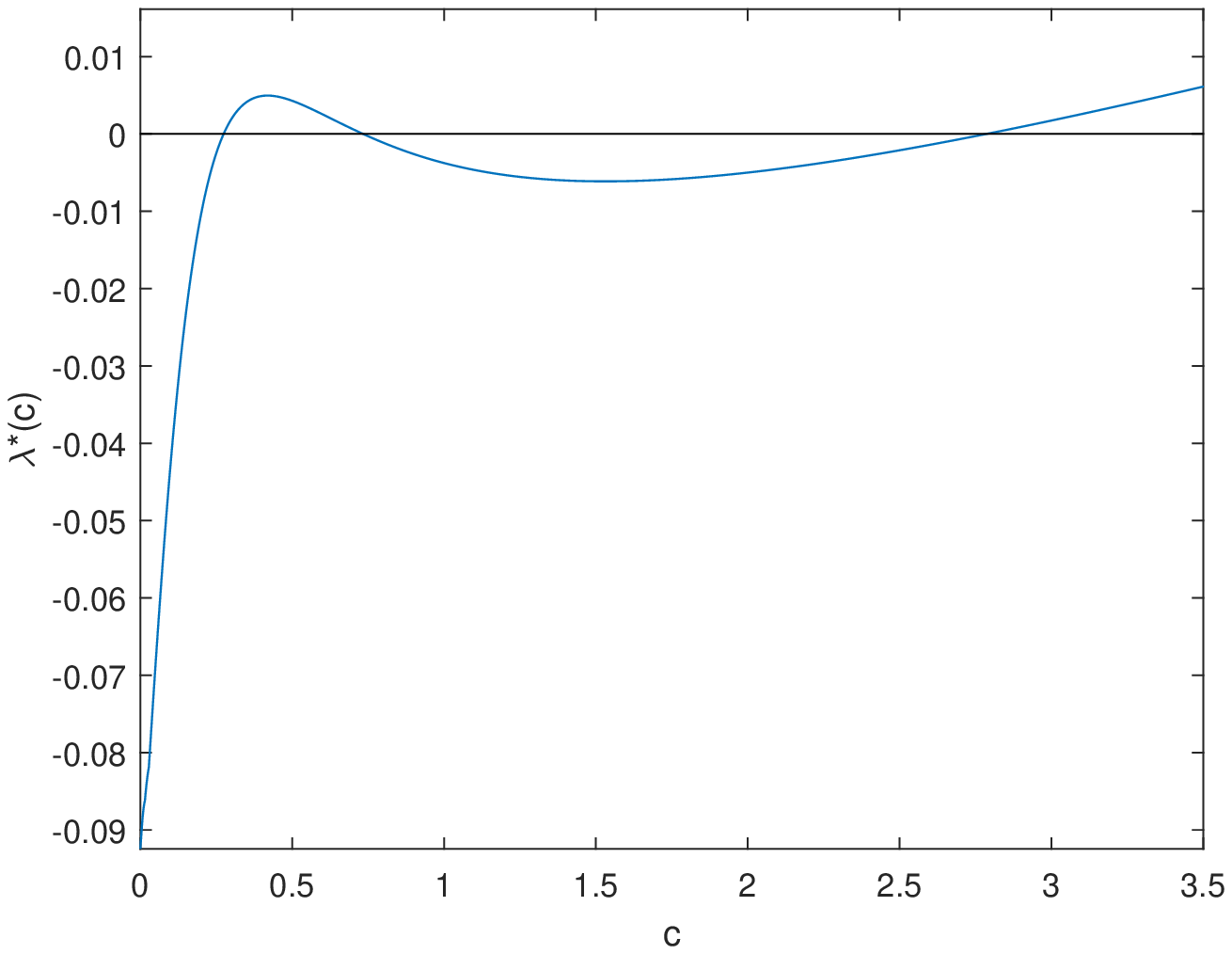}
\includegraphics[trim={0.3cm 0cm 0.9cm -0.2cm},clip, width=0.49\textwidth]{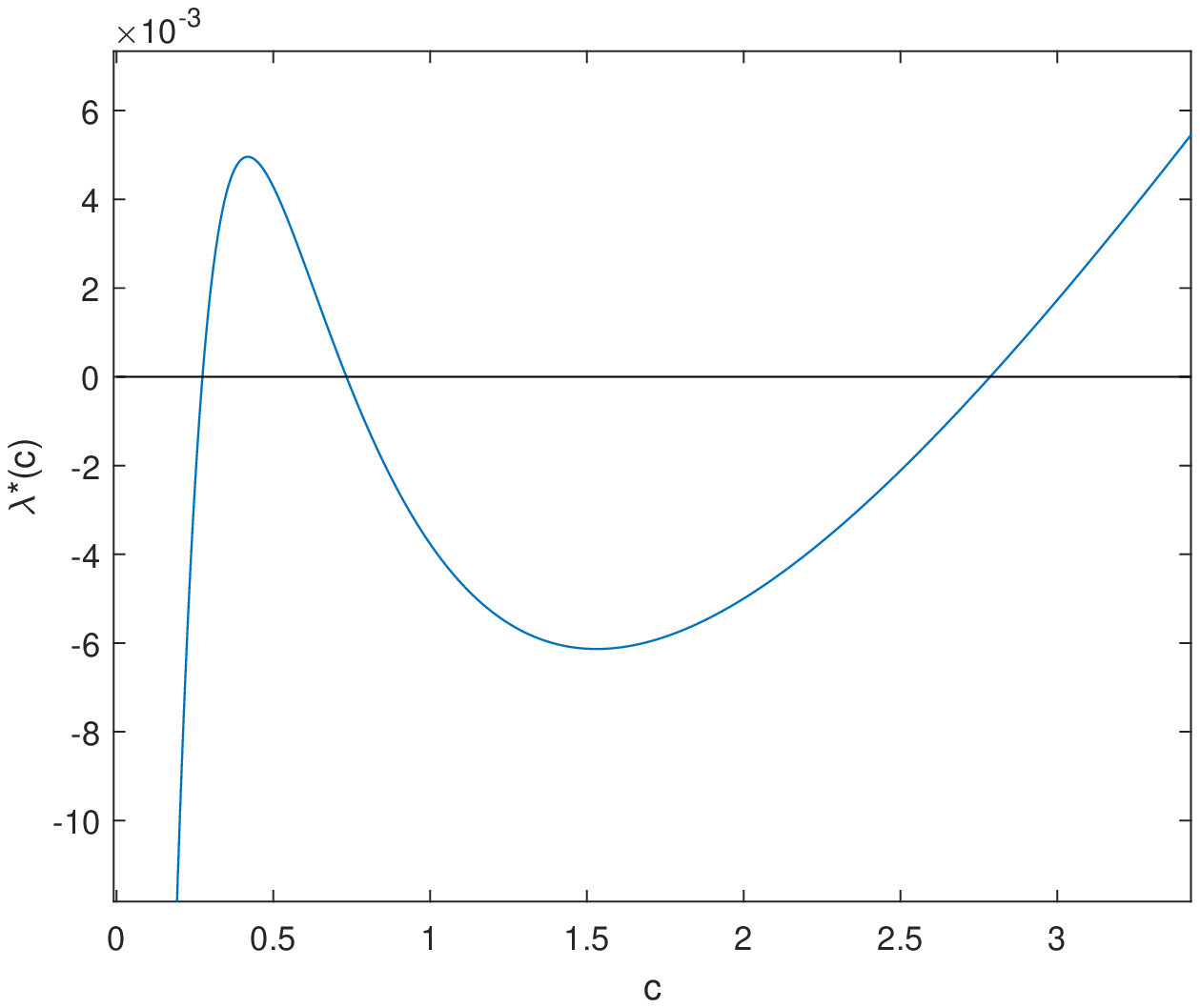}
\label{fig:lambda}
\end{figure}
\begin{figure}
\caption{Approximation of the solutions  $\mac$ (in solid red)
and $\mrc$ (in dashed blue) of \eqref{4.ecuini} with
$p$ given as in \eqref{4.eqpt}, in the vicinity of the first
tipping point in Figure \ref{fig:lambda}.}
\includegraphics[trim={3.5cm 0cm 2.9cm -0.6cm},clip,
width=\textwidth]{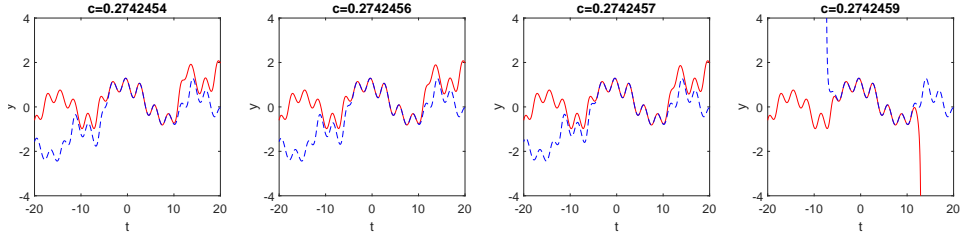}
\label{fig:point1}
\end{figure}

The observed transitions \hyperlink{CA}{{\sc Cases A}} $\rightarrow$
\hyperlink{CC}{{\sc C}} $\rightarrow$ \hyperlink{CA}{{\sc A}} $\rightarrow$
\hyperlink{CC}{{\sc C}}
(which can only occur if the situation fits the critical
\hyperlink{CB}{{\sc Case B}} at least for three values of
$c$) suggest that the critical tipping rate
may neither be just related to its magnitude nor unique.
In order to better understand the occurrence of the different instances of tipping, in Figure \ref{fig:lambda},
we show the graph of the map~$\lb_*$, which has been introduced before
Theorem~\ref{4.teorcont} (see also Remark~\ref{4.notaequiv}).
The mapping $\lb_*$ is computed using
6000 evenly-spaced values of $c$ in the interval $[0,3.5]$.
Each computed value of $\lb_*(c)$ is calculated using
a bisection method: we find suitable $\lb^{\text{\sc top}}$ and
$\lb^{\text{\sc bot}}$ with $\lb^{\text{\sc top}}>\lb_*(c)>\lb^{\text{\sc bot}}$
by checking if equations \eqref{4.ecuclini} are in \hyperlink{CA}{\sc Case~A}
or \hyperlink{CC}{\sc C}, respectively. Then we make the
same test for $\lb=(1/2)(\lb^{\text{\sc top}}+\lb^{\text{\sc bot}})$
to update the values of $\lb^{\text{\sc top}}$ and $\lb^{\text{\sc bot}}$.
We iterate this process until
$\lb^{\text{\sc top}}-\lb^{\text{\sc bot}}\le 10^{-8}$.
As explained in Theorem \ref{4.teortip}, one has $\lb_*(c)=0$ at a
tipping points $c$. Our numerical results suggest that for the
particular choice of $p$ in \eqref{4.eqpt},
the function $\lb_*$ seems to have three zeros,
although we do not have a theoretical justification for
this. 
Figure~\ref{fig:point1}
provides a closer look at the first of the three tipping points
appearing in Figure~\ref{fig:lambda}. The aim is to illustrate the
uniform convergence of $\mac$ and $\mrc$ to the unique bounded
solution at the tipping point, as proved in Theorem \ref{4.teortip}.
\par
The relation $\lb^*(0,p+\lambda)=\lb^*(0,p)-\lambda$ (proved
in Theorem \ref{3.teorlb*}(v)) and the fact that $\lb_*$ is bounded
(proved in Theorem~\ref{4.teorcont}(i)),
show that we can modify the function
$p$ in order to get examples of equations \eqref{4.ecuini}
for which no tipping occurs. For instance, replacing $p$ by $p_1:=p+\lb_*(0)-1$,
the corresponding equations \eqref{4.ecuini}$_c$ are in \hyperlink{CC}{{\sc Case C}}
for any $c\ge 0$: using Theorem \ref{4.teorcont}(i), we observe that
$\lb^*(0,p_1-q_c)=\lb_*(c)-\lb_*(0)+1\ge 1$ for any $c\ge 0$; of course,
this function $p_1$ does not satisfy Hypothesis~\ref{4.hipo}.
For $p_2:=p+\sup_{c\ge 0}\lb_*(c)+1$, we are always
in \hyperlink{CA}{{\sc Case A}}, since
$\lb^*(0,p_2-q_c)=\lb_*(c)-\sup_{c\ge 0}\lb_*(c)-1\le-1$.
Assuming the accuracy of the
representation of $\lb_*$, we can also get functions
$p_3$ and $p_4$ for which the corresponding $\lb_*$ takes the value
0 at a local maximum
or at a local minimum, so that a punctual \hyperlink{CB}{{\sc Case B}}
\lq\lq interrupts\rq\rq~\hyperlink{CA}{{\sc Cases A}} or
\hyperlink{CC}{{\sc C}}.
\par
\begin{figure}
\caption{Global dynamics of \eqref{4.ecupm}. The behavior is analogous at any
interval of integration.}
\includegraphics[trim={2.5cm 0cm 2.3cm -0.3cm},clip,
width=0.9\textwidth]{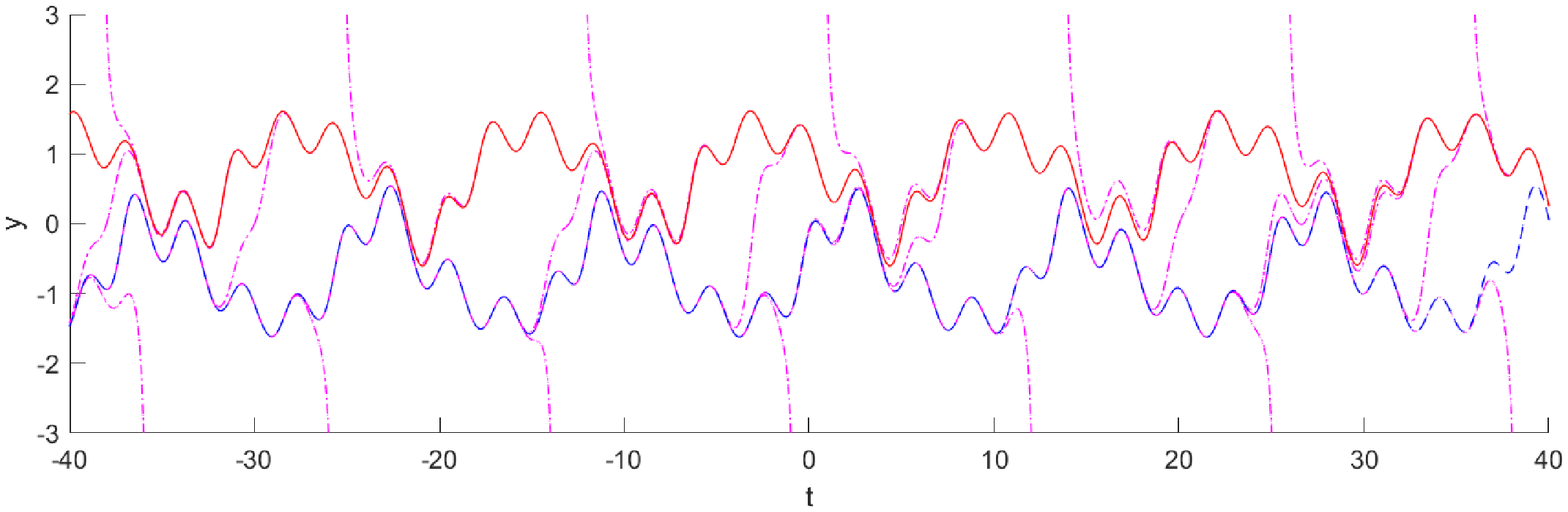}
\label{fig:assump}
\end{figure}

We conclude by explaining the numerical evidences we have mentioned at the beginning
of this subsection. We obtain them by representing solutions of \eqref{4.ecupm}.
Independently of the initial time, the numerical approximation of every
solution starting in an initial value greater than $2$ eventually falls
onto the graph of the function $\wit a_*$, which we represent in solid red in
Figure \ref{fig:assump}.
The analogous behavior is observed backwards in time when approximating
solutions with
initial value less than $-2$, which are eventually mapped on the graph
of $\wit r_*$, represented in dashed blue in Figure \ref{fig:assump}.
In addition, the solution corresponding to
any initial pair (initial time, initial value) between the graphs of $\wit a_*$ and
$\wit r_*$ falls onto the red curve as time increases and onto the blue curve as
time decreases. In other words, we observe numerically that
$(\wit a_*,\wit r_*)$ is an attractor-repeller pair for \eqref{4.ecupm},
which is our first assumption.
Finally, the time of collision of solutions starting at a distance less than 5
to one of the hyperbolic solutions is never greater than 20, which is much less than
$950$ and justifies the validity of our second assumption \eqref{4.cond2}.
%%%%%%%%%%%%%%%%%%%%%%%%%%%%%%%%%%%%%%%%%%%%%%%%%%%%%%%%%%%%%%%%%%%%%%%%%%%%%%%%%%%%%
%%%%%%%%%%%%%%%%%%%%%%%%%%%%%%%%%%%%%%%%%%%%%%%%%%%%%%%%%%%%%%%%%%%%%%%%%%%%%%%%%%%%%
%%%%%%%%%%%%%%%%%%%%%%%%%%%%%%%%%%%%%%%%%%%%%%%%%%%%%%%%%%%%%%%%%%%%%%%%%%%%%%%%%%%%%
%%%%%%%%%%%%%%%%%%%%%%%%%%%%%%%%%%%%%%%%%%%%%%%%%%%%%%%%%%%%%%%%%%%%%%%%%%%%%%%%%%%%%
\section*{Conclusions and outlook}
The first part of this work gives a complete mathematical description of the possible dynamical scenarios for  concave scalar differential equations of the form $x'= -x^2+q(t)x+p(t)$, where $q$ and $p$ are bounded and uniformly continuous functions. The obtained results are used to  study  rate-induced tipping for differential equations that are asymptotically non-autonomous. Specifically, they allow to characterize the occurrence of a rate-induced tipping in terms of  a non-autonomous saddle-node bifurcation for the class of differential equations  
\begin{equation}\label{eq:28/10-17:01}
 y' =-\Big(y-\frac{2}{\pi}\,\arctan(ct)\Big)^2+p(t)\,,
\end{equation}
where $p:\mathbb R \to \mathbb R$ is a bounded and uniformly continuous function, and $c\ge0$. The numerical analysis carried out in this paper shows that an intrinsically non-autonomous system subjected to a time-dependent variation of a parameter can admit more than one tipping point upon the variation of $c>0$. \par\smallskip

This work shades light on an important connection between rate-induced tipping and non-autonomous bifurcation theory. Moreover, it aims to  prompt further studies on this subject which, in spite of its inherent complexity, appears to be the natural framework for rate-induced tipping. Natural followups include for example the following points.
\begin{itemize}[leftmargin=*,itemsep=1pt]
\item In regard to \eqref{eq:28/10-17:01}, it seems natural to ask if $\arctan(ct)$ can be replaced by any other odd asymptotically constant function $\Gamma(ct)$. One can show that if $\Gamma$ is also $C^1$ then the change of variable \eqref{4.changeofvar} is still possible. Appropriate extensions of the results are, however, far from straightforward. In fact, an in-depth analysis shows that some properties only require that $\Gamma$ is uniformly continuous and bounded (such as the continuity of $\lambda^*$ provided by Theorem \ref{4.teorcont} for example). Other properties hold if $\Gamma$ is increasing, and some of our results seem to depend also on the properties of convexity and concavity of $\arctan$ (such as the number of tipping points for example).
\item The study of analogous results for higher-order (instead of quadratic) polynomials  is very interesting, but requires further significant analysis. Although the dynamical scenarios presented in this paper are still possible for some polynomials of order three or higher, different types of bifurcations or dynamical behaviors can emerge. We note that for our work, it is important that the model was globally concave and coercitive. For instance, some properties could remain valid for third-order polynomials, but important restrictions may need to be added to the equation. Moreover, it can be expected that the results will not be global in this case.   
\item Alkhayuon and Ashwin \cite{alas} introduce the concept of partial and total tipping in the context of asymptotically autonomous systems whose attractors are compact sets, each  given by the trajectory of an orbitally asymptotically stable solution. The current state of the art does not allow to pose the same question in the context of asymptotically nonautonomous systems. Nevertheless, a different interpretation of partial and total tipping seems plausible.  If we consider our past limit system $y’=-(y+1)^2 +p(t)$ with attractor $a(t)$ and define the time-translation of $a$ at time $s \in \R$ as $a_s(t)=a(t+s)$, then the translated solution $a_s$ is not a solution of the system anymore (translation invariance only holds for autonomous systems in general), but it is a solution of the translated system $y’=-(y+1)^2 +p_s(t)$, with $p_s(t)=p(t+s)$. It seems natural that the existence and number of  tipping points for the family of past limit equations obtained for $\{p_s\mid s\in\R\}$ would depend on $s$. In some cases and for some values of $c$, the existence of separate hyperbolic solutions persists for all values of $s$, but in other cases, this feature persists only for some values of $s$ and may fail for others, and finally it is also possible that a critical transition occurs for all values of $s$. 
\end{itemize}

\appendix
\section{Proof of Theorem \ref{3.teorhyp}}\label{a.appendix}
Let $q$ and $p$ be the bounded and uniformly continuous functions of
differential equation \eqref{3.ecucon}, and
let $\W_{q,p}$ be its {\em hull\/}; that is, the closure
in $C(\R,\R\times\R)$ of the set $\{(q,p)_t\,|\;t\in\R\}$, where
$C(\R,\R\times\R)$ is endowed
with the compact-open topology and $(q,p)_t(s):=(q(t+s),p(t+s))$.
It is well-known that $\W_{q,p}$ is a compact metric space, and that the map
$\sigma\colon\R\times\W_{q,p}\to\W_{q,p},\,(t,\w)\mapsto\w_t$, with $\w_t(s):=\w(t+s)$,
defines a continuous flow on it (see, e.g., \cite{sell2}). And it is obvious that
the operator $(q_*,p_*)\colon\W_{q,p}\to\R\times\R$, $\w\mapsto\w(0)$, is continuous
and satisfies $(q_*,p_*)(\w_t)=\w(t)$. Note that if $\w=(\w_1,\w_2)$ then
$q_*(\w)=\w_1(0)$ and $p_*(\w)=\w_2(0)$. We will
represent $\w^{q,p}:=(q,p)\in\W_{q,p}$, so that $q_*((\w^{q,p})_t)=q(t)$ and
$p_*((\w^{q,p})_t)=p(t)$. Now we can consider the family of scalar equations
\begin{equation}\label{A.ecuw}
 x'=-x^2+q_*(\w_t)\,x+p_*(\w_t)
\end{equation}
for $\w\in\W_{q,p}$, which includes \eqref{3.ecucon}. Let
$t\mapsto u(t,\w,x_0)$ represent
the maximal solution of \eqref{A.ecuw} with $u(0,\w,x_0)=x_0\in\R$.
The continuity of the flow on the hull, the uniqueness of solutions of
initial value problems for \eqref{A.ecuw}, and standard results on
continuous dependence for ordinary differential equations, ensure that
\begin{equation}\label{A.deftau}
 \tau\colon\mU\subseteq\R\times\W_{q,p}\times\R\to\W_{q,p}\times\R\,,\quad
 (t,\w,x_0)\mapsto (\w_t,u(t,\w,x_0))\,,
\end{equation}
defines a (local) continuous flow on $\W_{q,p}\times\R$. The set $\mU$
is obviously composed by those points $(t,\w,x_0)$ for which
$u(t,\w,x_0)$ exists. Clearly, the (scalar) flow $\tau$ is
monotone with respect to its state variable; i.e.,
if $x_1<x_2$, then $u(t,\w,x_1)<u(t,\w,x_2)$ as long as both
solutions are defined. In addition, the flow $\tau$ is $C^1$ and strictly concave
with respect to the state variable.
The strict concavity means that
\begin{equation}\label{A.conc1}
 u(t,\w,\rho\,x_1+(1-\rho)\,x_2)>\rho\,u(t,\w,x_1)+(1-\rho)\,u(t,\w,x_2)
\end{equation}
whenever $t>0$, $\w\in\W$, $x_1,x_2\in\R$ and $\rho\in(0,1)$,
and as long as all the involved terms are defined.
This property can be proved as for \eqref{3.concx}.
Taking $t<0$ reverts the sign of the inequalities.
Note also that
\begin{equation}\label{A.eqsol}
 x(t,0,x_0)=u(t,\w^{q,p},x_0) \qquad\text{and}\qquad
 x(t,s,x_0)=u(t-s,(\w^{q,p})_s,x_0)
\end{equation}
whenever the right (or left) terms of these equalities are defined.
\par
\begin{nota}\label{A.notarec}
The pair of functions $(q,p)$ (or the differential equation \eqref{3.ecucon})
is said to be {\em recurrent\/} if the flow
on its hull $\W_{q,p}$ is minimal. This is for instance the
case if $q\colon\R\to\R$ and $p\colon\R\to\R$ are almost periodic functions
(as deduced from the results in \cite[Chapter~1]{fink},
combined with \cite[Proposition~IV.2.3]{vrie}).
\end{nota}
Recall that Theorem~\ref{3.teorhyp}
refers to the solutions $a$ and $r$ of equation \eqref{3.ecucon} provided
by Theorem~\ref{3.teoruno}, which are assumed to be globally defined
and uniformly separated.\hspace{-1cm}~
\smallskip\par
\noindent{\it Proof of Theorem \ref{3.teorhyp}}.~
(i) The proof of the assertions in (i) is made in six steps,
making use of the flow $\tau$ defined by \eqref{A.deftau}.
We define
\[
 \delta:=\inf_{t\in\R}(a(t)-r(t))>0\,.
\]
\par
\hypertarget{S1}{{\sc Step 1}}. Let us consider the family of equations \eqref{A.ecuw}.
We prove that there exist globally defined and bounded
functions $a_*\colon\W_{q,p}\to\R$ and $r_*\colon\W_{q,p}\to\R$
such that $u(t,\w,x_0)$ is globally defined and bounded
if and only if $r_*(\w)\le x_0\le a_*(\w)$,  with
$\inf_{\w\in\W_{q,p}}(a_*(\w)-r_*(\w))\ge\delta$, and such that
\begin{equation}\label{A.solar}
 \!\!\!\!u(t,\w,a_*(\w))=a_*(\w_t) \;\text{ and }\;
 u(t,\w,r_*(\w))=r_*(\w_t)\;\text{ for $\w\in\W_{q,p}$ and $t\in\R$}\,.
\end{equation}
\par
Let us define
\[
 \mB_\w:=\{x_0\in\R\,|\;u(t,\w,x_0)\;\text{is globally defined and bounded}\}\,.
\]
We fix $m>0$ such that $-m^2+q_*(\w)\,m+p_*(\w)<-1$ for any $\w\in\W_{q,p}$.
We check that $\mB_\w$ is nonempty and contained in $[-m,m]$
for any $\w\in\W_{q,p}$.
Let us take any $\w^0\in\W_{q,p}$, choose $(t_n)$ such that $\w^0=
\lim_{n\to\infty}(\w^{q,p})_{t_n}$, and
assume without restriction that there exist $a^0=\lim_{n\to\infty}a(t_n)$ and
$r^0=\lim_{n\to\infty}r(t_n)$. Then,
\begin{align*}
 u(t,\w^0,a^0)&=\lim_{n\to\infty} u(t,(\w^{q,p})_{t_n},a(t_n))=
 \lim_{n\to\infty} u(t,u(t_n,\w^{q,p},a(0)))\\
 &=\lim_{n\to\infty}a(t+t_n)\,,
\end{align*}
so that the solution $u(t,\w^0,a^0)$ is bounded on its domain
and hence globally defined. This ensures that
$a^0\in\mB_{\w^0}$. Similarly, $r^0\in\mB_{\w^0}$. Now we define
$a_*(\w^0):=\sup\mB_{\w^0}$ and $r_*(\w^0):=\inf\mB_{\w^0}$.
Theorem~\ref{3.teoruno}(i)\&(ii) ensure that $a_*(\w^0)$ and $r_*(\w^0)$
belong to $[-m,m]$: they are the functions provided by that theorem
for the equation $x'=-x^2+p_*((\w^0)_t)$ evaluated at $t=0$,
and the corresponding condition \eqref{3.defm} is satisfied.
In addition,
$a_*(\w^0)-r_*(\w^0)\ge a^0-r^0=\lim_{n\to\infty}(a(t_n)-r(t_n))\ge\delta$.
To complete \hyperlink{S1}{{\sc Step 1}},
note that
\begin{equation}\label{A.tap}
 a_*((\w^{q,p})_t)=a(t)\quad\text{and}\quad
 r_*((\w^{q,p})_t)=r(t)\quad\text{for all $t\in\R$}\,,
\end{equation}
that the analogous equalities hold for any $\w\in\W$, and that they
guarantee \eqref{A.solar}.
\smallskip\par
\hypertarget{S2}{{\sc Step 2}}. We prove that, if $\mM\subset\W_{q,p}$ is a minimal
set, then the maps $\mM\to\R\,,\;\w\mapsto a_*(\w)$, and
$\mM\to\R\,,\;\w\mapsto r_*(\w)$, are continuous; and that
given $\rho>0$, there exist
$\beta_{\mM}>0$ and $k_{\rho,\mM}\ge 1$ such that
\begin{equation}\label{A.ar*}
\begin{split}
 &|a_*(\w_t)-u(t,\w,x_0)|\le k_{\rho,\mM}e^{-\beta_{\mM} t}|a_*(\w)-x_0|\\
 &\qquad\qquad\qquad
 \quad\text{for $t\ge 0$, $\w\in\mM$ and $x_0\ge r_*(\w)+\rho$}\,,\\
 &|r_*(\w_t)-u(t,\w,x_0)|\le k_{\rho,\mM}e^{\beta_{\mM} t}|r_*(\w)-x_0|\\
 &\qquad\qquad\qquad
 \quad\text{for $t\le 0$, $\w\in\mM$ and $x_0\le a_*(\w)-\rho$}\,.
\end{split}
\end{equation}
\par
For the map $a_*$, the assertions follow from the strict concavity of
the flow $\tau$ defined by \eqref{A.deftau}, and
\cite[Theorems 3.12 and 3.8(iv)]{nuos4}. To prove the assertions
concerning $r_*$, we proceed as follows:
\begin{itemize}
\item[-] We first define a new flow on the hull $\W_{q,p}$ by reversion of time:
$\sigma^-(t,\w)=\w_{-t}$.
\item[-] We consider the family of equations $z'=-z^2-q_*(\w_{-t})\,z+
p_*(\w_{-t})$, which induces the new flow
\[
 \qquad\tau^-\colon\mU\subseteq\R\times\W_{q,p}\times\R\to\W_{q,p}\times\R\,,\quad
 (t,\w,z_0)\mapsto (\w_{-t},w(t,\w,z_0))\,,
\]
and check that $w(t,\w,z_0)=-u(-t,\w,-z_0)$. Consequently,
\[
 \mB_\w^-:=\{z_0\,|\;{-}a_*(\w)\le z_0\le -r_*(\w)\}
\]
are the sets of initial conditions giving rise to bounded solutions.
\item[-] We observe that the differential equation corresponding to $\w^{q,p}$ is given by
$z'=-z^2-q(-t)\,z+p(-t)$, and hence that
the roles of the solutions $\bar a(t)=-a(-t)$ and $\bar r(t)=-r(-t)$ of
this differential equation correspond to the roles played by $r$ and $a$
for the original differential equation, respectively; and we have
$\inf_{t\in\R}(\bar r(t)-\bar a(t))=\delta>0$. Therefore, the results obtained
so far for $a$ and $a_*$ apply to $\bar r$ and $-r_*$, leading us to conclude that
given $\rho>0$, there exist $\beta_{\mM}>0$ and $k_{\rho,\mM}\ge 1$ such that
\[
\begin{split}
 &\qquad|{-}r_*(\w_{-t})-w(t,\w,z_0)|\le k_{\rho,\mM}e^{-\beta_{\mM} t}|{-}r_*(\w)-z_0|\\
 &\qquad\qquad\qquad\qquad
 \quad\text{for $t\ge 0$, $\w\in\mM$ and $z_0\ge -a_*(\w)+\rho$}\,,
\end{split}
\]
This fact and the previous equalities prove the second assertion in \eqref{A.ar*},
which completes this step.
\end{itemize}
\smallskip\par
\hypertarget{S3}{{\sc Step 3}}. We prove that any point $\w^1$ belonging
to a minimal subset of $\W_{q,p}$
is a continuity point for $a_*\colon \W_{q,p}\to\R$ and $r_*\colon \W_{q,p}\to\R$.
\par
We take a sequence $(\w_n)$ in $\W_{q,p}$ with limit $\w^1$.
We check that $a_*(\w^1)=\lim_{n\to\infty}a_*(\w_n)$ and
$r_*(\w^1)=\lim_{n\to\infty}r_*(\w_n)$.
Note that it is enough to check the following: if for a subsequence $(\w_m)$,
there exist
$a^1:=\lim_{m\to\infty} a_*(\w_m)$ and $r^1:=\lim_{m\to\infty}r_*(\w_m)$,
then $a^1=a_*(\w^1)$ and $r^1=r_*(\w^1)$.
\par
So that let us take such a subsequence. It follows from
\eqref{A.solar} that $u(s,\w^1\!,a^1)=\lim_{m\to\infty} a_*((\w_m)_s)$ and
$u(s,\w^1\!,r^1)=\lim_{m\to\infty} r_*((\w_m)_s)$ for $s\in\R$,
which, as seen in \hyperlink{S1}{{\sc Step 1}}, ensures that
$u(s,\w^1\!,a^1)$ and $u(s,\w^1\!,r^1)$ are bounded solutions.
Hence, $r_*(\w^1)\le r^1\le a^1\le a_*(\w^1)$ and, for any $s\in\R$,
\[
\begin{split}
 u(s,\w^1\!,a^1)-r_*((\w^1)_s)
 & =u(s,\w^1\!,a^1)-u(s,\w^1\!,r_*(\w^1))\ge
 u(s,\w^1\!,a^1)-u(s,\w^1\!,r^1)\\
 & =\lim_{m\to\infty}\big(a_*((\w_m)_s)-r_*((\w_m)_s)\big)\ge \delta\,.
\end{split}
\]
According to \hyperlink{S2}{{\sc Step 2}}, there exist
$\beta_{\mM}>0$ and $k_{\delta,\mM}\ge 1$
such that
\[
 |a_*((\w^1)_t)-u(t-s,(\w^1)_s,u(s,\w^1\!,a^1))|
 \le k_{\delta,\mM}\,e^{-\beta_{\mM}(t-s)}
 |a_*((\w^1)_s)-u(s,\w^1\!,a^1)|
\]
if $t\ge s$. Therefore, taking $t=0$,
\[
 |a_*(\w^1)-a^1| \le k_{\delta,\mM}\,e^{\beta_{\mM} s}
 |a_*((\w^1)_s)-u(s,\w^1\!,a^1)|
\]
for $s\le 0$. Since $\sup_{s\in\R}|a_*((\w^1)_s)-u(s,\w^1\!,a^1)|<\infty$,
we conclude that $a_*(\w^1)=a^1$, as asserted. The proof of
$r^1=r_*(\w^1)$ can be done in an analogous way.
\smallskip\par
\hypertarget{S4}{{\sc Step 4}}. We prove that the maps $a_*$ and $r_*$ are
continuous on $\W_{q,p}$.
\par
We take $\w^2\in\W_{q,p}$ and a sequence
$(\w_n)$ in $\W_{q,p}$ with limit $\w^2$.
We will check that $a_*(\w^2)=\lim_{n\to\infty}a_*(\w_n)$ and
$r_*(\w^2)=\lim_{n\to\infty}r_*(\w_n)$. Again,
we will check that if for a subsequence $(\w_m)$, there exist
$a^2:=\lim_{m\to\infty} a_*(\w_m)$ and $r^2:=\lim_{m\to\infty}r_*(\w_m)$,
then $a^2=a_*(\w^2)$ and $r^2=r_*(\w^2)$.
\par
The hypothesis on uniform separation ensures that
$r^2\le a^2-\delta$. As seen in \hyperlink{S2}{{\sc Step~2}}, the solutions
$u(s,\w^2\!,a^2)$ and $u(s,\w^2\!,r^2)$ are bounded,
and hence $r_*(\w^2)\le r^2<a^2\le a_*(\w^2)$.
This and \eqref{A.solar} ensure that, for any $s\in\R$,
\[
\begin{split}
 a_*((\w^2)_s)-u(s,\w^2\!,r^2)
 & =u(s,\w^2\!,a_*(\w^2))-u(s,\w^2\!,r^2)\ge
 u(s,\w^2\!,a^2)-u(s,\w^2\!,r^2)\\
 & =\lim_{m\to\infty} a_*((\w_m)_s)-r_*((\w_m)_s)\ge \delta\,.
\end{split}
\]
Now let us
assume for contradiction that $r^2>r_*(\w^2)$, and look for $\rho\in(0,1)$
such that $r^2=\rho\,a_*(\w^2)+(1-\rho)\,r_*(\w^2)$. Then, the
concavity of the flow $\tau$ ensures that
$u(s,\w^2,r^2)\ge\rho\,a_*((\w^2)_s)+(1-\rho)\,r_*((\w^2)_s)$
for any $s\ge 0$ (see \eqref{A.conc1}),
and hence,
\[
 u(s,\w^2,r^2)-r_*((\w^2)_s)\ge \rho\,\big(a_*((\w^2)_s)
 -r_*((\w^2)_s)\big)\ge\rho\,\delta
\]
for any $s\ge 0$.
Let us take a point $(\w^1,x^1)$ in a minimal set contained in the
omega-limit set of $(\w^2,r^2)$ for the flow $\tau$, and note
that $\w^1$ belongs to a minimal subset $\mM$ of $\W_{q,p}$ for the flow on the hull.
It follows easily from the previous inequalities, the continuity of $\tau$,
and the continuity of $r_*$ and $a_*$ at $\w^1\in\mM$ established on
\hyperlink{S4}{{\sc Step 4}} that
\[
 r_*((\w^1)_t)+\rho\,\delta\le u(t,\w^1,x^1)\le a_*((\w^1)_t)-\delta
\]
for any $t\ge 0$. But this contradicts the information regarding
the asymptotic behavior of the solutions provided by \hyperlink{S2}{{\sc Step 2}}:
$\lim_{t\to\infty}(a_*((\w^1)_t)-u(t,\w^1,x^1))=0$,
since $x^1\ge r_*(\w^1)+\rho\,\delta$.
This contradiction shows that $r^2=r_*(\w^2)$.
\par
The proof for $a_*$ can
be done similarly, now working for negative values of time and for
the alpha-limit set of $(\w^2,a^2)$. This completes \hyperlink{S5}{{\sc Step 4}}.
\smallskip\par
\hypertarget{S5}{{\sc Step 5}}. We prove the absence of non-trivial
bounded solutions for all the linear differential equations
\begin{equation}\label{A.families}
 y'= (2\,a_*(\w_t)-q_*(\w_t))\,y\quad \text{and}\quad
 y'=(2\,r_*(\w_t)-q_*(\w_t))\,y\,.
\end{equation}
Let us work with $a_*$, assuming for contradiction the
existence of $\w^4$ such that
\[
 \sup_{t\in\R}\left(\exp\int_{0}^{t} (2\,a_*((\w^4)_l)-q_*((\w^4)_l))\,\rmd l\right)=:
 \kappa<\infty\,,
\]
and let us take any $x^4\in(r_*(\w^4),a_*(\w^4))$.
%According to \eqref{3.conc2}
By repeating the argument leading to \eqref{3.tmays}, we get
\[
\begin{split}
 a_*((\w^4)_t)-u(t,\w^4,x^4)%\ge u_x(t,\w^4,a_*(\w^4))(a_*(\w^4)-x^4)
 &\ge(a_*(\w^4)-x^4)\exp\int_0^t(-2\,a_*((\w^4)_l)+q_*((\w^4)_l))\,\rmd l\\
 &\ge (1/\kappa)(a_*(\w^4)-x^4)=:\delta_3>0
\end{split}
\]
for any $t>0$. By reasoning as in \hyperlink{S4}{{\sc Step 4}}, we
prove that there exists $\rho\in(0,1)$ with
\[
  u(t,\w^4,x^4)-r_*((\w^4)_t)\ge \rho\,\delta
\]
for any $t\ge 0$, and we reach the required contradiction by repeating the
final argument of \hyperlink{S4}{{\sc Step 4}}.
\smallskip\par
\hypertarget{S6}{{\sc Step 6}}.
We prove that the solutions $a_*(\w_t)$ and $r_*(\w_t)$
of equation \eqref{A.ecuw}$_\w$ (see \eqref{A.solar}) are hyperbolic,
which combined with \eqref{A.tap} proves (i).
\par
We first take a point $\w$ in a minimal subset $\mM$
of $\W_{q,p}$, and $x_0\in[a_*(\w),r_*(\w)]$.
By repeating the argument leading to \eqref{3.tmays} and
\eqref{3.tmens}, and
using the information from \hyperlink{S2}{{\sc Step 2}},
we find
\[
\begin{split}
 &\exp\int_0^t(-2\,a_*(\w_l)+q_*(\w_l))\,\rmd l\le
 \frac{a_*(\w_t)-u(t,\w,x_0)}{a_*(\w)-x_0}\le
 k_{\rho,\mM}e^{-\beta_{\mM} t}\quad \text{if $t\ge 0$}\,,\\
 &\exp\int_0^t(-2\,r_*(\w_l)+q_*(\w_l))\,\rmd l\le
 \frac{u(t,\w,x_0)-r_*(\w_t)}{x_0-r_*(\w)}\le
 k_{\rho,\mM}e^{\beta_{\mM} t}\quad \text{if $t\le 0$}\,.
\end{split}
\]
It follows easily that the equations
\begin{equation}\label{A.fam}
 y'=(-2\,a_*(\w_t)+q_*(\w_t))\,y\quad \text{and}\quad
 y'=(-2\,r_*(\w_t)+q_*(\w_t))\,y
\end{equation}
have an exponential dichotomy for each element $\w$ of any minimal subset $\mM$
of $\W_{q,p}$, where the first ones are of Hurwitz type at $+\infty$ and
the second ones are of Hurwitz type at $-\infty$.
Therefore, also the equations \eqref{A.families}
have an exponential dichotomy for each element $\w$ of each
minimal subset $\mM$ of $\W_{q,p}$, where the first ones are
of Hurwitz type at $-\infty$ and the second ones are of Hurwitz
type at $+\infty$,
see e.g.~\cite[Proposition 1.73 and Theorem 1.60]{jonnf}.
This fact combined with the absence of bounded solutions established
in \hyperlink{S5}{{\sc Step 5}}
guarantees the exponential dichotomy of the equations
\eqref{A.families} for any $\w\in\W_{q,p}$, see \cite[Theorem 2]{sase3},
and applying
\cite[Proposition 1.73 and Theorem 1.60]{jonnf} again ensures that
every differential equation in \eqref{A.fam}
has an exponential dichotomy, which proves our assertion.
\smallskip\par
(ii) Let us concentrate on the hyperbolic solution $a$ first.
We fix a dichotomy constant pair for $a$ (see Section~\ref{2.sec}),
and now check that
$(k_a,\beta_a)$ is a dichotomy constant pair for the hyperbolic
solution $a_*(\w_t)$ of equation \eqref{A.ecuw}$_\w$
(see \eqref{A.solar}), for any $\w\in\W_{q,p}$.
So, we take $\w\in\W_{q,p}$ and obtain it as limit
$\w=\lim_{n\to\infty}(\w^{q,p})_{t_n}$
for a suitable sequence $(t_n)$. The continuity of the flow on
$\W_{q,p}$, the continuity of $a_*$ (proved in
\hyperlink{S4}{{\sc Step 4}} of (i)) and the first equality in
\eqref{A.tap} show that
\[
\begin{split}
 &\exp\int_s^t (-2\,a_*(\w_l)+q_*(\w_l))\,\rmd l
 =\lim_{n\to\infty}\exp\int_s^t (-2\,a(t_n+l)+q(t_n+l))\,\rmd l\\
 &\qquad\qquad=\lim_{n\to\infty}\exp\int_{s+t_n}^{t+t_n} (-2\,a(l)+q(l))\,\rmd l
 \le k_a\,e^{-\beta_a(t-s)}\quad\text{if $t\ge s$}\,,
\end{split}
\]
which proves the assertion.
\par
Let us fix $\bar\beta_a\in(0,\beta_a)$.
Now we reason as in the proof of Proposition \ref{3.prophyp}.
By reviewing the proof of \cite[Theorem III.2.4]{hale},
we conclude that there exists $\rho>0$ (which depends just
on the choice of $\bar\beta_a$) such that
\begin{equation}\label{A.asi}
\begin{split}
 &|a_*(\w_t)-u(t,\w,x_0)|\le k_ae^{-\bar\beta_a t}|a_*(\w)-x_0|\\
 &\qquad\qquad\qquad
 \quad\text{for $\,t\ge 0\,$,  $\,\w\in\W_{q,p}\,$ and
 $\,|x_0-a_*(\w)|<\rho$}\,.
\end{split}
\end{equation}
And we also check (as in Proposition \ref{3.prophyp}) that
\begin{equation}\label{A.asin}
\begin{split}
 &|a_*(\w_t)-u(t,\w,x_0)|\le k_a\,e^{-\beta_a\, t}\,|a_*(\w)-x_0|\\
 &\qquad\qquad\qquad
 \quad\text{for $\,t\ge 0\,$,  $\,\w\in\W_{q,p}\,$ and
 $\,x_0\ge a_*(\w)$}\,.
\end{split}
\end{equation}
Now we fix $\ep>0$ and check that for each $\w^0\in\W_{q,p}$,
there exists a time $t_{\w^{\!0}}$ such that
\[
 |u(t_{\w^{\!0}},\w^0,r_*(\w^0)+\ep)-a_*((\w^0)_{t_{\w^{\!0}}})|< \rho\,.
\]
If $r(\w^0)+\ep\ge a_*(\w^0)-\rho$, this inequality follows from \eqref{A.asi}
and \eqref{A.asin}.
Assume hence that $r_*(\w^0)+\ep\le a_*(\w^0)-\rho$.
The non-existence of $t_{\w^0}$ would ensure that
$a_*((\w^0)_t)-u(t,\w^0,r_*(\w^0)+\ep)\ge\rho$ for
every $t\ge 0$. In that case, we take $(\w^1,x^1)$ in
a minimal set contained in the omega-limit set of $(\w^0,r_*(\w^0)+\ep)$
and conclude from the continuity of $\tau$
and that of $a_*$ (proved in (i)) that
$a_*((\w^1)_t)-u(t,\w^1,r_*(\w^1)+\ep)\ge\rho$ for
every $t\ge 0$. But this contradicts the
first inequality in \eqref{A.ar*}, so that our assertion is proved.
We point out that $t_{\w^{\!0}}$ depends on $\w^0$,
$\ep$ and $\rho$, and hence on $\w^0$, $\ep$ and $\bar\beta_a$.
\par
Note now that for any $\w^0\in\W_{q,p}$
there exists an open neighborhood $\mU_{\w^0}$
such that $|u(t_{\w^{\!0}},\w,r_*(\w)+\ep)-a_*(\w_{t_{\w^{\!0}}})|<\rho$
for any $\w\in\mU_{\w^0}$.
Therefore the compactness of $\W_{q,p}$ provides a finite
number of times $t_1,\ldots,t_n$ such that
for any $\w\in\W_{q,p}$ there exists $j=j(\w)\in\{1,\ldots,n\}$
with $|u(t_j,\w,r_*(\w)+\ep)-a_*(\w_{t_j})|<\rho$.
We define $T:=\max(t_1,\ldots,t_n)$ and note that $T$ depends on
the choices of $\ep$ and $\bar\beta_a$.
\par
Let us now fix $\w\in\W_{q,p}$ and $x_0\ge r_*(\w)+\ep$.
If $x_0\ge a_*(\w)-\rho$, then \eqref{A.asi} and \eqref{A.asin}
ensure that
\begin{equation}\label{A.acal}
 |a_*(\w_t)-u(t,\w,x_0)|\le k_ae^{-\bar\beta_a t}|a_*(\w)-x_0|
 \quad\text{if $t\ge 0$}\,.
\end{equation}
So that we assume that $x_0\in[r_*(\w)+\ep, a_*(\w)-\rho]$.
We choose $j=j(\w)$ as above, and note that
the monotonicity of the flow ensures that
$|u(t_j,\w,x_0)-a_*(\w_{t_j})|<\rho$. Therefore,
if $t\ge T$ (and hence $t\ge t_j$),
\begin{equation}\label{A.aca2}
\begin{split}
 |a_*(\w_t)-u(t,\w,x_0)|&=
 |a_*((\w_{t_j})_{(t-t_j)})-u(t-t_j,\w_{t_j},u(t_j,\w,x_0))|\\
 &\le k_a\,e^{-\bar\beta_a(t-t_j)}|a_*(\w_{t_j})-u(t_j,\w,x_0)|\\
 &\le k_a\,e^{-\bar\beta_a t}e^{\bar\beta_a T}\rho
 \le k_a\,e^{\bar\beta_a T}e^{-\bar\beta_a t}|a_*(\w)-x_0|\,.
\end{split}
\end{equation}
Now we define
\[
 \kappa:=\sup_{t\in[0,T],\;\w\in\W,\;x_0\in[r_*(\w)+\ep,a_*(\w)-\rho]}
 \frac{|a_*(\w_t)-u(t,\w,x_0)|}{e^{-\bar\beta t}|a_*(\w)-x_0|}\,,
\]
which is finite since the function on the right is a continuous
map on a compact metric space. While being independent of $\w$ and $x_0$,
the quantity $\kappa$ depends on $\ep$ and on $\bar\beta_a$ (as $\rho)$.
Then, if $t\in[0,T]$,
\begin{equation}\label{A.aca3}
\begin{split}
 |a_*(\w_t)-u(t,\w,x_0)|\le \kappa\,e^{-\bar\beta_a t}|a_*(\w)-x_0|\,.
\end{split}
\end{equation}
Summing up, assertions \eqref{A.acal}, \eqref{A.aca2} and \eqref{A.aca3} lead to
\[
\begin{split}
 &|a_*(\w_t)-u(t,\w,x_0)|\le k_{a,\ep}\,e^{-\bar\beta_a t}|a_*(\w)-x_0|\\
 &\qquad\qquad\qquad
 \quad\text{for $t\ge 0$,  $\w\in\W_{q,p}$ and $x_0\ge r_*(\w)+\ep$}\,,
\end{split}
\]
where $k_{a,\ep}:=\max(k_a,k_a\,e^{\bar\beta_a T},\kappa)$.
Note that $k_{a,\ep}$ depends on $\ep$ and $\bar\beta_a$, but neither on $\w$
nor on $x_0$.
\par
The first statement concerning $a$ in (ii) follows from the
previous assertion and
the equalities \eqref{A.tap} and \eqref{A.eqsol}, and the second
one from these equalities and \eqref{A.asin}.
\par
In order to prove the result for $r$, we make the change of variable
$z(t)=-x(-t)$, which transforms the equation $x'=-x^2+q(t)\,x+p(t)$ in
$z'=-z^2-q(-t)\,z+p(-t)$,
with solutions $z(t,s,z_0)=-x(-t,-s,-z_0)$.
We apply the results so far proved to the
transformed equation, for which the functions
$\bar r(t)=-r(-t)$ and $\bar a(t)=-a(-t)$ play
the same roles as $a$ and $r$ for the
initial one, respectively. From here the proof is easily completed:
see the end of the proof of \hyperlink{S2}{{\sc Step 2}} in (i).
\smallskip\par
(iii) Let us check that $a$ and $r$ are uniformly separated,
reasoning by contradiction: it is easy to deduce from (i) that
otherwise there exist $k\ge 1$ and
$s_0$ such that $\exp\int_{s_0}^t(-2\,r(l)+q(l))\,\rmd l\le k\,e^{-(\beta_a/2)(t-s_0)}$
whenever $t\ge s_0$. But this implies that the solutions
$z(t,s_0,z_0)$ of the linearized equation $z'=(-2\,r(t)+q(t))\,z$
tend to $0$ as $t\to\infty$, which is impossible, see the comments after
the definition of hyperbolicity in Section~\ref{2.sec}.
\par
Let $b$ be a bounded solution that is different from $a$ and $r$, which implies $a<b<r$.
Then, (ii) ensures that $\lim_{t\to\infty}|a(t)-b(t)|=0$ and
$\lim_{t\to-\infty}|r(t)-b(t)|=0$. Therefore, $a$ and $r$
are the unique bounded and uniformly separated solutions.
As above, the same
property shows that, for the equation
$z'=(-2\,b(t)+q(t))\,z$, there coexist solutions bounded for $t\to\infty$
and solutions bounded for $t\to -\infty$, which precludes the
hyperbolicity of $b$ (see Section~\ref{2.sec}).
The proof is complete.\hfill{$\Box$}
%%%%%%%%%%%%%%%%%%%%%%%%%%%%%%%%%%%%%%%%%%%%%%%%%%%%%%%%%%%%%%%%%%%%%%%%%%%%%%%%%%%%%
%%%%%%%%%%%%%%%%%%%%%%%%%%%%%%%%%%%%%%%%%%%%%%%%%%%%%%%%%%%%%%%%%%%%%%%%%%%%%%%%%%%%%
%%%%%%%%%%%%%%%%%%%%%%%%%%%%%%%%%%%%%%%%%%%%%%%%%%%%%%%%%%%%%%%%%%%%%%%%%%%%%%%%%%%%%
%%%%%%%%%%%%%%%%%%%%%%%%%%%%%%%%%%%%%%%%%%%%%%%%%%%%%%%%%%%%%%%%%%%%%%%%%%%%%%%%%%%%%

\subsection*{Acknowledgments}
We thank two anonymous referees for their comments and suggestions, which have been  included in the current version of the paper.


\begin{thebibliography}{99}

\bibitem{alas} {\rm H.M.~Alkhayoun, P.~Ashwin},
        Rate-induced tipping from periodic attractors: partial
        tipping and connecting orbits,
        {\em Chaos\/} {\bf 28} (3) (2018), 033608, 11 pp.
\bibitem{aaJQW} {\rm H.M.~Alkhayoun, P.~Ashwin, L.C.~Jackson, C.~Quinn,R.A.~Wood},
        Basin bifurcations, oscillatory instability and rate-induced thresholds for Atlantic meridional overturning circulation in a global oceanic box model,
        {\em  Proc. R. Soc. A} {\bf 475} (2225) (2019), 20190051.
        
\bibitem{alob3} {\rm A.I. Alonso, R.~Obaya},
        The structure of the bounded trajectories set of a scalar convex differential equation,
        {\em Proc. Roy. Soc. Edinburgh} {\bf 133 A} (2003), 237--263.
\bibitem{aloo} {\rm A.I. Alonso, R.~Obaya, R. Ortega},
        Differential equations with limit-periodic forcings,
        {\em Proc. Amer. Math. Soc.} {\bf 131} (3) (2002), 851-–857.
\bibitem{anja} {V.~Anagnostopoulou, T.~J\"{a}ger},
        Nonautonomous saddle-node bifurcations: Random and deterministic forcing,
        {\em J. Differential Equations\/} {\bf 253} (2) (2012), 379--399.
\bibitem{aspw} {P.~Ashwin, C.~Perryman, S.~Wieczorek},
        Parameter shifts for nonautonomous systems in low dimension:
        bifurcation and rate-induced tipping,
        {\em Nonlinearity\/} {\bf 30} (6) (2017), 2185--2210.
\bibitem{ashwintippinclass} {\rm P.~Ashwin, S.~Wieczorek, R.~Vitolo, P.~Cox},
        Tipping points in open systems: bifurcation, noise-induced and
        rate-dependent examples in the climate system,
        {\em Phil. Trans. R. Soc. A} \textbf{370} (2012), 1166--1184.
\bibitem{carigimres} {\rm G. Carigi},
        Rate-induced tipping in nonautonomous dynamical systems with bounded noise,
        MRes Thesis, University of Reading, 2017.
 \bibitem{calr} {\rm A. Carvalho, J. Langa, J. Robinson},
        Attractors for Infinite-dimensional Non-autonomous Dynamical Systems,
        {\em Appl. Math. Sci.} {\bf 182}, Springer, 2013.
\bibitem{copp2} {\rm W.A. Coppel},
        Disconjugacy,
        {\em Lecture Notes in Math.} {\bf 220},
        Springer-Verlag, Berlin, Heidelberg, New York, 1971.
\bibitem{copp1} {\rm W.A. Coppel},
        Dichotomies in Stability Theory,
        {\em Lecture Notes in Math.} {\bf 629},
        Springer-Verlag, Berlin, Heidelberg, New York, 1978.
\bibitem{fink} {\rm A.M. Fink},
        Almost Periodic Differential Equations,
        {\em Lecture Notes in Math.} {\bf 377},
        Springer-Verlag, Berlin, Heidelberg, New York, 1974.
\bibitem{fuhr} {\rm G.~Fuhrmann},
        Non-smooth saddle-node bifurcations III:
        Strange attractors in continuous time,
        {\em J. Differential Equations\/} {\bf 261} (3) (2016), 2109--2140.
\bibitem{hale} {\rm J.K. Hale}
        {\em Ordinary Differential Equations\/},
        Wiley-Interscience, New York, 1969.
\bibitem{hartlthesis} {\rm M. Hartl},
        Non-autonomous Random Dynamical Systems:
        Stochastic Approximation and Rate-Induced Tipping,
        PhD Thesis, Imperial College London, 2019.

\bibitem{hill} {\rm A.v.~Hill},
        Excitation and accommodation in nerve,
        {\em  Proc. R. Soc. B} {\bf 119} (814) (1936), 305--355.
\bibitem{kato} {\rm T.~Kato},
        On the Adiabatic Theorem of Quantum Mechanics,
        {\em J. Phys. Soc. Jpn.} {\bf 5} (6) (1950), 435--439.

\bibitem{klra} {\rm P.~Kloeden, M.~Rasmussen},
        Nonautonomous Dynamical Systems,
        {\em Mathematical Surveys and Monographs},
        Amer. Math. Soc., 2011.
%\bibitem{john12} {\rm R. Johnson},
%        The recurrent Hill's equation,
%        {\em J. Differential Equations\/} {\bf 46} (1982), 165--193.
\bibitem{jonnf} {\rm R.~Johnson, R.~Obaya, S.~Novo, C.~N\'{u}\~{n}ez, R.~Fabbri},
        {\em Nonautonomous Linear Hamiltonian Systems:
        Oscillation, Spectral Theory and Control}
        Developments in Mathematics 36, Springer, 2016
\bibitem{jnot} {\rm \`{A}.~Jorba, C.~N\'{u}\~{n}ez, R.~Obaya, J.C.~Tatjer},
        Old and new results on Strange Nonchaotic Attractors,
        {\em Int. J. Bifurcation Chaos\/}, {\bf 17} (11) (2007), 3895--3928.
\bibitem{mill} {\rm V.M. Million\u{s}\u{c}ikov},
        Proof of the existence of irregular systems of linear differential
        equations with almost periodic coefficients,
        {\em Differ. Uravn.} {\bf 4} (3) (1968), 391--396.
\bibitem{mill2} {\rm V.M. Million\u{s}\u{c}ikov},
        Proof of the existence of irregular systems of linear differential
        equations with quasi periodic coefficients,
        {\em Differ. Uravn.} {\bf 5} (11) (1969), 1979--1983.
%\bibitem{noos2} {\rm S. Novo, R.~Obaya, A.M. Sanz},
%        Almost periodic and almost automorphic dynamics for scalar convex
%        differential equations,
%        {\em Israel J. Math.} {\bf 144} (2004), 157--189.
%\bibitem{noos4} {\rm S.~Novo, R.~Obaya, A.M.~Sanz},
%        Stability and extensibility results for abstract skew-product semiflows,
%        {\em J. Differential Equations} {\bf 235} (2007), 623--646.
%\bibitem{noos5} {\rm S. Novo, R.~Obaya, A.M. Sanz},
%        Exponential stability in non-autonomus delayed equations with applications
%        with neural networks, DEJARLO???
%        {\em Discrete Contin. Dynam. Systems\/} {\bf 18} (2\&3) (2007), 517-–536
\bibitem{nuob6} {\rm C.~N\'{u}\~{n}ez, R.~Obaya},
        A nonautonomus bifurcation theory for deterministic scalar
        differential equations,
        {\em Discrete Contin. Dynam. Systems\/}, {\bf 9} (3\&4) (2008), 701--730.
\bibitem{nuos4} {\rm C.~N\'{u}\~{n}ez, R.~Obaya, A.M. Sanz},
        Minimal sets in monotone and concave skew-product semiflows I:
        a general theory,
        {\it J. Differential Equations\/} {\bf 252} (2012), 5492--5517.
\bibitem{nuos5} {\rm C.~N\'{u}\~{n}ez, R.~Obaya, A.M. Sanz},
        Minimal sets in monotone and concave skew-product semiflows II:
        two-dimensional systems of differential equations,
        {\it J. Differential Equations\/} {\bf 252} (2012), 3575--3607.
\bibitem{potz} {\rm C.~P\"{o}tzsche},
        Nonautonomous continuation of bounded solutions,
        {\em Commun. Pure Appl. Anal.} {\bf 10} (3) (2011), 937--961.
\bibitem{rasmln} {\rm M.~Rasmussen},
        Attractivity and bifurcation for nonautonomous dynamical systems,
        {\em Lecture Notes in Math.} {\bf 1907}, Springer-Verlag, Berlin, 2007.
\bibitem{sase3} {\rm R.J.~Sacker, G.R.~Sell},
        Existence of dichotomies and invariant splittings for linear differential systems II,
        {\em J. Differential Equations\/} {\bf 22} (1976), 478--496.
\bibitem{scheffer} {\rm M.~Scheffer},
        {\em Critical Transitions in Nature and Society},
        Princeton University Press, 2009.
\bibitem{svNhh} {\rm M.~Scheffer, E.H.~Van Nes, M.~Holmgren, T.~Hughes},
        Pulse-driven loss of top-down control: the critical-rate hypothesis,
        {\em  Ecosystems } {\bf 11} (2008), 226--237.
\bibitem{sell2} {\rm G.R. Sell},
        {\em Topological Dynamics and Differential Equations\/},
        Van Nostrand Reinhold, London, 1971.
\bibitem{vwf} {\rm A.~Vanselow, S.~Wieczorek, U.~Feudel},
        When very slow is too fast-collapse of a predator-prey system,
        {\em J. Theor. Biol.} {\bf 479} (2019), 64--72.
\bibitem{vino} {\rm R.E. Vinograd},
        A problem suggested by N.P. Erugin,
        {\em Differ. Uravn.} {\bf 11} (4) (1975), 632--638.
\bibitem{vrie} {\rm J.~de Vries},
        {\em Elements of Topological Dynamics},
        Kluwer Academic Publishers, Dordrecht, 1993.

\bibitem{walc} {\rm S.~Wieczorek, P.~Ashwin, C.M.~Luke, P.M.~Cox},
        Excitability in ramped systems: the compost-bomb instability, {\em Proc. R. Soc. A} {\bf 467} (2011), 1243--1269.    
        
\bibitem{wieczcompact} {\rm S.~Wieczorek, C.~Xie, C.K.R.T.~Jones},
        Compactification for asymptotically autonomous dynamical systems:
        theory, applications and invariant manifolds, Preprint (2020).
\bibitem{xiethesis} {\rm C.~Xie},
        Rate-Induced Critical Transitions,
        PhD Thesis, University College Cork, 2020.
\end{thebibliography}
\end{document}